\documentclass[11pt]{article}

\pdfoutput=1

\RequirePackage{my_packages}
\RequirePackage{my_theorems}
\RequirePackage{my_macros}
\RequirePackage{my_tikz}

\usepackage[backend = biber,        % biber, bibtex
            language = english ,
            style = alphabetic-verb ,   % numeric-comp , alphabetic , authoryear , authoryear-comp , authoryear-icomp, alphabetic-verb
            giveninits = true,
            isbn = false,
            url = false,
            doi = false,
            sorting = nyt,
            maxnames = 6,
            backref=true
            ]{biblatex}

\renewbibmacro{in:}{%
  \ifentrytype{article}{}{%
  \printtext{\bibstring{in}\intitlepunct}}}
\AtEveryBibitem{%
  \clearfield{note}%
  \clearfield{edition}%
  \clearlist{language}
}

\DeclareNameAlias{labelname}{given-family}

\title{A list of applications of Stallings automata}
\author[1]{Jordi Delgado\thanks{\url{jdelgado@crm.cat}}}
\author[2]{Enric Ventura\thanks{\url{enric.ventura@upc.edu}}}

\affil[1]{Department of Mathematics, University of the Basque Country (EHU)}
\affil[2]{Departament de Matem\`atiques, Universitat Polit\`ecnica de Catalunya, and Institut de Matem\`atiques de la UPC-BarcelonaTech, Catalonia}

\begin{document}

\maketitle
\begin{abstract}
\noindent
This survey is intended to be a fast (and reasonably updated) reference for the theory of Stallings automata and its applications to the study of subgroups of the free group, with the main accent on algorithmic aspects. Consequently, results concerning finitely generated subgroups have greater prominence in the paper. However, when possible, we try to state the results with more generality, including the usually overlooked non-(finitely-generated) case. 
\end{abstract}

%\vspace{30pt}

\vfill
\hfill
\textsf{Version 1.01}
\newpage
\tableofcontents

\newpage

\section{\bf Introduction}
The classical \emph{theory of Stallings automata}
provides a neat geometric representation of subgroups of the free group, and constitutes the modern --- and probably the most natural and fruitful --- approach to their study. Moreover, if the involved subgroups are finitely generated, then this description is finitary and very well suited for algorithmic treatment.

The original result (hinted by the work of Serre in \parencite{serre_arbres_1977}, and precisely stated in the seminal paper~\parencite{stallings_topology_1983} by Stallings) 
interprets the subgroups of the free group $\Free[A]$ (with basis $A$ of cardinal $n$) as covering spaces of the bouquet of $n$ circles.
Despite of this (mainly topological) original viewpoint, when restricted to finitely generated subgroups, the given bijection was easily computable, and it soon became clear that the whole approach admitted an appealing restatement in terms of automata (see~\mbox{\parencite{kapovich_stallings_2002,bartholdi_rational_2010}}). 
One of our goals in this paper is to complement both viewpoints
in order to gain insight on the underlying ideas.

In the algorithmic case, the basic idea is to start with the flower automaton (each of its petals reading one of the initial generators) of the given \emph{finite} generating set for the subgroup $H \leqslant \Free[A]$, and then keep making identifications (now called \emph{Stallings foldings}) between incident arcs with the same label until reaching a unique final deterministic object (which we call the \emph{reduced} or \emph{Stallings automaton} of $H$).
As we will see in \Cref{sec: list}, many classic results on free groups --- such as the Nielsen--Schreier theorem,
the Schreier index formula, Howson's property,
Marshall Hall Theorem, residual finiteness, etc. --- can be neatly deduced from this graphical interpretation.

Moreover, if the starting generating set is finite, then so is the final automaton. Hence, the former bijection makes correspond finitely generated (f.g.) subgroups with finite Stallings automata, and it is fully algorithmic %
in this case, turning it into a very powerful tool for attacking algorithmic problems on subgroups of free groups:
the subgroup membership problem (\MP),
the finite index problem (\FIP),
the computation of bases for intersections,
or
checking the normality and malnormality of a given f.g.~subgroup (among many other algorithmic problems) can be successfully solved using Stallings' machinery.

This is a fairly well known topic, and several different presentations of its construction and main applications can be found in the literature. After the 
(topologically flavored) original paper \parencite{stallings_topology_1983} by Stallings, we should mention the papers \parencite{bartholdi_rational_2010}, and especially \parencite{kapovich_stallings_2002}, where a detailed account (including proofs) of the construction and some of the results in this survey can be found. Our idea is to complement and update these presentations in a somewhat lighter but still formal way, keeping the precise definitions and statements, but only sketching the ideas behind the proofs.

More precisely, our main goals are listed below:
\begin{itemize}[beginpenalty=10000]
    \item state in a clear and precise way the main ingredients (definitions, theorems, \ldots) of the theory, and set a convenient notation for them. This is done in the preliminary sections for the general (common) concepts, and in the corresponding subsections of \Cref{sec: list} for the specific ones;
    
    \item keep the discussion lighthearted but still rigorous, making it a quick reference for experts, but also a suitable --- and hopefully enjoyable and stimulating --- read for students;
    
    \item list a reasonably comprehensive\footnote{We must emphasize that, at the time of writing, a search in a standard research engine returns dozens of results for this topic; therefore, this selection will necessarily be biased towards the research interests (and knowledge) of the authors.} and updated selection of applications of this theory together with a (somewhat informal) explanation of the ideas behind the results;
    \item provide some examples helping to develop the intuition behind the results.
\end{itemize}

It is worth mentioning that the tremendous success of the theory of Stallings automata has led to many attempts to extend it to broader contexts; first by Stallings himself, to groups acting non-freely on graphs and trees   \parencite{stallings_foldings_1991}, and later successively refined and extended in the graph of groups context by multiple authors, see \eg 
\parencite{bestvina_bounding_1991, dunwoody_folding_1998,dunwoody_small_1999, rips_cyclic_1997, sela_acylindrical_1997, sela_diophantine_2001,
 dicks_equalizers_1999,dunwoody_groups_1997,guirardel_approximations_1998,guirardel_reading_2000,bowditch_peripheral_2001,bogopolski_uniqueness_2002}. 

Other (mostly algorithmically-oriented) generalization targets of Stallings automata include 
inverse monoids and semigroups \parencite{margolis_free_1993,margolis_closed_2001,delgado_combinatorial_2002,steinberg_inverse_2002}, groups satisfying certain small-cancellation properties~\parencite{arzhantseva_class_1996}, fully residually free groups~\parencite{myasnikov_fully_2006,kharlampovich_subgroups_2004,nikolaev_finite_2011}, free products~\parencite{ivanov_intersection_1999}, certain graphs of groups ~\parencite{kapovich_foldings_2005}, amalgams of finite groups~\parencite{markus-epstein_stallings_2007}, groups acting freely on $\ZZ^n$-trees~\parencite{nikolaev_membership_2012}, virtually free groups~\parencite{silva_finite_2016}, quasi-convex subgroups~\parencite{kharlampovich_stallings_2017}, free-abelian by free groups~\parencite{delgado_algorithmic_2013,delgado_extensions_2017,delgado_stallings_2022}, CAT(0) cube complexes~\parencite{beeker_stallings_2018},
certain relatively hyperbolic groups~\parencite{kharlampovich_generalized_2020}, right-angled Coxeter groups~\parencite{dani_subgroups_2021}, etc. We remark, however, that in the present document we restrict ourselves to applications of the classical Stallings theory to subgroups of free groups.

\medskip

In \Cref{sec: words} we introduce the basic terminology and notation regarding formal languages, which is developed in order to precisely define free groups, for which standard definitions and  some of its main properties are presented.

A similar approach is used in \Cref{sec: automata}, this time starting with the concept of (di)graph, to arrive at a precise definition of Stallings automata and describe their relation with the subgroups of the free group. 
Finally, \Cref{sec: list} contains a summary of results (together with explanations) on the applications of Stallings theory on (different aspects of) subgroups of the free group. This is the main section of the survey.

Some of the sections are marked with an asterisk. These correspond to rather advanced or specific topics, for which little more than the results and necessary definitions are given; we refer the reader to the original work for details.

\subsection{\bf Notation, terminology, and conventions}
Most of the notation and terminology used throughout this survey is standard;  however, we clarify below some conventions that may cause confusion.

The set of natural numbers, denoted by $\NN$, is assumed to contain zero, and we specify conditions on it using subscripts; for example, we denote by $\NN_{\geq 1}$ the set of strictly positive integers.
The cardinal of a set $S$ is denoted by $\card S$, whereas the notation $|\,\cdot\,|$ is reserved to denote length (in different contexts).
Latin letters $m,n,r,
\ldots$ are used to denote finite cardinals, whereas Greek letters $\kappa,\mu,\ldots$ denote arbitrary cardinals.

We denote by $\Free$ a generic free group, whereas the notations $\Free[A]$ and $\Free[\kappa]$ (\resp $\Fn$) are used to emphasize a chosen basis $A$ and its rank $\kappa$ (\resp finite rank $n$), respectively. The first lowercase letters of the Latin alphabet ($a,b,c,\ldots$) are commonly used to denote symbols in our formal alphabets, whereas the latter ones ($u,v,w,\ldots$) usually denote formal words or elements in the free group. 

Functions are assumed to act on the right. That is, we denote by $(x)\varphi$ (or simply by $x \varphi$) the image of the element $x$ by the homomorphism $\varphi$, and we denote by  $\varphi \psi$ the composition \smash{$A \xto{\varphi\phantom{\psi}\hspace{-8pt}} B \xto{\psi} C$}. Accordingly, we write, for example, $g^{h} = h^{-1} gh$ and $[g,h] = g^{-1} h^{-1} g h$. We sometimes write $A \into B$ (\resp $A \onto B$) to express that a function $A \to B$ is injective (\resp surjective),
and we denote by $(C)\varphi\preim$ the full preimage of $C \subseteq B$ by the map $\varphi\colon A \to B$.
We denote by $X \pto Y$ a partial function from $X$ to $Y$, and by $Y^{\subseteq X}$ the set of such partial functions.

Throughout the paper we write $H \leqslant\fg G$, $H \leqslant\ff G$, $H \leqslant\fin G$, $H \leqslant\alg G$ to denote that the subgroup $H$ is finitely generated, a free factor, of finite index, and algebraic in $G$, respectively.

For the different variants of graphs appearing in the paper, the following convention applies: uppercase Greek letters ($\Gamma, \Delta, \ldots$) denote unlabelled undirected graphs, the arrow accentuated versions ($\dGri,\dGrii,\ldots$) denote unlabelled digraphs, and the boldface versions $(\Ati,\Atii, \ldots)$ denote labelled digraphs (including automata). In a similar vein, we distinguish between (undirected) \emph{paths} and \emph{cycles}, and (directed and allowing repetitions) \emph{walks} and \emph{rounds} respectively.

Regarding computability, the terms computable (or decidable) and computably enumerable are preferred to the also very common recursive (or solvable) and recursively enumerable, respectively.

\section{\bf Letters, words, languages, \ldots\ and groups} \label{sec: words}

Throughout the paper, we denote by $A = \set{a_i}_{i}$ a
set called \defin{alphabet}, whose elements are called  \defin[letter]{letters}. Depending on the context, this set may be finite or not. 
A \defin{word} or \emph{string} over $A$ is a finite sequence of elements in $A$. If $w = (a_{i_1},\ldots ,a_{i_l}) \in A^l$ then we write $w = a_{i_1} \cdots a_{i_l}$, and  we say that $w$ has \defin{length} $l$, denoted by $|w| =l$. In particular, the (unique) word of length $0$ is called the \defin{empty word} and it is denoted by $\emptyword$. We denote by $A^{\!*}$ the \defin{free monoid} on $A$, consisting of all the words on $A$ (including the empty word, working as the neutral element) endowed with the concatenation operation. 

\begin{rem}
If the alphabet $\Alfi$ is finite, then the free monoid $\Alfi^{\!*}$ has countable cardinal $\aleph_0$; otherwise,  $\Alfi$ and $\Alfi^{\!*}$ have the same cardinal. So, a nontrivial free monoid is always infinite.
\end{rem}

We say that a nonempty word $u\in A^*$ is a \defin{factor} (\resp \defin{proper factor}) of a word $w\in A^*$ if $w=vuv'$ for some words $v,v'\in A^*$ (\resp at least one of them being non-trivial); specifically, if $v=1$ (\resp  $v'=1$), $u$ is called a \defin{prefix} (\resp \defin{suffix}) of~$w$.

If $A$ is an alphabet, the subsets of $A^{\!*}$ are called \defin[language]{languages} (over $A$), or $A$-languages. A language over a finite alphabet $A$ is said to be \defin[rational language]{rational} if it can be  obtained from the (finitely many) letters in $A$ using the operators of union, product and star (\ie submonoid generated by a language) finitely many times. A direct consequence of a fundamental theorem of Kleene (see for example \parencite{sakarovitch_elements_2009}) is that the set of rational $A$-languages is closed under finite intersection and complement.

In the context of groups (and free groups in particular), however, a very particular kind of alphabets have a leading role. Given an alphabet $A$,
we  denote by~$A^{-1}$ the \defin{set of formal inverses} of~$A$. Formally, $A^{-1}$ can be defined as a new set $A\!'$ equipotent and disjoint with $A$, together with a bijection $ A \to A\!'$. Then, for every $a \in A$, we write $a^{-1}$ the image of~$a$, and we call it the \defin{formal inverse} of $a$. So, $A^{-1}=\{ a^{-1} \st a \in A \}$ and $A \cap A^{-1} = \varnothing$. Then, the set $A^{\pm} = A \sqcup A^{-1}$, called the \defin{involutive closure} of $A$, is equipped with an involution $^{-1} 
$ (by defining $(a^{-1})^{-1} = a$), which can be extended to
$(A^{\pm})^*$ in the natural way: $(a_{i_1}\cdots a_{i_l})^{-1} = a_{i_l}^{-1}\cdots a_1^{-1}$, for all $a_{i_1}, \ldots,a_{i_l} \in A^{\pm}$. An alphabet is called \defin[involutive alphabet]{involutive} if it is the involutive closure of some other alphabet.

In such an involutive ambient, a word in $(A^{\pm})^{*}$ (which we just refer to as a \defin{word on $A$}) is said to be \defin[(freely) reduced word]{(freely) reduced} if it contains no consecutive mutually inverse letters (\ie it has no factor of the form $a a^{-1}$, where $a \in A^{\pm}$). It is well known that
the word obtained from $w \in (A^{\pm})^{*}$ by successively removing pairs of consecutive inverse letters (in any possible order) is unique; we call it the \defin{free reduction} of $w$, and we denote it by $\red{w}$. Similarly, we write $\red{S} = \set{\red{w} \st w \in S}$, for any subset $S\subseteq (A^{\pm})^{*}$, and we denote by $\Red_A$ the rational language
 \begin{equation*}
\Red_A \,=\, \red{(A^{\pm})^{*}} \,=\, (A^{\pm})^{*} \setmin \textstyle{\bigcup_{a \in A^{\pm}} (A^{\pm})^{*} a a^{-1} (A^{\pm})^{*}, }
 \end{equation*}
called  the \defin{set of reduced words} on $A$. In a similar vein, a word $w$ on $A$ is said to be \defin[cyclically reduced word]{cyclically reduced} if all of its cyclic permutations are reduced (that is, if $w^2$ is reduced). The \defin{cyclic reduction} of a word $w \in \Free[A]$, denoted by $\cred{w}$, is obtained after iteratively removing from $\red{w}$ the first and last letters if they are inverses of each other. We also extend this notation to subsets; \ie we write $\cred{S} = \set{\cred{w} \st w \in S}$, for any subset $S\subseteq (A^{\pm})^{*}$, and denote by $\Cred_A$ (or $\Cred_n$) the set of all cyclically reduced words in $(A^{\pm})^{*}$, that is the rational language
 \begin{equation*}
\Cred_A \,=\, \cred{(A^{\pm})^{*}} \,=\, (A^{\pm})^{*} \setmin \textstyle{\bigcup_{a \in A^{\!\pm}} \left( (A^{\pm})^{*} a a^{-1} (A^{\pm})^{*}
\cup a (A^{\pm})^{*} a^{-1} \right)}.
 \end{equation*}
If $\card A = n$, we usually write $\Red_A = \Red_n$ and $\Cred_A = \Cred_n$.

\subsection{\bf The free group}
The free group can be introduced from several points of view: combinatorial, categorical, geometric, \ldots\ From the combinatorial point of view used so far, it can be seen as the quotient of the free involutive monoid $(A^{\pm})^*$ by the congruence generated by cancellations between respectively inverse letters.

\begin{defn}
Let $A$ be an arbitrary set (alphabet). 
The \defin[${\Free[A]}$]{free group on $A$}, denoted by $\Free[A]$, is the quotient
 \begin{equation} \label{eq: free group def}
\Free[A] \,=\, (A^{\pm})^* / \mathcal{C} \,,
 \end{equation}
where $\mathcal{C}$ is the congruence generated by the pairs $\set{(a a^{-1}, \emptyword) \st a \in A^{\pm}}$, with the product inherited from concatenation in $(A^{\pm})^*$. We say that a group is \emph{free on $A$} (or \defin[free basis]{free with basis $A$}) if it is isomorphic to $\Free[A]$, and we say that a group is \defin[free group]{free} if it is isomorphic to $\Free[A]$, for some set $A$.

Furthermore, it is easy to see that two free groups $\Free[A]$ and $\Free[B]$ are isomorphic if and only if $\card A = \card B$. That is, the cardinal of $A$ is an algebraic invariant of $\Free[A]$; we will see that this cardinal is precisely the rank of $\Free[A]$, \ie $\rk(\Free[A])=\card A$. A proof (for both facts) based on Stallings machinery will be given in \Cref{thm: isom iff card}.
\end{defn}

\begin{rem}
Although most of them can be overcome by assuming the axiom of choice, in order to avoid set-theoretical subtleties, throughout the paper, we will assume that ${\card A = \kappa \leq \aleph_0}$. That is, we will consider groups $\Free = \Free[A]$ of at most countable rank.
\end{rem}

\begin{rem}
Note that $\Free[\varnothing]$ is the trivial group, and $\Free[\hspace{1pt}1]$ is isomorphic to the group~$\ZZ$ of integers. It is not hard to see (see \eg \cite[Ch.\,1]{johnson_presentations_1997}) that these are the only abelian (and somewhat unrepresentative) free groups.
\end{rem}

\begin{lem}
Let $u,v \in \Free$. Then, $uv=vu$ if and only if $u=w^p$ and $v=w^q$, for some $w\in \Free$ and some $p,q\in \ZZ$.
\end{lem}

It is a folklore result (see \eg \parencite{bogopolski_introduction_2008})
that the set $\Red_A$ (of reduced words on~$A$) is a complete set of representatives for $\Free[A]$ (as defined in \eqref{eq: free group def}). Therefore, we will take $\Red_A$ as a set of (computable) normal forms for the elements in~$\Free[A]$. That is, if we identify every reduced word $w\in \Red_A$ with the element $[w] \in \Free[A]$ then the free group $\Free[A]$ can be thought as the set $\Red_A$ with the operation consisting of ``concatenation followed by reduction'', \ie for every $u,v\in \Red_A$, $uv=\red{uv}$.  We shall use this interpretation throughout the paper.

In particular, a word $w \in (A^\pm)^*$ represents the trivial element in $\Free$ if and only if $\red{w} = \emptyword$. Since word-reduction is obviously algorithmic, this remark immediately solves the word problem for free groups.

We further recall that, with this identification, Benois' Theorem (see~\parencite{benois_parties_1969}) allows us to understand the rational subsets of $\Free[A]$ as reductions of rational $A^{\pm}$-languages.

The universal property below turns out to be fundamental as it characterizes (group) freeness, and can consequently be taken as an alternative (categorical) definition for a free group (see \eg \parencite{bogopolski_introduction_2008}).

\begin{thm} \label{thm: free cat}
A group $F$ is free on $A \subseteq F$ if and only if every function from $A$ to an arbitrary group~$G$ can be extended to a homomorphism $F \to G$ in a unique possible way.
\end{thm}

We note that the extendability in \Cref{thm: free cat} corresponds to the set $A$ being \defin[freely independent set]{freely independent} in~$F$ (\ie a basis of $\gen{A} \leqslant F$), and the uniqueness of the extension corresponds to $A$ 
generating the whole~$F$.

Several alternative characterizations for free groups do exist. We state below two important geometrical ones along the lines of the present article (details can be found in~\parencite{bogopolski_introduction_2008}).  

\begin{thm} \label{thm: free geom}
Let $F$ be a group. Then, the following statements are equivalent:
\begin{enumerate}[dep]
\item $F$ is a free group;
\item \label{item: free iff graph} $F$ is the fundamental group of a connected graph; %
\item \label{item: free iff tree} $F$ acts freely and without inversion of edges on a tree.%
\end{enumerate}
\end{thm}
An important corollary of these graphical characterizations is the classical Nielsen--Schreier Theorem claiming that every subgroup of a free group is again free. This result is immediate from \ref{item: free iff tree} and also easily derivable from \ref{item: free iff graph} using the Stallings construction (see \Cref{thm: Nielsen-Schreier}).

\begin{defn}\label{defn: extension}
Let $H,K$ be subgroups of a free group $\Free$. If $H \leqslant K$ then we say that $K$ is an \emph{extension} (or \defin{overgroup}) of $H$ and that $H \leqslant K$ is an \defin{extension of subgroups}. An extension of subgroups $H\leqslant K\leqslant \Free[\Alfi]$ is called \defin[free extension]{free} (we also say that $H$ is a \defin{free factor} of $K$, denoted by $H\leqff K$), if some (and hence any) basis of $H$ can be extended to a basis for $K$. In particular, $H\leqff K$ implies $\rk(H)\leqslant \rk(K)$. A kind of dual notion called algebraicity  will be defined in~\Cref{ssec: algebraic extensions}.
\end{defn}

\subsection{Groups and presentations}

Note that given a (any) group $G$ and a generating set $A \subseteq G$ for~$G$, we can apply \Cref{thm: free cat} to the inclusion map $A \into G$ to conclude that there exist a unique extension epimorphism $\Free[A] \onto G$, thus obtaining the crucial property below, which, in part, explains the importance and ubiquity of free groups in group theory. 

\begin{cor} \label{cor: G = F/N}
Every group is (isomorphic to) a quotient of some free group.
\end{cor}

More concretely, if we call $N$ the kernel of the map $\Free[A] \onto G$,
we obtain that $G \isom \Free[A]/N$, which (after choosing a subset $R\subseteq N$ generating $N$ \emph{as a normal subgroup}, \ie such that $\normalcl{R} = N$) leads us to the fundamental concept of \emph{group presentation}.

\begin{defn}
Let $G$ be a group. A \defin[group presentation]{presentation} for $G$ is a pair $(A,R)$ such that $A$ is a set, $R \subseteq \Free[A]$, and $G\isom \Free[A] / \normalcl{R}$. Then, we say that the group $G$ admits the presentation $\pres{A}{R}$, and we write $G \isom \pres{A}{R}$, or $G = \pres{A}{R}$ (slightly abusing the notation).
\end{defn}

Note that from \Cref{cor: G = F/N}, every group admits a presentation (in fact infinitely many). 
It is also clear from the definition that if $G = \pres{A}{R}$, then $A$ is (corresponds to) a generating set for the group $G$, and $R$ is (corresponds to) a complete set of \defin{relations}
in $G$. Accordingly, we call $A$ and $R$ the sets of \defin[presentation generators]{generators} and \defin[presentation relators]{relators} of the presentation $\pres{A}{R}$, respectively. Also, we say that a group $G$ is \emph{finitely generable} or \defin{finitely generated}\footnote{This more ambiguous terminology is standard.} (f.g.) --- \resp \emph{finitely presentable} or \defin{finitely presented} (f.p.) --- if $G$ admits a presentation  $\pres{A}{R}$, where $A$ is --- \resp both $A$ and $R$ are --- finite.

\begin{rem}
A group is free if and only if it admits a presentation with no relators. Concretely $\Free[A]$ admits the presentation $\pres{A}{-}$, which is called the \defin{standard presentation} for $\Free[A]$, and will be assumed by default throughout the paper.
\end{rem}

From the previous discussion it should be clear that every group $G$ admits the following sequence of epimorphisms (the first one, $\mu$, of monoids; and the second one, $\rho$, of groups)
 \begin{equation}
\begin{array}{rcccl}
(A^\pm)^*& \xonto{\mu\,} & \Free[A] & \xonto{\rho\,} & G  \\ w & \longmapsto &  \red{w} & \longmapsto & [w]_G = [\red{w}]_G  \,,
\end{array}
 \end{equation}
where $A$ is a generating set for $G$. To avoid verbosity, if the context is clear, we shall often abuse language and use words in the generators as if they were the elements they represent in the groups. For example, we may write $bbab^{-1}ba \in \Free_{\set{a,b}}$; or, if $S \subseteq (A^\pm)^*$, we may write $\gen{S}_G$ to refer to the subgroup generated in $G$ by the elements (in $G$) represented by the words in $S$.
 
\begin{defn}
    The \defin{rank of a group} $G$, denoted by~$\rk (G)$, is the smallest cardinality of a generating set for $G$.
\end{defn}

\section{\bf Graphs, digraphs, automata, \ldots\ and subgroups} \label{sec: automata}

As we have mentioned, we shall use automata (\ie essentially labelled digraphs) to describe the subgroups of a free group. For technical convenience, we define digraphs in the sense of Serre.

\begin{defn} \label{def: digraph}
A \emph{directed graph} (also called \defin{digraph}, for short) is a tuple $\dGri =(\Verts, \Edgs, \init,\term)$, where $\Verts$ and $\Edgs$ are disjoint sets (called the set of \defin{vertices} of~$\dGri$ and the set of \defin{arcs} (or \emph{directed edges}) of $\dGri$, respectively), and $\init,\term \colon \Edgs \to \Verts$ are maps assigning to each arc in $\dGri$ its \defin[arc origin]{origin} (or \emph{initial vertex}) and \defin[arc end]{end} (or \emph{terminal vertex}), respectively.  
\end{defn}

Note that both the possibility of arcs having the same vertex as origin and end (called \defin[loop]{loops}), and of multiple arcs sharing the same origin and end (called \defin[parallel arc]{parallel arcs}) are included in \Cref{def: digraph}.

We denote by $\Verts \dGri$ and $\Edgi \dGri$ the respective sets of vertices and arcs of a digraph~$\dGri$. A digraph $\dGri$ is called \defin[finite digraph]{finite}, \defin[countable digraph]{countable}, etc.\ if the cardinal $\card{(\Verts \dGri \sqcup \Edgi \dGri)}$
is so.

\begin{defn}
A \defin{(directed) walk} in a digraph $\dGri$ is a finite alternating sequence $\walki = \verti_0 \edgi_1 \verti_1 \cdots \edgi_{l} \verti_{l}$, where $\verti_i \in \Verts \dGri$, $\edgi_i \in \Edgs \dGri$,  $\init{\edgi_i} = \verti_{i-1}$ and $\term{\edgi_i}=\verti_{i}$ for $i = 1,\dots,l$. Then, $\verti_0$ and $\verti_l$ are called the \defin[walk origin]{origin} (or \emph{initial vertex}) and \defin[walk end]{end} (or \emph{terminal vertex}) of $\walki$, respectively, and we say that $\walki$ is a walk from $\verti_0$ to $\verti_l$ (a \defin{$(\verti_0,\verti_l)$-walk} for short), and we write $\walki \colon \smash{\verti_0 \xwalk{\ } \verti_l}$. We write $\verti \xwalk{\ } \vertii$ if there exists a walk from $\verti$ to $\vertii$. If the first and last vertices of $\walki$ coincide then we say that $\walki$ is a \defin{closed walk} or a \emph{round}. A closed walk from $\verti$ to $\verti$ is called a \defin[round]{$\verti$-round}. The \defin{length of a walk} $\walki$, denoted by $|\walki|$, is the number of arcs in~$\walki$ (counting possible repetitions). The walks of length $0$ are called \defin[trivial walk]{trivial walks}, and correspond %
to the vertices in $\dGri$. %
\end{defn}

\begin{defn}
Let $A$ be an alphabet.
An \defin[labelled digraph]{$\Alfi$-labelled digraph} (an \defin{$\Alfi$-digraph}, for short) is a pair $\Ati = (\dGri,\lab)$, where $\dGri = (\Verts, \Edgs, \init,\term)$ is a digraph, and $\lab \colon \Edgs \to \Alfi$ is a map assigning to every arc in~$\dGri$ a label from  $\Alfi$.
\end{defn}

If $\edgi$ is an arc from $\verti$ to $\vertii$ and $\lab(\edgi) = \alfi$, then we write $\verti \xarc{\alfi\,} \vertii$ and 
we say that $\alfi$ is the \defin[arc-label]{label} of $\edgi$, and that $\edgi$ is an \defin{$a$-arc}.

An $\Alfi$-labelling on arcs extends naturally to an $\Alfi^{\!*}$-labelling on walks by concatenating the corresponding arc labels; \ie if $\walki = \verti_0 \edgi_1 \verti_1 \cdots \edgi_{l} \verti_{l}$ is a walk on an \mbox{$\Alfi$-digraph}, then $\lab(\walki) = \lab(\edgi_1) \cdots \lab(\edgi_{l}) = w \in \Alfi^{\!*}$, and the label of any trivial walk is the empty word~$\emptyword$. Then, we say that the walk $\walki$ \defin{reads} (or \defin{spells}) the word $\wordi$, that the word~$\wordi$ \defin{labels} the walk $\walki$,
and we write $\walki\colon \verti_0 \xwalk{_{\scriptstyle{w}}} \verti_l$.

\begin{defn}
Let $\Ati$ be an $\Alfi$-digraph and let $P,Q$ be subsets of vertices in~$\Ati$. Then, the set of words read by walks from vertices in $P$ to vertices in $Q$ is called the \defin{language from $P$ to $Q$} (in $\Ati$), and is denoted by $\Lang_{P,Q}(\Ati) $. For simplicity, if $\verti,\vertii \in \Verts\Ati$, then the languages from $\{ \verti \}$ to $\{ \vertii \}$, and from $\{ \verti \}$ to $\{ \verti \}$, are denoted by $\Lang_{\verti,\vertii}(\Ati)$ and~$\Lang_{\verti}(\Ati)$, respectively.  
\end{defn}

When distinguished sets of initial and terminal vertices are fixed in a labelled digraph $\Ati$, we say that $\Ati$ is an automaton (see \Cref{def: automaton}).

A vertex $\verti$ in an $\Alfi$-digraph $\Ati$ is said to be \defin[saturated vertex]{saturated} if for every letter $a \in \Alfi$ there is (at least) one $a$-arc with origin $\verti$.
Otherwise,
we say that both $\verti$ and $\Ati$ are unsaturated (or \defin[deficient vertex]{$a$-deficient}, if we want to allude a missing label).
The \defin{$a$-deficit} of $\Ati$, denoted by $\defc[a]{\Ati}$ is the number (cardinal) of $a$-deficient vertices in $\Ati$. 
An $\Alfi$-digraph is said to be \defin[saturated $\Alfi$-digraph]{saturated} (or \emph{complete}\footnote{The term `complete' is standard in automata theory; we use the term `saturated' instead in order to avoid confusion with the notions of `complete graph' and `complete digraph'.}) if all its vertices are so, that is if $\defc[a]{\Ati} = 0$, for all $a \in A$.

An $A$-digraph $\Ati$ is said to be \defin[deterministic vertex]{deterministic} at a vertex $\verti \in \Verts \Ati$ if no two arcs with the same label depart from $\verti$. An $A$-digraph is said to be \defin[deterministic $\Alfi$-digraph]{deterministic} if it is deterministic at every vertex; that is, if for every vertex $\verti \in \Verts \Ati$, and every pair of arcs $\edgi,\edgi'$ leaving $\verti$, $\lab(\edgi) = \lab(\edgi')$ implies $\edgi = \edgi'$.

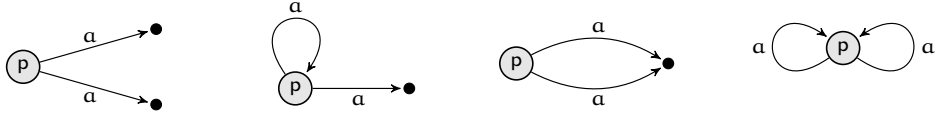
\begin{figure}[H]
\centering
  \begin{tikzpicture}[shorten >=1pt, node distance=.5cm and 1.75cm, on grid,auto,>=stealth']
   \node[state,semithick, fill=gray!20, inner sep=2pt, minimum size = 10pt] (1) {$\scriptstyle{\verti}$};
   \node[state] (2) [above right = of 1] {};
   \node[state] (3) [below right = of 1] {};

   \path[->]
        (1) edge[]
            node[pos=0.45,above] {$\alfi$}
            (2)
            edge[]
            node[pos=0.45,below] {$\alfi$}
            (3);
\end{tikzpicture}
\hspace{20pt}
\begin{tikzpicture}[shorten >=1pt, node distance=.5cm and 1.5cm, on grid,auto,>=stealth']
   \node[state,semithick, fill=gray!20, inner sep=2pt, minimum size = 10pt] (1) {$\scriptstyle{\verti}$};
   \node[state] (2) [right = 1.5 of 1] {};

   \path[->]
        (1) edge[loop above,min distance=12mm,in= 90-35,out=90+35]
            node[] {\scriptsize{$a$}}
            (1)
            edge
            node[pos=0.5,below] {$\alfi$}
            (2);
\end{tikzpicture} 
\hspace{25pt}
\begin{tikzpicture}[shorten >=1pt, node distance=.5cm and 2cm, on grid,auto,>=stealth']
   \node[state,semithick, fill=gray!20, inner sep=2pt, minimum size = 10pt] (1) {$\scriptstyle{\verti}$};
   \node[state] (2) [right = of 1] {};

   \path[->]
        (1) edge[bend right]
            node[pos=0.5,below] {$\alfi$}
            (2)
        (1) edge[bend left]
            node[pos=0.5,above] {$\alfi$}
            (2);
   \end{tikzpicture}
\hspace{20pt}
\begin{tikzpicture}[shorten >=1pt, node distance=.5cm and 1.5cm,
 on grid,auto,>=stealth']
   \node[state,semithick, fill=gray!20, inner sep=2pt, minimum size = 10pt] (1) {$\scriptstyle{\verti}$};

   \path[->]
        (1) edge[loop right,min distance=12mm,in=35,out=-35]
            node[] {\scriptsize{$a$}}
            (1)
             edge[loop left,min distance=12mm,in=180-35,out=180+35]
            node[] {\scriptsize{$a$}}
            (1);
\end{tikzpicture}   
\caption{Nondeterministic situations at vertex $\verti$}
\end{figure}
If $\Ati$ is deterministic, then for every vertex $\verti$ in $\Ati$ and every word $w\in \Alfi^*$ there is at most one walk in $\Ati$ reading $w$ from $\verti$; we denote by $\verti w$ its final vertex in case it exists (otherwise, $\verti w$ is undefined). That is, we have a (partial) \defin{transition function} $\tau \colon \Verts \Ati \times \Alfi^{\!*} \pto \Verts \Ati$, $(\verti,w) \mapsto \verti w$, whose restrictions to particular $w$'s define partial functions $\tau_{w} \colon \Verts \Ati \pto \Verts \Ati$, for every $w \in \Alfi^{\!*}$. It is not difficult to see that the corresponding function $w \mapsto \tau_w$ is indeed a homomorphism (of monoids) $\tau_{*} \colon \Alfi^{\!*} \to \Verts\Ati^{\subseteq\Verts \Ati}$ (from the free monoid over $\Alfi$, to the monoid of partial functions from $\Verts \Ati$ to $\Verts \Ati$ with composition), whose image $\tau_{*}(\Alfi^*)$ is called the \defin{transition monoid} of~$\Ati$, denoted by $\TM{\Ati}$.

Note also that $\Ati$ is fully recoverable from the set $\set{\tau_{\alfi} \st \alfi \in \Alfi}$ of partial functions corresponding to the action on $\Verts \Ati$ of the letters from the alphabet $\Alfi$ (this alternative way of defining a deterministic automaton will be useful in \Cref{sec: asymptotic}). Finally, note that if $\Ati$ is both deterministic and saturated then the function $\tau_w\colon \Verts \Ati \to \Verts \Ati$ is \emph{total} for every word $w \in \Alfi^{\!*}$. In this case, $\Ati$ is characterized by the set of total functions $\set{\tau_{\alfi}: \alfi \in \Alfi}$.

\subsection{Involutive digraphs and automata} \label{ssec: involutive automata}

Recall that, for an arbitrary alphabet $\Alfi$, we write $\Alfi^{-1} = \set {\alfi^{-1} \st \alfi \in \Alfi}$ the \defin[formal inverse]{set of formal inverses} of~$\Alfi$, and $\Alfi^{\pm} = \Alfi \sqcup \Alfi^{-1}$ the \defin{involutive closure} of $\Alfi$. 

\begin{defn} \label{def: involutive digraph}
An \defin{involutive $\Alfi$-digraph} is an $\Alfi^{\pm}$-digraph
with a labelled involution
${\edgi \mapsto \edgi^{-1}}$
on its arcs; \ie to every arc
$\smash{\edgi \equiv \verti \xarc{\,\alfi\ }\vertii}$ we associate a unique arc $\smash{\edgi^{-1} \equiv \verti \xcra{\,\alfi^{\text{-}1}\!} \vertii}$
(called the \defin[inverse arc]{inverse} of $\edgi$) such that $(\edgi^{-1})^{-1} = \edgi$.
\end{defn}

That is, in an involutive $\Alfi$-digraph, labelled arcs appear by (mutually inverse) pairs.

\begin{figure}[H]
\centering
\begin{tikzpicture}[shorten >=3pt, node distance=.3cm and 2cm, on grid,auto,>=stealth']
   \node[state] (0) {};
   \node[state] (1) [right = of 0] {};
   \path[->]
        ([yshift=0.3ex]0.east) edge[]
            node[pos=0.5,above=-.2mm] {$a$}
            ([yshift=0.3ex]1)
   ([yshift=-0.3ex]1.west) edge[dashed]
            node[pos=0.5,below=-.2mm] {$a^{-1}$}
            ([yshift=-0.3ex]0);
\end{tikzpicture}
\vspace{-5pt}
\caption{An arc in an involutive labelled digraph with its negative part dashed}
\label{fig: involutive arc}
\end{figure}
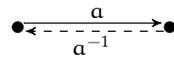

An arc in an involutive $A$-digraph $\Ati$ is said to be \defin[positive arc]{positive} (\resp \defin[negative arc]{negative})
if it is labelled by a letter in~$A$ (resp.~$A^{-1}$).
We respectively denote by~$\Edgi^{+} (\Ati)$ and $\Edgi^{-} (\Ati)$
the sets of positive
and negative
arcs in an involutive digraph~$\Ati$.
The \defin{positive part} $\Ati^{+}$ (\resp \defin{negative part} $\Ati^{-}$) of $\Ati$ is the $A$-digraph obtained after removing all the negative (\resp positive) arcs from $\Ati$.
It is clear that $\Ati$ is fully characterized by $\Ati^{+}$ (or $\Ati^{-}$), which we will typically use to describe involutive automata. That is, we will usually represent involutive $\Alfi$-automata through their positive (\ie $\Alfi$-labelled) part, with the convention that a positive arc $\verti \xarc{\alfi\,} \vertii$ reads the inverse label $\alfi^{-1}$ when crossed backwards (\ie from $\vertii$ to $\verti$); see \Cref{fig: involutive arc}.

The \defin[inverse walk]{inverse of a walk} $\walki = \verti_0 \edgi_1 \verti_1 \cdots \edgi_{l} \verti_{l}$  (reading $\lab(\edgi_1)\cdots\lab(\edgi_l)$) in an involutive digraph is the walk 
$\walki^{-1} = \verti_l \edgi_l^{-1} \verti_{l-1} \cdots \edgi_{1}^{-1} \verti_{0}$ (reading $\lab(\edgi_l)^{-1}\cdots\lab(\edgi_1)^{-1}$).

A walk in an involutive digraph is said to present \defin{backtracking} if it has two successive labelled arcs inverse of each other. A walk without backtracking is said to be \defin[reduced walk]{reduced}.
Note that if $\walki$ is a walk in an involutive $\Alfi$-digraph~$\Ati$, then the label $\lab(\walki)$ is a word in the involutive monoid~$(\Alfi^{\pm})^*$. We define the \defin{reduced label} of a walk $\walki$ in a labelled digraph $\Ati$ as $\rlab(\walki) = \red{\lab(\walki)} \in \Free[A]$.
It is clear that if the label $\lab(\walki)$ is reduced then the walk $\walki$ is reduced, but the converse is only true if $\Ati$ is deterministic.

Note that if $\Ati$ is deterministic and involutive, then the partial functions $\tau_{w}\colon \Verts \Ati \pto \Verts \Ati$ are indeed \defin[partial injections $\tau_{w}$]{partial injections},
and, for each $\alfi \in \Alfi$, $\tau_{\alfi^{\text{-}1}}$ and $\tau_{\alfi}$ are inverse of each other in the inverse semigroup sense: $ \tau_{\alfi} \tau_{\alfi^{\text{-}1}} \tau_{\alfi} = \tau_{\alfi}$ and $\tau_{\alfi^{\text{-}1}} \tau_{\alfi} \tau_{\alfi^{\text{-}1}} = \tau_{\alfi^{\text{-}1}}$. Then, the transition monoid $\TM{\Ati}$ is an inverse semigroup, and we say that $\Ati$ is an \defin{inverse automaton}. If $\Ati$ is also saturated, then the partial injections $\tau_{\alfi}$ (and~$\tau_{w}$) are indeed total bijections $\Verts \Ati \to \Verts \Ati$, the corresponding transition monoid is a subgroup of the symmetric group $\operatorname{Sym}(\Verts \Ati)$ called \defin{automaton group} (denoted by $G_{\Ati}$), and $\Ati$ is called a \defin{group automaton}.

Also, in the case of $\Ati$ being deterministic and finite, a straightforward adaptation of the \emph{handshaking lemma} for digraphs provides the following important property.

\begin{lem} \label{lem: =deficient}
Let $\Ati$ be an involutive, deterministic and finite $\Alfi$-digraph. Then, for every $\alfi \in \Alfi$, the number of
$\alfi$-deficient vertices in $\Ati$ equals the number of $\alfi^{-1}$-deficient vertices in $\Ati$, \ie $\defc[\alfi]{\Ati} = \defc[\alfi^{\text{-}1}]{\Ati}$.
\end{lem}

Note that the finiteness condition is essential in~\Cref{lem: =deficient}, as the following example shows.
\begin{figure}[H] 
  \centering
  \begin{tikzpicture}[shorten >=1pt, node distance=1.2 and 1.2, on grid,auto,>=stealth']
   \node[state] (0) {};
   \node[state] (1) [right = of 0]{};
   \node[state] (2) [right = of 1]{};
   \node[state] (3) [right = of 2]{};
   \node[state] (4) [right = of 3]{};
   \node[] (5) [right = 0.8 of 4]{};
   \node[] (dots) [right = 0.2 of 5]{$\cdots$};

   \path[->]
        (0) edge[loop above,red,min distance=10mm,in=55,out=125]
            node[left = 0.1] {\scriptsize{$b$}}
            (0)
            edge[blue]
            node[below] {\scriptsize{$a$}}
            (1);

    \path[->]
        (1) edge[loop above,red,min distance=10mm,in=55,out=125]
            node[left = 0.1] {}
            (1)
            edge[blue]
            (2);

    \path[->]
        (2) edge[loop above,red,min distance=10mm,in=55,out=125]
            (2)
            edge[blue]
            (3);

    \path[->]
        (3) edge[loop above,red,min distance=10mm,in=55,out=125]
            (3)
            edge[blue]
            (4);
    \path[->]
        (4) edge[loop above,red,min distance=10mm,in=55,out=125]
            (4)
            edge[blue]
            (5);

\end{tikzpicture}
\caption{An infinite deterministic automaton $\Ati$ with $\defc[\alfi]{\Ati}=0$ but $\defc[\alfi^{\text{-}1}]{\Ati}=1$}
\label{fig: Infinte automaton with different deficit}
\end{figure}

If we ignore the labelling and identify all the mutually inverse pairs in an involutive digraph $\Ati$, we obtain an undirected graph called the \defin{underlying graph} of $\Ati$, usually denoted by the same letter unboldfaced (\eg $\Gri$). It is clear that every \defin{undirected graph} can be obtained in this way. Involutive automata inherit terminology from its underlying graph; for example, the rank $\rk(\Ati)$
of $\Ati$ is defined to be the rank of its underlying graph, and we say that an involutive digraph is \emph{connected}, \emph{vertex transitive}, a \emph{path}, a \emph{tree}, etc. if its underlying graph is so. Similarly, the \defin[vertex degree]{degree} of a vertex in an \emph{involutive} digraph is defined to be its degree --- \ie the number of edges incident to it --- in the underlying undirected graph.

In our context, an automaton is an involutive labelled digraph where certain sets of initial and terminal vertices are distinguished.

\begin{defn} \label{def: automaton}
An \defin[automaton]{$\Alfi$-automaton} is a tern $\att{\Ati\!}{P}{Q}
= (\Ati,P,Q)$, where $\Ati$ is an $\Alfi$-digraph, and $P$ and $Q$ are distinguished nonempty sets of vertices of~$\Ati$, called the sets of \defin[initial vertex]{initial} and \defin[terminal vertex]{terminal} vertices of $\att{\Ati\!}{P}{Q}$, respectively. A vertex in an automaton is called \defin[distinguished vertex]{distinguished} if it is either initial or terminal. A walk in $\att{\Ati\!}{P}{Q}$ is said to be \defin[successful walk]{successful} if it starts at a vertex in $P$ and ends at a vertex in $Q$.
A word $w \in A^{\!*}$ is said to be \defin[successful word]{successful} if it is the label of some successful walk.
The \defin[recognized language]{language recognized} by  $\att{\Ati\!}{P}{Q}$,
denoted by $\Lang(\att{\Ati\!}{P}{Q})$,
is the set of successful words in $\att{\Ati\!}{P}{Q}$; \ie $\Lang(\att{\Ati\!}{P}{Q}) = \Lang_{P,Q}(\Ati)$. If an automaton $\att{\Ati\!}{P}{Q}$ has a unique initial vertex, this vertex is called the basepoint of $\Ati$, usually denoted by~$\bp$ (\ie $P=\set{\bp}$). An automaton is said to be \defin[pointed automaton]{pointed} if it has a unique common initial and terminal vertex (\ie if~$P=Q=\set{\bp}$); in this case, we write $\Ati_{\!P,Q} = \Ati_{\!\bp}$ or even $\Ati_{\!\bp} = \Ati$ if the basepoint is clear from the context.
\end{defn}

Most of the automata appearing in this paper are pointed and involutive.
The reason for this relevance is their neat link with subgroups, which is stated below, along with its generalization to cosets.

\begin{lem} \label{lem: recognized subgroup}
Let $\Ati$ be an involutive $\Alfi$-digraph, let $\verti, \vertii$ be two vertices of $\Ati$, and let $u$ be the label of a walk from $\verti$ to $\vertii$. Then,
\begin{enumerate}[ind]
\item the free reduction of the language recognized by $\Ati_{\!\verti}$ is a subgroup of $\Free[A]$ called the \defin[recognized subgroup]{subgroup recognized} by $\Ati_{\!\verti}$, and denoted by~$\gen{\Ati_{\!\verti}}$; \ie $\gen{\Ati_{\!\verti}} = \red{\Lang} (\Ati_{\!\verti}) \leqslant \Free[A]$.
If $\Ati$ is pointed, then we write $\gen{\Ati} = \gen{\Ati_{\!\bp}}$.

\item the free reduction of the language recognized by $\Ati_{\!\verti,\vertii}$ is a coset of $\gen{\Ati_{\!\verti}} \leqslant \Free[A]$ called the \defin[recognized coset]{coset recognized} by $\Ati_{\!\verti,\vertii}$, and denoted by~$\gen{\Ati_{\!\verti ,\vertii}}$; more precisely, $\gen{\Ati_{\!\verti,\vertii}} = \red{\Lang} (\Ati_{\!\verti,\vertii}) =\gen{\Ati_{\! \verti}}u=u\gen{\Ati_{\! \vertii}}$.   (Note that then, $\gen{\Ati_{\!\verti}}^{u} = \gen{\Ati_{\!\verti u}} = \gen{\Ati_{\!\vertii}}$.) 
\end{enumerate}
\end{lem}

The notion of recognized subgroup of a pointed $\Alfi$-automaton $\Ati$ is closely related to the topological notion of fundamental group, but they do not coincide in general. The way they are related is crucial for one of our later purposes. As the standard definition in general topology, the \defin{fundamental group} of (the underlying graph of)~$\Ati$ at $\bp$, denoted by $\pi_{\bp}\Ati$, is the set of classes of $\bp$-rounds modulo path homotopy (which, in the case of graphs, is exactly the same as saying modulo backtracking) together with the operation of concatenation of \mbox{$\bp$-rounds}. We emphasize that (the free group) $\pi_{\bp}\Ati$ does not use the labelling but only the underlying graph of $\Ati$.

\begin{prop}
	Let $\Ati$ be an involutive pointed $A$-automaton, and let $\pi_{\bp} \Ati$ be the fundamental group at $\bp$ of the underlying graph of $\Ati$. Then, the map
	\begin{equation} \label{eq: pi -> gen}
		\begin{array}{rcl}
			\widetilde{\ell} \colon \pi_{\bp} \Ati & \to & \Free[A]\\
			\textnormal{[$\walki$]} & \mapsto & \ell(\walki) 
		\end{array}
	\end{equation}
	is a well defined group homomorphism (neither injective nor surjective in general) with image $\gen{\Ati} \leqslant \Free[A]$.
\end{prop}

\begin{defn}
	Let $\Ati$ be an involutive pointed $A$-automaton. The \defin{loss} of~$\Ati$, denoted by $\loss(\Ati)$, is the minimum number of generators of $\ker( \widetilde{\ell})$ \emph{as normal subgroup}.
\end{defn}

\begin{rem} \label{rem: loss det}
	Note that $\loss(\Ati)$ is precisely the loss in rank between the free groups $\pi_{\bp}\Ati$ and $\gen{\Ati}$. Note also that if $\Ati$ is deterministic then $\widetilde{\ell}$ is injective and so $\loss(\Ati)=0$ and $\pi_{\bp}(\Ati) \isom \gen{\Ati}$.
\end{rem}

\begin{defn}
	If a graph (or digraph) $\Ati$ can be obtained by identifying a vertex $\verti$ of a \emph{nontrivial} tree $\Treei$ with a vertex of some graph $\Grii$ disjoint with~$\Treei$, then we say that $\Treei$ is a \defin{hanging tree} of~$\Ati$, and that $\Atii$ is obtained from $\Ati$ after \defin[hanging tree removal]{removing the hanging tree}~$\Treei$.
\end{defn}

\begin{defn}
An involutive automaton $\Ati=\Ati_{\!P,Q}$ is said to be \defin[core automaton]{core} if every vertex
in $\Ati$ appears in some reduced successful walk. The \defin{core} of an involutive automaton~$\Ati$, denoted by $\core(\Ati)$, is the maximal core subautomaton of $\Ati$ containing the distinguished vertices $P$ and $Q$.
\end{defn}

Note that an automaton is core if and only if every hanging tree contains a distinguished vertex. In particular, if $P=\set{\verti}$ then $\core(\Ati)$ is the automaton obtained after taking the connected component of $\Ati$ containing $\verti$, and removing from it all the hanging trees not containing distinguished vertices (if any). It is clear that $\gen{\core (\Ati)} = \gen{\Ati}$. 

As suggested by the previous characterization, the core of a pointed automaton $\Ati$ can still have a hanging tree, which must be a path containing the basepoint as its unique leaf. The maximal such path is called the \defin[automaton tail]{tail} of~$\Ati$.

\begin{defn}
The \defin{restricted core} of an automaton (or digraph) $\Ati$, denoted by $\core^*(\Ati)$, is the labelled digraph obtained after ignoring the distinguished vertices and removing all the hanging trees from $\Ati$. It is clear that $\core^*(\Ati) \subseteq \core{\Ati}$ (after ignoring the distinguished vertices).
\end{defn}
\begin{defn}
A pointed and involutive automaton is called \defin[reduced automaton]{reduced} if it is both deterministic and core.
\end{defn}

\subsection{\bf Homomorphisms between automata}\label{ssec: hom}

\begin{defn}
Let $\Ati = (\Verts, \Edgs, \init,\term,\lab,P,Q)$ and $\Ati'=(\Verts', \Edgs', \init',\term',\lab',P',Q')$ be two $\Alfi$-automata. A \defin[homomorphism of automata]{homomorphism (of automata)} between $\Ati$ and $\Ati'$ is a function $\theta \colon \Verts \to \Verts'$ preserving the automaton structure; that is:
\begin{enumerate}[ind]
\item\label{item: homomorpphism arcs} for every $\verti,\vertii \in \Verts$ and every $\alfi\in \Alfi$, if $\verti \xarc{\alfi\,} \vertii$ then $\verti \theta \xarc{\alfi\,} \vertii \theta$;
\item\label{item: homomorpphism distinguished} for every $\verti \in P$ and every $\vertii \in Q$, $\verti \theta \in P'$ and $\vertii \theta \in Q'$. 
\end{enumerate}
\end{defn}
Note that if $\Ati'$ is deterministic, then condition~\ref{item: homomorpphism arcs} allows us to define a function $\theta_\Edgs: \Edgs \to \Edgs'$ sending every arc $\edgi \equiv \verti \xarc{\alfi\,} \vertii$ to the corresponding arc $\edgi \theta \equiv \verti \theta \xarc{\alfi\,} \vertii \theta$ in $\Ati$. 
We usually abuse language and denote the previous action of the homomorphism on the arcs also by $\theta$, and write $\theta \colon \Ati \to \Ati'$ to denote succinctly the homomorphism of automata. If necessary, we shall denote by $\theta_{\Verts}$ and $\theta_{\Edgs}$ the respective  `restrictions' to vertices and arcs of a homomorphism $\theta = (\theta_\Verts,\theta_\Edgs) \colon \Ati \to \Ati'$.

\begin{rem}
The respective versions of homomorphism for labelled digraphs, digraphs and undirected graphs are obtained after interpreting them as (partial) automata and neglecting the conditions on the missing parts. For example, a homomorphism of (unlabelled) digraphs is a map $\theta \colon \Verts \to \Verts'$ such that for every~$\verti,\vertii \in \Verts$, if $\verti \xarc{\ \,} \vertii$ then $\verti \theta \xarc{\ \,} \vertii \theta$.
\end{rem}

The primordial link between homomorphisms of automata and the corresponding recognized subgroups is summarized in the straightforward statement below.

\begin{lem} \label{lem: sgp iff homom}
Let $\Ati$ and $\Atii$ be a (pointed and involutive) deterministic $\Alfi$-automata. There is a homomorphism of automata $\theta \colon \Ati \to \Atii$ if and only if $\gen{\Ati}\leq \gen{\Atii}\leq \Free[A]$; and in this case, the homomorphism $\theta$ is unique.
\end{lem}

\begin{defn}
We say that a homomorphism $\theta\colon \Ati \to \Ati'$ is \defin{locally injective} (or an \defin{immersion}) if distinct incident arcs always have distinct images by $\theta$; that is if, for every vertex $\verti\in \Verts\Ati$, the restriction of $\theta$ to the neighborhood of $\verti$, say $\theta_{\verti}$, is injective; similarly, we say that $\theta$ is \defin{locally bijective} (or a \defin{cover}) if, for every $\verti\in \Verts\Ati$, $\theta_{\verti}$ is bijective. Also, we say that $\theta$ is an \defin{epimorphism of automata} if both $\theta_\Verts$ and $\theta_\Edgs$ are surjective. 
\end{defn}

\subsection{Stallings automata} \label{ssec: Stallings}
Below, we define the Schreier and the Stallings automaton of a given subgroup ${H\leq \Free}$, two related automata which will be essential along the rest of the paper.

\begin{defn}
	Let $\Free[A]$ be a free group, and let $H$ be a subgroup of $\Free[A]$. Then, the \defin[Schreier automaton]{(right) Schreier automaton} of $H$ \wrt $A$, denoted by $\schreier(H,A)$, is the labelled digraph having $H\backslash \Free[A]$ (the set of right cosets of $\Free[A]$ modulo $H$) as vertex set, an arc $Hw \xarc{a\,} Hwa$ for every coset $Hw\in H\backslash \Free[A]$ and every element $a\in A^{\pm}$, and the coset $H$ as basepoint.
\end{defn}

\begin{rem}
When $H=1$ this is the standard \defin{Cayley digraph} of $\Free$ \wrt $A$; that is, $
\schreier(\Trivial,A) = \cayley(\Free,A) $. More generally, if $H$ is normal in $\Free$ and $\basis=\set{u_i}_i$ is a basis for $\Free$, then the Schreier automaton of $H$ \wrt $B$ is the Cayley digraph of the quotient $\Free/H$; more precisely, $\schreier(H,\set{u_i}_i) = \cayley(\Free/H,\set{H u_i}_i)$.
\end{rem}

The introduced notation for automata applies, \eg we write $\schreier_{\verti,\vertii}(H)$ or $\schreier_{\verti}(H)$  if we want to distinguish particular initial or terminal vertices in $\schreier(H)$.

Note that Schreier automata are involutive, deterministic connected and saturated, but not necessarily core.

\begin{defn}
Let $\Free[A]$ be a free group, let $u \in \Free$, and let $H$ be a subgroup of $\Free$. Then, the \defin{Stallings automaton} of the coset~$Hu$ \wrt $A$ is $\stallings(Hu,A) = \core(\schreier_{H,Hu}(H,A))$. Note that $\gen{\stallings(Hu,A)}=Hu$. 
In particular, the Stallings automaton of $H$ \wrt $A$ is $\stallings(H,A) = \core (\schreier(H,A))$.
Finally, the restricted core of $\stallings(H,A)$ is called the \defin[restricted Stallings]{restricted Stallings automaton} of $H$ \wrt $A$, and denoted by~$\rstallings(H,A)$.
\end{defn}

In both cases, if the basis $A$ is clear from the context (\eg the one given by the presentation), we usually drop the reference to it and write simply $\schreier(H)$, $\stallings(H)$, $\rstallings(H)$, etc. 

\begin{exm}
The (positive part of the) Stallings automaton $\stallings(\Free[A],A)$ has one single vertex and one loop labelled by each $a \in A$. Such an automaton is called a (directed) \defin{bouquet} or \defin{rose} (of $\card \Alfi$ petals), and denoted by $\bouquet_A$ or $\bouquet_n$, where $n = \card \Alfi$.
\begin{figure}[H]
    \centering
    \begin{tikzpicture}[shorten >=1pt, node distance=.5cm and 1.5cm, on grid,auto,>=stealth']
   \node[state,accepting] (1) {};

   \path[->]
        (1) edge[blue,loop right,min distance=12mm,in=330+35,out=330-35]
            node[right] {\scriptsize{$a$}}
            (1)
            edge[red,loop right,min distance=12mm,in=90+35,out=90-35]
            node[above] {\scriptsize{$b$}}
            (1)
            edge[loop right,min distance=12mm,in=210+35,out=210-35]
            node[left] {\scriptsize{$c$}}
            (1);
\end{tikzpicture} 
    \vspace{-15pt}
    \caption{$\bouquet_{3}$, a bouquet of three petals}
    \label{fig: bouquet}
\end{figure}
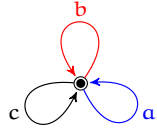
\end{exm}
It is convenient to keep in mind the next remark, which follows immediately from the fact that $\schreier(H)$ is always saturated, $\stallings(H)$ is always core, and $\stallings(H) = \core(\schreier(H))$.
\begin{rem} \label{rem: full}
Let $H$ be a subgroup of a free group $\Free[\Alfi]$. Then, the following statements are equivalent:
\begin{enumerate}[dep,beginpenalty=10000]
    \item $\schreier(H,\Alfi)$ is core.
    \item $\stallings(H,\Alfi)$ is saturated.
    \item $\stallings(H,\Alfi) = \schreier(H,\Alfi)$.
\end{enumerate}
\end{rem}
\medskip

It is obvious from \Cref{lem: recognized subgroup} that there exists a map:
 \begin{equation}\label{eq: Ati ->> <Ati>}
\begin{array}{rcl}
\Set{\text{pointed involutive $A$-automata}} & \to & \Set{\text{subgroups of $\Free[A]$}} \\
\Ati & \mapsto & \gen{\Ati}
\end{array}
 \end{equation}
It is also easy to see that this map is indeed surjective. That is, every subgroup  $H \leqslant \Free[\Alfi]$ is recognized by some pointed involutive $\Alfi$-automaton. Namely, if $S=\set{w_i}_i \subseteq \Free[A]$ is a set of reduced words generating~$H$, then it is enough to consider, for every $w_i = a_{i_1} a_{i_2} \cdots a_{i_p} \in S$ ($a_i \in A^{\pm}$), the (pointed involutive) oriented cycle $\flower(w_i)$ spelling~$w_i$ (or $w_i^{-1}$ if read in the opposite direction),
 \begin{figure}[H]
\centering
\begin{tikzpicture}[shorten >=1pt, node distance=0.2 and 1.5, on grid,auto,>=stealth']
    \node[state, accepting] (0) {};
    \node[state] (a) [above right =  of 0] {};
    \node[state] (1) [above right =  of a] {};
    \node[state] (b) [right =  of 1] {};
    \node[state] (d) [below right =  of 0] {};
    \node[state] (4) [below right =  of d] {};
    \node[state] (c) [right =  of 4] {};

    \path[->]
        (0) edge[]
            node[above] {$\alfi_{i_1}$}
            (a);

     \path[->]
        (a) edge[]
            node[above left] {\scriptsize{$\alfi_{i_2}$}}
            (1);

     \path[->]
        (1) edge[]
            node[above] {$\alfi_{i_3}$}
            (b);

     \path[->,dashed]
        (b) edge[bend left,out=90,in=90,min distance=10mm]
            (c);

    \path[->]
        (c) edge[]
            node[midway,below] {$\alfi_{i_{p-2}}$}
            (4);

    \path[->]
        (4) edge[]
            node[below] {$\alfi_{i_{p-1}}$}
            (d);

    \path[->]
        (d) edge[]
            node[below] {$\alfi_{i_p}$}
            (0);
\end{tikzpicture}
\caption{The petal automaton $\flower(w_i)$} \label{fig: wedge petal}
 \end{figure}
\noindent and then define the \defin{flower automaton} $\flower(S)$, to be the automaton obtained after identifying the basepoints of the petals of the elements in $S$.
\vspace{-10pt}
\begin{figure}[H]
\centering
\begin{tikzpicture}[shorten >=1pt, node distance=2cm and 2cm, on grid,auto,>=stealth',
decoration={snake, segment length=2mm, amplitude=0.5mm,post length=1.5mm}]
  \node[state,accepting] (1) {};
  \path[->,thick]
        (1) edge[loop,out=160,in=200,looseness=8,min distance=25mm,snake it]
            node[left=0.2] {$w_1$}
            (1);
            (1);
  \path[->,thick]
        (1) edge[loop,out=140,in=100,looseness=8,min distance=25mm,snake it]
            node[above=0.15] {$w_2$}
            (1);
  \path[->,thick]
        (1) edge[loop,out=20,in=-20,looseness=8,min distance=25mm,snake it]
            node[right=0.2] {$w_p$}
            (1);
\foreach \n [count=\count from 0] in {1,...,3}{
      \node[dot] (1\n) at ($(1)+(45+\count*15:0.75cm)$) {};}
\end{tikzpicture}
\vspace{-20pt}
\caption{The flower automaton $\flower(w_1,w_2,\ldots,w_p)$}
\label{fig:  flower automaton}
\end{figure}
It is clear that the labels of the $\bp$-rounds of $\flower(S)$ precisely describe the subgroup~$H$; \ie $\gen{\flower(S)} = \gen{S} = H \leqslant \Free[A]$, and hence the map \eqref{eq: Ati ->> <Ati>} is surjective. Observe also that it is far from injective since, for example, different sets of generators for $H$ give different preimages of $H$. Note however that,
$\flower(S)$ is core, and deterministic
\emph{except maybe at the basepoint}.
A key result due to
\citeauthor{stallings_topology_1983}
is that \emph{full} determinism is essentially the only missing condition in order to make this representation unique. We state the language-theoretic version of this result below.

\begin{prop}\label{prop: reduced automata}
Two reduced $\Alfi$-automata $\Ati,\Ati'$ are isomorphic if and only if they recognize the same subgroup; that is, $\Ati \isom \Ati' \Biimp \gen{\Ati} = \gen{\Ati'}$.
\end{prop}

The implication to the right is obvious. And the other one follows easily from \Cref{lem: sgp iff homom}: assuming $\gen{\Ati} = \gen{\Ati'}$ we have the two inclusions and so two morphisms $\theta\colon \Ati \to \Ati'$ and $\theta'\colon \Ati'\to \Ati$; therefore $\theta\theta'$ (\resp $\theta'\theta$) is a homomorphism from $\Ati$ to $\Ati$ (\resp from $\Ati'$ to $\Ati'$) which must be the identity by unicity; hence, $\theta$ and $\theta'$ are inverse of each other and~${\Ati \isom \Ati'}$.

To obtain a reduced representative for a given subgroup $H = \gen{\Ati}$, we define equivalence relations on the set of vertices and arcs of $\Ati$ as follows:
for every $\verti,\vertii \in \Verts \Ati$, we say that $\verti$ and $\vertii$  are \defin[vertex equivalence]{equivalent}, denoted by $\verti \equiv_{\Verts} \vertii$, if and only if there exists a walk from $\verti$ to $\vertii$ recognizing the trivial element in~$\Free[\Alfi]$. Also, for each pair of arcs $\edgi,\edgii \in \Edgs \Ati$, we say that $\edgi$ and $\edgii$ are \defin[arc equivalence]{equivalent}, denoted by $\edgi \equiv_{\Edgs} \edgii$ if and only if they have the same origin, end and label.
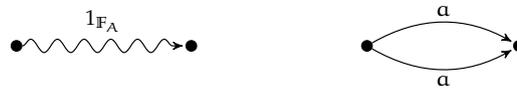
\begin{figure}[H]
\centering
\begin{tikzpicture}[shorten >=1pt, node distance=.5cm and 1.5cm, on grid,auto,>=stealth',anchor=base,baseline]
   \node[state] (1) {};
   \node[state] (2) [right = 2.3 of 1] {};

   \path[->]
        (1) edge[snake it]
            node[pos=0.5,above = 0.1] {$\trivial_{\Free[\Alfi]}$}
            (2);
\end{tikzpicture} 
\hspace{55pt}  
\begin{tikzpicture}[shorten >=1pt, node distance=.5cm and 2cm, on grid,auto,>=stealth',anchor=base,baseline]
   \node[state] (1) {};
   \node[state] (2) [right = of 1] {};

   \path[->]
        (1) edge[bend right]
            node[pos=0.5,below] {$\alfi$}
            (2)
        (1) edge[bend left]
            node[pos=0.5,above] {$\alfi$}
            (2);
   \end{tikzpicture}  
\caption{Equivalent vertices (left) and equivalent arcs (right)}
\end{figure}

It is easy to see that both $\equiv_{\Verts}$ and $\equiv_{\Edgs}$ are indeed equivalence relations compatible with the automaton structure. That is, for every pointed automaton $\Ati$, there exist epimorphisms (of pointed automata) $\smash{\Ati \xonto{\varphi_{\Verts}\,} \Ati / \!\equiv_{\Verts}}$ and
$\smash{\Ati \xonto{\varphi_{\Edgs}\,} \Ati /\!\equiv_{\Edgs}}$ (where equivalent vertices and arcs respectively get identified). If we apply these two epimorphisms successively on a pointed automaton $\Ati$, we obtain $\smash{\varphi \colon \Ati \xonto{\varphi_{\Verts} \varphi_{\Edgs}\,} (\Ati / \!\equiv_{\Verts}) / \!\equiv_{\Edgs}}$.
In this context, the composite function $\varphi = \varphi_{\Verts} \varphi_{\Edgs}$ is called the \defin{(total) folding} of~$\Ati$, and the quotient (folded) automaton is denoted by $\Ati / \varphi$.

Note that $\Ati/\varphi$ is, by construction, a deterministic automaton recognizing the same subgroup as $\Ati$, but not necessarily core. It is not difficult to see that after taking the core of $\Ati / \varphi$ we recover the notion of Stallings automata.

\begin{cor} \label{cor: core folded = Stallings}
Let $\Ati$ be a pointed $\Alfi$-automaton. Then, $\core(\Ati/\varphi) = \stallings(\gen{\Ati},\Alfi)$.
\end{cor}

That is, determinism (followed by core) is enough to obtain a unique representative for the fibers of \eqref{eq: Ati ->> <Ati>}, as we wanted. This is the seminal result of this theory.

\begin{thm}[\citenr{stallings_topology_1983}]\label{thm: Stallings bijection}
Let $\Free[\Alfi]$ be a free group with basis $\Alfi$. Then, 
 \begin{equation}\label{eq: Stallings bijection}
\begin{array}{rcl} \operatorname{St}\colon \set{\,\text{subgroups of } \Free[A] \,} & \leftrightarrow & \set{\,\text{(isom.\ classes of) reduced $\Alfi$-automata}\,} \\ H & \mapsto & \stallings(H,\Alfi) \\ \gen{\Ati} & \mapsfrom & \Ati \end{array}
 \end{equation}
is a bijection. Furthermore, finitely generated subgroups correspond precisely to finite automata and, in this case, the bijection is computable.
\end{thm}

\begin{figure}[H] 
\centering
\begin{tikzpicture}[shorten >=1pt, node distance=1cm and 1.5cm, on grid,>=stealth']
\begin{scope}
   \node[state] (1) {};
   \node[state] (2) [above right = 0.5 and 1.3 of 1] {};
   \node[] (21) [above right = 0.3 and 0.7 of 2] {};
   \node[] (22) [below right = 0.3 and 0.7 of 2] {};
   \node[state] (3) [below right = 0.5 and 1.3 of 1] {};
   \node[] (31) [above right = 0.3 and 0.7 of 3] {};
   \node[] (32) [below right = 0.3 and 0.7 of 3] {};

   \path[->]
        (1) edge[]
            node[pos=0.52,above] {$a$}
            (2)
            edge[]
            node[pos=0.5,below] {$a$}
            (3)
        (2) edge[gray]
            (21)
            edge[gray]
            (22)
        (3) edge[gray!80]
            (31)
            edge[gray!80]
            (32);

    \node[] (i) [right = 2.1 of 1]{};
    \node[] (f) [right = 0.9 of i] {};
    \path[->]
        (i) edge[bend left] (f);

   \node[state] (1') [right= 0.3 of f] {};
   \node[state] (N) [right =  1.3 of 1'] {};
   \node[] (21') [above right = 0.7 and 0.3 of N] {};
   \node[] (22') [above right = 0.3 and 0.7 of N] {};
   \node[] (32') [below right = 0.3 and 0.7 of N] {};
   \node[] (33') [below right = 0.7 and 0.3 of N] {};
   \path[->]
        (1') edge
            node[pos=0.48,above] {$a$}
            (N)
        (N) edge[gray]
            (21')
            edge[gray]
            (22')
            edge[gray!80]
            (32')
            edge[gray!80]
            (33');

\foreach \n [count=\count from 0] in {1,...,3}{
       \node[dot,gray] (2d\n) at ($(2)+(-10+\count*10:0.4cm)$) {};}

\foreach \n [count=\count from 0] in {1,...,3}{
       \node[dot,gray!80] (2d\n) at ($(3)+(-10+\count*10:0.4cm)$) {};}

\foreach \n [count=\count from 0] in {1,...,3}{
       \node[dot,gray] (2d\n) at ($(N)+(34+\count*10:0.4cm)$) {};}

\foreach \n [count=\count from 0] in {1,...,3}{
       \node[dot,gray!80] (2d\n) at ($(N)+(-53+\count*10:0.4cm)$) {};}
  \end{scope}
  
 \begin{scope}[xshift=7.5cm]
  \node[state] (1)  {};
   \node[state] (2) [right = of 1] {};
  \node[] (21) [above right = 0.7 and 0.3 of 2] {};
  \node[] (22) [above right = 0.3 and 0.7 of 2] {};
  \node[] (32) [below right = 0.3 and 0.7 of 2] {};
  \node[] (33) [below right = 0.7 and 0.3 of 2] {};
  \path[->]
        (2) edge[gray]
            (22)
            edge[gray]
            (32);

  \path[->]
        (1) edge[bend left]
            node[pos=0.5,above] {$a$}
            (2);
  \path[->]
        (1) edge[bend right]
            node[pos=0.5,below] {$a$}
            (2);

    \node[] (i) [right = 2.3 of 1] {};
    \node[] (f) [right = 0.9 of i] {};
    \path[->]
        (i) edge[bend left] (f);

  \node[state] (1') [right= 0.3 of f] {};
  \node[state] (N) [right =  1.3 of 1'] {};
  \node[] (21') [above right = 0.7 and 0.3 of N] {};
  \node[] (22') [above right = 0.3 and 0.7 of N] {};
  \node[] (32') [below right = 0.3 and 0.7 of N] {};
  \node[] (33') [below right = 0.7 and 0.3 of N] {};
  \path[->]
        (1') edge
            node[pos=0.48,above] {$a$}
            (N)
        (N) edge[gray]
            (22')
            edge[gray]
            (32');

\foreach \n [count=\count from 0] in {1,...,3}{
      \node[dot,gray] (2d\n) at ($(2)+(-12+\count*10:0.4cm)$) {};}

\foreach \n [count=\count from 0] in {1,...,3}{
      \node[dot,gray] (2d\n) at ($(N)+(-12+\count*10:0.4cm)$) {};} 
 \end{scope}
\end{tikzpicture}
 \vspace{-10pt}
\caption{An open folding (left) and a closed folding (right)}
\label{fig: foldings}
\end{figure}
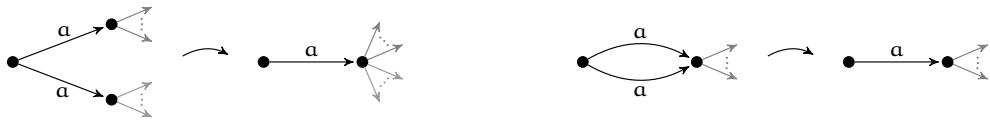

The computability of the finitely generated case follows from the finiteness of the flower automaton when the generating set $S$ for $H$ is finite. Roughly speaking, the algorithm consists of drawing the flower automaton $\flower(S)$ and then computing the folded automaton $\flower(S)/\varphi$ by successively applying a \emph{finite} sequence of local (binary) versions of the identifications $\equiv_{\Verts}$ and $\equiv_{\Edgs}$ defined above. These are traditionally encapsulated in the two kind of transformations --- called \defin[elementary foldings]{elementary} or \defin{Stallings foldings} --- described below (see \Cref{fig: foldings}).

\begin{description}
  \item [\mbox{\defin[open foldings]{Open foldings:}}] identifications of two nonparallel arcs with the same origin and labels (and its corresponding inverses).

  \item [\mbox{\defin[closed foldings]{Closed foldings:}}] \label{item: closed foldings}
  identifications of two parallel arcs with the same label (and its corresponding inverses).
\end{description}

Note that closed foldings are nothing more than partial (binary) versions of the arc identification~$\varphi_{\Edgs}$, whereas open foldings consist of a binary version of the vertex identification $\varphi_{\Verts}$ followed by the corresponding closed folding. More precisely, the algorithm to compute $\stallings(H)$ given a finite generating set $S$ for $H$ consists of building $\flower(S)$ and then successively detecting and performing possible (open or closed) foldings until no more foldings are available (so reaching the folded automaton $\flower(S)/\varphi$). It is not difficult to see that if the words in $S$ are already reduced then none of the folding transformations breaks the coreness of $\flower(S)$; therefore, the obtained automaton $\flower(S)/\varphi$ is already reduced, and hence isomorphic to $\stallings(H)$. Such a procedure is called a \defin[folding sequence]{Stallings folding sequence} for the subgroup $H$.
 \begin{figure}[H]
\begin{equation} \label{eq: folding process}
\flower(S) = \Ati
\xtransf{\!\!\!\! \varphi^{(1)} \!\!\!}\Ati^{^{(1)}}
\xtransf{\!\!\!\! \varphi^{(2)} \!\!\!} \cdots
\xtransf{\!\!\!\! \varphi^{(p)} \!\!\!} \Ati^{^{(p)}}
\! = \stallings(H)
\end{equation}
\vspace{-15pt}
\caption{A (finite) Stallings folding sequence computing $\stallings(H)$}
 \end{figure}
In particular, both kinds of Stallings foldings preserve the recognized subgroup; and, more importantly,
the final object being $\stallings(H)$ means that the result of the Stallings folding sequence depends  neither  on  the  chosen  sequence  of  foldings  nor  on  the  starting (finite) generating set~$S$ for~$H$. Note also that
\begin{enumerate*}[ind]
\item\label{item: lost arc} exactly one arc is lost in each Stallings folding and
\item\label{item: loss} the loss of a closed (\resp open) Stallings folding is~$1$ (\resp $0$).
\end{enumerate*}
Since the number of arcs at the beginning is finite, \ref{item: lost arc} guarantees reaching the final reduced automaton $\stallings(H)$ in finite time, and from \ref{item: loss} one sees that, \emph{if $\Ati$ is finite}, then the loss of the total folding $\smash{\Ati \xtransf{\varphi} \Ati/\varphi}$ is the total (finite) number of closed foldings in the (any) folding sequence.

Remarkably enough, it has been proven that this Stallings folding procedure can be performed in almost linear time (see~\parencite{touikan_fast_2006}), a result emphasizing the algorithmic friendliness of the process. 

Conversely, for the computation of the inverse map in \eqref{eq: Stallings bijection}, a well known topological argument applies.
It can be summarized as follows.

\begin{nott}
If $\verti$ and $\vertii$ are vertices in a tree $T \leqslant \Ati$,
then we denote  by $\verti \xwalk{\scriptscriptstyle{T}} \vertii$ the unique reduced walk in $T$ from $\verti$ to $\vertii$.
\end{nott}

\begin{prop}
\label{prop: generators from tree}
Let $\Ati$ be a connected pointed $\Alfi$-automaton and let $T$ be a spanning tree of $\Ati$. Then
the set
\begin{equation*} \label{eq: generating set}
S_T
\,=\,
\big{\{}
\,
\red{\lab}(\bp \xwalk{\scriptscriptstyle{T}} \! \bullet \! \arc{\edgi} \!\bullet\! \xwalk{\scriptscriptstyle{T}} \bp)\st \edgi\in \Edgs^+\Ati \setmin \Edgs T \,
\big{\}}
\,\subseteq\, \Free[A]
\end{equation*}
is a generating set for
$\gen{\Ati}$. Furthermore, if $\Ati$ is deterministic then
$S_T$ is a (Nielsen reduced) free basis
for~$\gen{\Ati}$, and hence $\rk(\gen{\Ati}) = \rk(\Ati)$. 
Furthermore, if $\Ati$ is finite then $S_T$ is computable and $\rk(\gen{\Ati}) = 1 - \card \Verts(\Ati) + \card \Edgs^{+}(\Ati)$.
\end{prop}

An immediate (but relevant) consequence of \Cref{prop: generators from tree} is stated below.
\begin{cor}\label{cor: incl => ff}
Let $\Ati$ and $\Atii$ be pointed, deterministic and involutive $\Alfi$-automata. If $\Ati$ is a subautomaton of~$\Atii$, then ${\gen{\Ati} \leqslant\ff \gen{\Atii}}$.
\end{cor}

\begin{rem}
The converse of  \Cref{cor: incl => ff} is far from true, since the number of subatomata of a finite automaton is finite, whereas the number of free factors of a non-cyclic subgroup is not.
\end{rem}

For full proofs and details on the results in this section, we refer the reader to the (topologically oriented) original work~\parencite{stallings_topology_1983} by Stallings, and to the more combinatorial and language-theoretic approaches~\parencite{kapovich_stallings_2002,bartholdi_rational_2010} (closer to ours).

\section{\bf List of applications} \label{sec: list}
In this section we summarize a selection of applications of Stallings theory, many of them of algorithmic nature. We recall that the seminal Dehn problems --- namely, the \emph{Word Problem}~(\WP), the \emph{Conjugacy Problem}~(\CP), and the \emph{Isomorphism Problem}~(\IP) --- are all easily decidable for finitely generated free groups (see \eg~\parencite{lyndon_combinatorial_2001}).
\subsection{\bf First applications}

In this subsection we survey several basic results about free groups and their lattice of subgroups, which can be easily derived from Stallings theory. Many of them are classical, with original proofs more combinatorially oriented, and more technical and involved than the modern ones using Stallings machinery.

\subsubsection{\bf Structure of subgroups and Nielsen--Schreier theorem}

The fundamental theorem below, was initially proved by Nielsen for finitely generated subgroups, and later extended by Schreier to full generality. The original proofs (using somewhat elaborated combinatorial arguments) become transparent when translated into automata theoretic language.

\begin{thm}[Nielsen--Schreier] \label{thm: Nielsen-Schreier}
Subgroups of free groups are again free.
\end{thm}

As we have seen in \Cref{ssec: Stallings}, if $H$ is a subgroup of a free group, then $H$ is recognized by $\Ati = \stallings(H)$.  Since $\Ati$ is deterministic, we know (\Cref{rem: loss det}) that $H=\gen{\Ati}$ is isomorphic (as group) to the fundamental group of its (undirected) underlying graph; \ie $H = \gen{\Ati} \isom \pi_{\bp}(\Gamma)$. By \Cref{thm: free geom}, the claimed result follows.

Let us point out that the Stallings bijection \eqref{eq: Stallings bijection} allows us not only to guarantee that every subgroup of a free group is again free (\Cref{thm: Nielsen-Schreier}) but also to immediately describe which (isomorphic classes of) free groups embed in a given finitely generated free group: since Stallings automata of any countable rank are possible over an alphabet of two letters, the description below follows.

\begin{cor} \label{cor: ranks of subgroups}
Let $\Free[\kappa]$ be a free group of rank $\kappa \in [2,\aleph_0]$. Then, for every cardinal $\mu \in [0,\aleph_0]$ there exists a subgroup of $\Free[\kappa]$ isomorphic to $\Free[\mu]$.
\end{cor}

We provide two classical examples of infinite rank subgroups of $\Free_2$ through their corresponding Stallings automata.

\begin{exm}
Let $\Free_{\set{a,b}} = \pres{a,b}{-}$ be a free group of rank $2$. Then, the normal closure of $b$, denoted by $\normalcl{b}$, is a subgroup of $\Free_{\set{a,b}}$ of infinite rank, with (infinite rank) Stallings automaton $\stallings(\normalcl{b})$ depicted in \Cref{fig: Stallings <<b>>}: 
\begin{figure}[H] \label{fig: Stallings <<b>>}
\centering
  \begin{tikzpicture}[shorten >=1pt, node distance=1.2 and 1.2, on grid,auto,>=stealth']
   \node[state,accepting] (0) {};
   \node[state] (1) [right = of 0]{};
   \node[state] (2) [right = of 1]{};
   \node[state] (3) [right = of 2]{};
   \node[] (4) [right = 1.5 of 3]{$\cdots$};
   \node[] (c) [right = .5 of 4]{;};

   \node[state] (-1) [left = of 0]{};
   \node[state] (-2) [left = of -1]{};
   \node[state] (-3) [left = of -2]{};
   \node[] (-4) [left = 1.5 of -3]{$\cdots$};

   \path[->]
        (0) edge[loop above,red,min distance=10mm,in=55,out=125]
            node[] {\scriptsize{$b$}}
            (0)
            edge[blue]
            node[below] {\scriptsize{$a$}}
            (1);

    \path[->]
        (1) edge[loop above,red,min distance=10mm,in=55,out=125]
            (1)
            edge[blue]
            (2);

    \path[->]
        (2) edge[loop above,red,min distance=10mm,in=55,out=125]
            (2)
            edge[blue]
            (3);

    \path[->]
        (3) edge[loop above,red,min distance=10mm,in=55,out=125]
            (3)
            edge[blue]
            (4);

    \path[->]
        (-1) edge[loop above,red,min distance=10mm,in=55,out=125]
            (-1)
            edge[blue]
            (0);

    \path[->]
        (-2) edge[loop above,red,min distance=10mm,in=55,out=125]
            (-2)
            edge[blue]
            (-1);

    \path[->]
        (-3) edge[loop above,red,min distance=10mm,in=55,out=125]
            (-3)
            edge[blue]
            (-2);

    \path[->]
        (-4) edge[blue]
            (-3);
\end{tikzpicture}
\caption{The Stallings automaton for $\normalcl{b} \leqslant \Free_{\set{a,b}}$}
\end{figure}

Also, the commutator subgroup
$H=\Comm{\Free_{\set{a,b}}}$ has infinite rank and its Stallings automaton $\stallings(H)$ is the Cayley digraph of $\Free_{\set{a,b}} /H\simeq \mathbb{Z}^2$ \wrt the canonical free abelian basis, \ie the standard infinite bidimensional grid (with horizontal edges labelled, say $a$, and vertical edges labelled $b$). 
\end{exm}

\subsubsection{\bf Basis and rank}
Many well known fundamental facts about free bases can be easily derived from Stallings construction. The lemma below is the starting point for most of them.

\begin{lem} \label{lem: free iff loss=0}
Let $\Free$ be a free group, and let $S \subseteq \Free$. Then,
$S$ is freely independent if and only if $\loss(\flower(S)) = 0$.
\end{lem}
Hence, a subset $S$ of a free group $\Free[A]$ is a basis of $\Free[A]$ if and only if $\loss(\flower(S)) = 0$ and $\stallings(\gen{S}) = \bouquet_{A}$. An alternative proof for the essential result below follows.

\begin{thm} \label{thm: isom iff card}
All the bases of a free group $\Free$ have the same cardinality, which is precisely the rank of~$\Free$. Hence, $\Free[A] \isom \Free[B]$ if and only if $\card A = \card B$.
\end{thm}

A geometric description (and computability, in the finitely generated case) of bases and ranks of subgroups of the free group is given in \Cref{prop: generators from tree}. The computability of the Stallings automaton of finitely generated subgroups immediately provides that of a basis and the corresponding rank.

\begin{cor}
If $S$ is a finite subset of a free group $\Free$, then a basis (and hence the rank) of $\gen{S}$ is computable. Concretely, ${\rk(\gen{S})=1-\card \Verts\Ati+\card \Edgs^+\Ati}$, where $\Ati =\stallings(\gen{S})$. 
\end{cor}

According to \Cref{thm: isom iff card}, bases of $\Free$ are generating sets of minimum cardinality (namely, the rank of $\Free$). The  converse also follows, as far as this cardinality is finite: indeed, if $S \subseteq \Fn$ has cardinal~$n$ and generates~$\Fn$, then
$%
    \rk (\flower(S))
    =
    n
    =
    \rk (\Fn)
    =
    \rk (\stallings(\gen{S}))
$. %
Therefore, $\loss(S) = 0$ and, from \Cref{lem: free iff loss=0}, $S$ is already a free basis for~$\Fn$.

\begin{prop} \label{prop: basis iff min gs}
A subset $S$ of a  free group $\Fn$ (of finite rank $n$) is  a basis for $\Fn$ if and only if $\gen{S} = \Fn$ and $\card S = n$. That is, finitely generated free groups are Hopf.
\end{prop}

Recall that a group $G$ is said to be \defin[Hopf group]{Hopf} (or \defin[Hopfian group]{Hopfian}) if every epimorphism $G \to G$ is indeed injective (and hence an isomorphism).

Note that the previous result mimics to certain extent the situation in linear algebra, were a (minimal) generating set with as many elements as the ambient dimension is already a basis. It is worth noting, however, that in free groups \emph{it is not true} that a minimal generating set of $\Fn$ is necessarily a basis of~$\Fn$.

\begin{lem}
In a nontrivial finitely generated free group of any rank, there exist minimal generating sets of arbitrarily big finite cardinality. 
\end{lem}

To see this fact in the free group of rank 1, namely $\mathbb{Z}$, we just need to use an elementary arithmetic property: take $n$ pairwise different prime numbers $p_1, \ldots ,p_n$ and consider the set $S=\{p_1\cdots p_{i-1} p_{i+1} \cdots p_n \st \allowbreak i=1,\ldots ,n\}$; clearly $\gen{S}=\mathbb{Z}$, whereas no strict subset of $S$ generates the full group~$\mathbb{Z}$. Graphically, this corresponds to a flower automaton with $n$ petals such that, after a sequence of foldings, collapses down to the bouquet, but not if we remove any single petal. Of course, one can construct many other such examples for non-cyclic free groups  based  in  geometric properties other than just arithmetic.

The \defin[subgroup corank]{corank} of a subgroup $H$ of  group $G$, denoted by $\operatorname{crk}(H)$, is the minimum number of elements that one needs to add to $H$ in order to generate the full ambient group~$G$. 

\begin{thm}[{\citenr{puder_primitive_2014}; \citenr{delgado_lattice_2020}}]\label{thm: join-corank is computable}
There exists an algorithm that, given a finite subset $S$ of $\Fn$, outputs the corank of $\gen{S}$ in~$\Fn$.
\end{thm}

Not surprisingly, one can also use Stallings automata to prove the previous result. For example, the proof in~\parencite{delgado_lattice_2020} is by induction, based on the fact that 
$\operatorname{crk}(\gen{\Ati}) \leq r$
if and only if
$\operatorname{crk}( \gen{\Ati_{\!\!/\verti \shorteq \vertii}}) \leq r-1$, where $\Ati_{\!\!/\verti \shorteq \vertii}$ denotes the automaton obtained from $\Ati$ after identifying  some pair of distinct vertices~$\verti,\vertii$.

\subsubsection{\bf Subgroup membership problem}

The (uniform) subgroup membership problem, often called simply \emph{membership problem} is formally stated below for a finitely presented group $G = \pres{A}{R}$.

\begin{named}[(Subgroup) Membership Problem, $\MP(G)$]
Decide\margin{Membership Problem (\MP)}, given a finite set of words $w,u_1,\ldots,u_k$ in $A^{\pm}$, whether $w$ represents an element in the subgroup $\gen{u_1,\ldots,u_k}_G$.
\end{named}

If $S \subseteq \Free$, it is clear by construction that the reduced labels of $\bp$-rounds in the flower automaton $\flower(S)$ describe exactly the elements of the subgroup $H =\gen{S}$ generated by $S$. Moreover, since the folding operation does not change the recognized elements, it is immediate that the reduced labels of the $\bp$-rounds in  $\stallings(H)$ are exactly the same as those in $\flower(S)$.

\begin{prop}
Let $H$ be a subgroup of a free group $\Free$ and let $w \in \Free$. Then, $w\in H$ if and only if $w=\rlab(\walki)$, where $\walki$ is a \bp-round in $\stallings(H)$. 
\end{prop}

Now, given $S=\{u_1, \ldots ,u_k\}$, compute $\stallings(S)$; since it is deterministic, the word $w$ can be read as the label of \emph{at most} one reduced $\bp$-round in $\stallings(H)$, and we can easily check whether this is the case or not by inspection. If not, then $w\not\in \gen{u_1,\ldots ,u_k}$; otherwise, $w\in \gen{u_1,\ldots ,u_k}$. 

\begin{thm}\label{thm: MP}
The subgroup membership problem $\MP(\Free)$ is decidable.
\end{thm}

Moreover, whenever the answer of $\MP(\Free)$ to an input $(w,S)$ is \yep, one can further obtain an expression of $w$ as a word in the original generators $S$ by a brute force algorithm: start a diagonal procedure enumerating all the formal words in $\gen{S}$ and successively comparing them with $w$ (using the decidabililty of the word problem) until reaching a granted match. However, this approach can be highly improved using the Stallings machinery: when the answer to \MP\ is \yep, we have $w$ realized as (the label of) a reduced $\bp$-round in $\stallings(H)$; lifting it back along the sequence of foldings, we obtain a reduced $\bp$-round in $\flower(S)$ whose label also equals (\ie reduces to) $w$: this is precisely an expression of $w$ as a word in the original generators~$\{u_1, \ldots ,u_k\}$. 

In this argument we have used that, along each elementary folding $\smash{\Ati \! \xtransf{\!\!\!} \!\Ati'}$, any reduced $\bp$-round in $\Ati'$ lifts to a reduced $\bp$-round in $\Ati$ spelling the same label (in $\Free[A]$). Note further that this lifting is \emph{unique} whenever the folding is open, and not unique (in fact, there are infinitely many) when the folding is closed. This corresponds to the fact that the original set of generators $S=\{u_1, \ldots ,u_k\}$ is a basis of $\gen{u_1, \ldots ,u_k}$ if and only if there is no loss in the sequence of foldings, \ie there is no closed folding: in this case, of course, $w$ can be expressed as a word in the $\{u_1, \ldots ,u_k\}$ in a unique possible way, while in the opposite case there are non-trivial relations among them and so the expression of $w$ is not unique. 

Moreover, we can obtain a full set of such relations, \ie a finite presentation of the (free) group $\gen{u_1, \ldots ,u_k}$ with the (non-free) generators $S=\{u_1, \ldots ,u_k\}$ by following a similar procedure (note that the number of relations must be precisely $\loss(\flower(S))$, \ie the difference between the number of original generators, $k$, and the real rank of $\gen{u_1, \ldots ,u_k}$): at each closed folding $\Ati \! \smash{\xtransf{\!\!\!}} \!\Ati'$ in the sequence, take the obvious non-trivial reduced $\bp$-round in $\Ati$ with trivial label (choose an arbitrary walk $\gamma$ from~$\bp$ to the origin of the parallel edges to be folded, cross one of them, return through the other, and then go back to $\bp$ along $\gamma^{-1}$) and lift it up to $\flower(S)$, obtaining in this way a word on $\{u_1, \ldots ,u_k\}$ with trivial label, \ie a relation in $\gen{u_1, \ldots ,u_k}$. It can be seen that, repeating this operation with each closed folding, one obtains a complete set of relations, \ie a presentation of  the (free) group $\gen{u_1, \ldots ,u_k}$ with the (non-free) generators~$S=\{u_1, \ldots ,u_k\}$.

\subsubsection{\bf Conjugation, normality, and normalizers}

Stallings theory encapsulates subgroup conjugation in the result below, which is clear from \Cref{lem: recognized subgroup}.

\begin{prop} \label{prop: St conj sgp}
Let $\Free$ be a free group, let $H \leqslant \Free$, and let $w \in \Free$. Then, $\stallings(H^{w}) \,=\, \core (\schreier_{Hw}(H))$.
\end{prop}

\begin{thm}
Let $H,K \leqslant \Free$. Then, $H$ and $K$ are conjugate if and only if their restricted Stallings digraphs coincide, \ie if and only if $\rstallings(H) = \rstallings(K)$. 
\end{thm}

Clearly, we can always algorithmically check this condition for finitely generated subgroups. The result below follows.

\begin{prop}
The (uniform) subgroup conjugacy problem $\SCP(\Free)$ is solvable.
\end{prop}

Recall that a subgroup $H$ of $G$ is normal (in $G$) if for every element $g \in G$, $H^g =H$. A graphical characterization of normality within free groups follows neatly from \Cref{prop: St conj sgp}.

\begin{prop} \label{prop: normal iff}
Let $H$ be a nontrivial subgroup of a free group $\Free$. Then the following statements are equivalent:
\begin{enumerate}[dep]
\item H is normal in $G$;
\item $\schreier(H)$ is vertex transitive;
\item $\schreier(H)$ is vertex transitive and core;
\item $\stallings(H)$ is vertex transitive and saturated.
\end{enumerate}
\end{prop}

Recall that if $H$ is a subgroup of a group $G$, then the \defin[subgroup normalizer]{normalizer of} $H$ (in~$G$) is $N_G(H)=\set{g\in G \st H^g=H}\leqslant G$; that is, the biggest subgroup of $G$ where $H$ is normal, $H\normaleq N_G(H)\leqslant G$.
Using the graphical idea behind~\Cref{prop: normal iff}, it is easy to see the following result for free groups. 

\begin{prop}
Let $H$ be a nontrivial subgroup of a free group $\Free$. If $H$ is finitely generated then so is its normalizer $N_{\Free}(H)$, with a basis being computable from any given finite set of generators for $H$. In particular, $N\normaleq\fin N_{\Free}(H)\leq \Free$.  
\end{prop}

To see this, compute the Stallings automaton $\stallings(H)$ for $H$, the length ${\ell\geq 0}$ of the tail going from $\bp$ to $\rstallings(H)$, and denote by $\estallings(H)$ the \mbox{$\ell$-neighborhood} of $\rstallings(H)$ within $\schreier(H)$. Now, since $\estallings(H)$ is still finite, we can list the vertices of $\estallings(H)$ which are images of $\bp$ by labelled-digraph automorphisms of $\estallings(H)$: if $\verti$ is such a vertex, then the label $w$ of any walk from $\bp$ to $\verti$ satisfies $H^w=H$ (and any such $w\in \Free$ is covered in this way). Therefore, identifying all those vertices to $\bp$ and folding, we obtain the Stallings automaton for the normalizer of $H$ in $\Free$; in particular, $N_{\Free}(H)$ is again finitely generated, and a basis for it is computable. Moreover, since $\trivial \neq H$ is normal in $N_{\Free}(H)$ by construction, and both are free and finitely generated, we deduce that ${N\normaleq\fin N_{\Free}(H)\leq \Free}$.

\subsubsection{\bf Whitehead's cut vertex lemma}

A classical result from the 1930's by Whitehead~\parencite{whitehead_certain_1936} has been of central historic importance to investigate the properties of free groups: the so-called \emph{cut vertex lemma}. As happened with many other classical results, there is a modern much simpler proof of it using Stallings automata. We survey it in this subsection. 

An element $w\in \Free[n]$ is said to be \defin[primitive element]{primitive} if it belongs to some basis for $\Free[n]$; more generally, a set of elements $S\subseteq \Free[n]$ is called \defin[separable subset]{separable} if there exists a non-trivial decomposition $\Free[n] = H * K$ such that any $s\in S$ has a conjugate either in $H$ or in $K$; of course, if $w$ is primitive then $\{w\}$ is separable (in general, the other implication is not true).

Whitehead defined the so-called (undirected) \defin{Whitehead graph}, $\Wh_A(S)$, of a set $S$ of cyclically reduced words in $\Free[A]$ as follows: the set of vertices is $A^{\pm 1}$ and, for any consecutive letters $x^{\epsilon}y^{\delta}$ in the reduced expression of each $s\in S$ considered cyclically, there is an undirected edge between $x^{\epsilon}$ and $y^{-\delta}$. So, $\Wh_A(w)$ has $2n$ vertices and $\sum_{s\in S} |s|$ edges, including possible parallel edges (and no loops because each $s \in S$ is assumed to be cyclically reduced). For $S=\{w\}$, we write $\Wh_A(w)=\Wh_A(\{w\})$.  

The classical Whitehead cut vertex lemma and a later improvement due to Stallings are stated below. 

\begin{thm}[\citenr{whitehead_certain_1936}]
Let $w\in \Free[A]$. If $w$ is primitive then $\Wh_A(w)$ is either disconnected or has a cut vertex.
\end{thm}
\begin{thm}[\citenr{stallings_whitehead_1999}]
Let $S\subseteq \Free[A]$. If $S$ is separable then $\Wh_A(S)$ is either disconnected or has a cut vertex.
\end{thm}

Here, a \defin{cut vertex} is a vertex whose removal disconnects the connected component where it is contained. From very early in the research history of free groups, this result provided an easy-to-check necessary condition for an element $w \in \Free[A]$ to be primitive; note that the converse is not true, in general. The original proof was not elementary and made strong use of \emph{handlebodies}, a special kind of 3-dimensional varieties useful to investigate free groups. The proof of Stallings extension was along the same lines. 

\begin{figure}[H]
\centering
\begin{tikzpicture}[shorten >=1pt, node distance=2cm and 2cm, on grid,auto,>=stealth',
decoration={snake, segment length=2mm, amplitude=0.5mm,post length=1.5mm}]
  \node[state] (1) {};
  \node[state] (2) [right = of 1]{};
            (1);
  \path[->,thick]
        (1) edge[blue,loop,out=110,in=160,looseness=8,min distance=15mm]
            node[above left] {$a_1$}
            (1);
  \path[->,thick]
        (1) edge[loop,out=200,in=250,looseness=8,min distance=15mm]
            node[below left] {$a_k$}
            (1);
    
    \path[->,thick]
        (1) edge[blue,bend left = 45] node[above] {$a_{1}$} (2)
        (1) edge node[above] {$a_{k+1}$} (2)
        (1) edge[bend right = 55] node[below] {$a_{\ell}$} node[above] {$\vdots$} (2);

    \path[->,thick]
        (2) edge[loop,out=20,in=70,looseness=8,min distance=15mm]
            node[above right] {$a_{\ell+1}$}
            (2);
  \path[->,thick]
        (2) edge[loop,out=-70,in=-20,looseness=8,min distance=15mm]
            node[below right] {$a_{n}$}
            (2);
\foreach \n [count=\count from 0] in {1,...,3}{
      \node[dot] (1\n) at ($(1)+(167 +\count*15:0.75cm)$) {};}
      
      \foreach \n [count=\count from 0] in {1,...,3}{
      \node[dot] (1\n) at ($(2)+(-13 +\count*15:0.75cm)$) {};}
\end{tikzpicture}
\vspace{-20pt}
\caption{An almost-rose}
\label{fig: almost-rose}
\vspace{-5pt}
\end{figure}
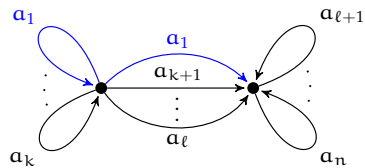

In 2019, M. Heusener and R. Wiedmann found a very elegant and short proof for these results using Stallings automata (see~\parencite{heusener_remark_2019}). The main idea is to generalize the notion of Whitehead graph, and to consider $A$-almost-roses. They define the \defin[automaton Whitehead graph]{Whitehead graph} of an $A$-automaton $\Ati$, denoted by $\Wh_A(\Ati)$, as follows: the vertex set is $A^{\pm 1}$ and, for any label $x^{\epsilon}y^{\delta}$ of a length two reduced walk in $\Ati$, we put an undirected edge between $x^{\epsilon}$ and $y^{-\delta}$. Note that, for $S\subseteq \Free[A]$, $\Wh_A(S)=\Wh_A(\Ati_{\!S})$, where $\Ati_{\!S}$ is the disjoint union of directed cycles labelled by the cyclic reductions of the $s\in S$. And they define an \defin{almost-rose} as a Stallings automaton which folds into the $A$-bouquet with a \emph{single open folding}. It is straighforward to see that an almost-rose has 2 vertices, $\card A+1$ edges and, up to permutation and inversion of labels, looks like \Cref{fig: almost-rose}.
\vspace{-10pt}
They see that the Whitehead graph of an almost-rose always has a cut vertex, and then prove a lemma saying that, for any separable set $S\subseteq \Free[A]$, there exists an almost-rose $\bm{\Theta}$ where the cyclic reduction of each $s\in S$ is readable as a closed walk; in particular, $\Wh_A(S)\subseteq \Wh_A(\bm{\Theta})$, both having the same set of vertices. Then, $\Wh_A(\bm{\Theta})$ has a cut vertex and so, either $\Wh_A(S)$ is disconnected or it has a cut vertex too. 

\subsection{\bf Cosets and index}

From the algorithmic point of view, maybe the first natural question on the index of subgroups is about its finiteness. As usual, we state it for a group $G$ given by a finite presentation $\pres{A}{R}$. 

\begin{named}[Finite Index Problem, $\FIP(G)$]
Decide\margin{Finite Index Problem (\FIP)}, given a finite subset $S \subseteq (A^{\pm})^*$, whether the subgroup $\gen{S}_{G}$ has finite index in $G$, and, if so, compute a full transversal (and hence the index).
\end{named}

The well-known Schreier's Lemma (see \eg \cite[\textnormal{1.6.11}]{robinson_course_1996}) implies that finite index subgroups of finitely generated groups are finitely generated as well. Let's see how the Stallings machinery allows us to recover (and tighten) this result in the realm of free groups. It is enough to recall that, by definition, the index of a subgroup $H\leqslant \Free$ is the number of vertices in $\schreier(H)$ so, $H$ is of finite index in $\Free$ if and only if $\card \Verts \schreier(H)<\infty$. Since $\stallings(H) = \core(\schreier(H))$, a characterization in terms of the Stallings automaton follows. 

\begin{prop} \label{prop: fi iff}
Let $H$ be a subgroup of a free group $\Free$ (of arbitrary rank). Then, $H$ has finite index in $\Free$ if and only if $\stallings(H)$ is saturated and $\card\Verts\stallings(H)<\infty$.
\end{prop}

Note that, in this case, the index $\ind{H}{\Free}$ is the number (cardinal) of vertices in $\stallings(H) = \schreier(H)$. In particular, if $H$ is a subgroup of finite index \emph{in a finitely generated ambient} $\Free =\Fn$ (of rank, say, $n$), then the index ${\ind{H}{\Fn}}$ is the number of vertices in the saturated \emph{$2n$-regular} automaton $\Ati =\stallings(H)$, which can be easily related to the rank of $\Ati$ (and hence of $H$) using elementary graph-theoretical arguments.

\begin{cor}[\emph{Schreier index formula}]
If\margin{Schreier index formula} $H$ is a subgroup of $\Free$ of finite index, then
\begin{equation}
    \rk(H) - 1 \,=\, \ind{H}{\Free}\cdot (\rank(\Free)-1) \,. 
\end{equation}
\end{cor}

Note that the Schreier index formula is still valid when $\Free =\Free[\,\aleph_0]$ (corresponding to the fact that a finite index subgroup of a group of infinite rank is always of infinite rank as well).

On the other side, if both the ambient group $\Free$ and the subgroup $H$ are finitely generated, then the Stallings automaton $\stallings(H)$ is computable, and the condition in \Cref{prop: fi iff} is immediately checkable. Finally, if $\stallings(H)$ is finite, it is enough to take an spanning tree $T$ of $\stallings(H)$, and it is clear that the set
\begin{equation}
	\big\{\, \rlab(\bp \xwalk{T} \verti) \st \verti \in \Verts\, \stallings(H)\,\big\}
\end{equation}
of reduced labels of the $T$-walks from the basepoint to all the vertices in $\stallings(H)$ constitutes a transversal of $H$ in $\Fn$.

\begin{cor} \label{cor: FIP}
The finite index problem $\FIP(\Fn)$ is algorithmically decidable.
\end{cor}

Note also that, combining the graphical characterizations of finite index (\Cref{prop: fi iff}) and normality (\Cref{prop: normal iff}), we reach a sort of reciprocal of Schreier's Lemma for nontrivial normal subgroups of the free group.

\begin{cor}
A nontrivial normal subgroup of a finitely generated free group $\Fn$ is finitely generated if and only if it has finite index.
\end{cor}

\begin{cor}
A nontrivial normal subgroup of $\Free[\aleph_0]$ has always infinite rank.
\end{cor}

\medskip

Furthermore, the characterization in \Cref{prop: fi iff} can be combined with \Cref{lem: =deficient} to prove the existence of finite index subgroups satisfying certain relevant properties: the key idea is that, in any finite Stallings automaton $\Ati$, for every $a\in A$, we can pair up the $a$-deficient vertices with the $a^{-1}$-deficient vertices and, joining them with new arcs (and no more vertices) we can `complete' $\Ati$ to a finite saturated automaton. Together with \Cref{prop: fi iff}, this proves a sequence of classical interrelated theorems proved by M. Hall Jr.\ in the middle of the XX century.

\begin{lem}[\citenr{hall_jr_coset_1949}] \label{lem: F rf}
For any nontrivial element $w$ in a free group $\Free$, there exists a subgroup $H\leqslant \Free$ of index $|w|+1$ not containing $w$.
\end{lem}

\begin{figure}[H] 
  \centering
  \begin{tikzpicture}[shorten >=1pt, node distance=1.2 and 1.2, on grid,auto,>=stealth']
   \node[state, accepting] (0) {};
   \node[state] (1) [right = of 0]{};
   \node[state] (2) [right = of 1]{};
   \node[state] (3) [right = of 2]{};
   \node[state] (4) [right = of 3]{};
   \node[] (5) [right = 0.8 of 4]{};

   \path[->]
        (0) edge[loop,red,dashed,min distance=10mm,in=55,out=125]
            node[left = 0.1] {\scriptsize{$b$}}
            (0)
            edge[blue,dashed, bend right]
            (2);

    \path[->]
        (1) edge[blue]
        node[above] {\scriptsize{$a$}}
            (0)
            edge[red,dashed, bend left]
            (2);

    \path[->]
        (2) edge[red]
            (1)
            edge[blue]
            (3);

    \path[->]
        (3) edge[red]
            (4)
            edge[blue,dashed, bend right]
            (4);
            
    \path[->]
        (4) edge[blue,dashed, bend right = 40]
            (1)
            edge[red,dashed, bend right]
            (3);
\end{tikzpicture}
\caption{An index 5 subgroup of $\Free_{\set{a,b}}$ not containing the commutator $a^{-1} b^{-1} a b$: $H=\langle b,\, aba,\, a^{-1}b^2a,\, a^{-1}b^{-1}a^2b^{-1}a^{-1}ba,$ $a^{-1}b^{-1}ab^2a^{-1}ba,\,\allowbreak a^{-1}b^{-1}aba^2\rangle$}
\label{fig: fi not containing w}
\end{figure}
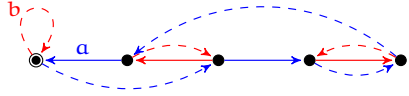

Since finite index subgroups in finitely generated groups have finitely many conjugates, the intersection $c(H)=\bigcap_{x\in \Fn} H^x$ is finite and provides a finite index normal subgroup of $\Fn$ not containing $w$ as well. This proves the well known result that (finitely generated) free groups are residually finite. Pushing the same idea a bit further, M.\ Hall reached a stronger result. 

\begin{thm}[\citenr{hall_jr_topology_1950}] \label{thm: Hall vff}
If $H$ is a finitely generated subgroup of a free group~$\Free$, then $H$ is a free factor of a finite-index subgroup of $\Free$; that is
 \begin{equation} \label{eq: Hall vff}
H\leqslant\fg \Free \Imp \exists K \st H \leqslant\ff K \leqslant\fin \Free.
 \end{equation}
\end{thm}

Finally, combining \Cref{thm: Hall vff} with \Cref{lem: F rf}, we see that free groups are \defin[fully LERF group]{fully LERF}, \ie we can choose the subgroup $K$
in \Cref{eq: Hall vff} avoiding any given finite subset~$S \subseteq \Free \setmin H$. 

\begin{thm}[\citenr{hall_jr_topology_1950}]
Finitely generated free groups are fully LERF. In particular, every finitely generated subgroup of $\Free$ is closed in the profinite topology.
\end{thm}

\begin{figure}[H] 
  \centering
  \begin{tikzpicture}[shorten >=1pt, node distance=1.2 and 1.2, on grid,auto,>=stealth']
   \node[state, accepting] (0) {};
   \node[state] (1) [right = of 0]{};
   \node[state] (2) [right = of 1]{};
   \node[state] (3) [right = of 2]{};
   \node[state] (4) [right = of 3]{};
   \node[state] (5) [below = 1 of 2]{};

   \path[->]
        (0) edge[red,dashed, bend left = 50] (2)
            edge[blue, bend right = 25] (5);

    \path[->]
        (1) edge[red, bend left= 25]
        node[below] {\scriptsize{$b$}}
            (0)
            edge[blue, bend right= 25]
         node[above] {\scriptsize{$a$}}
            (0);

    \path[->]
        (2) edge[blue]
            (1)
            edge[red]
            (3);

    \path[->]
        (3) edge[blue] (4)
            edge[red, bend left] (5);
            
    \path[->]
        (4) edge[loop,red,dashed,min distance=10mm,in=35,out=-35](4)
        (4) edge[blue,dashed, bend right = 50]
            (2);
            
    \path[->]
        (5) edge[red]
            (1)
            edge[blue,dashed, bend left]
            (3);
\end{tikzpicture}
\caption{The Stallings automaton of a finite index subgroup of $\Free_2$ containing $H = \gen{a b^2,a^{-2} b^{4}, a^{-1} b}$ (in solid lines) as free factor, and not containing the element $a^{-2} ba$.}
\label{fig: fi containing H not containing w}
\end{figure}
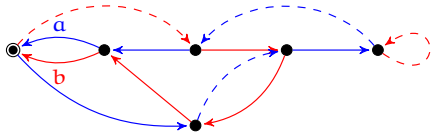

\medskip

To finish this subsection, we use the graphical characterization in \Cref{prop: fi iff}
to prove the classic Greenberg--Stallings \Cref{thm: Greenberg-Stallings}. Recall that two subgroups are called \defin[commensurable subgroups]{commensurable} if their intersection has finite index in both of them; 
and, for a subgroup $H$ of a group $G$, the \defin{commensurator} of $H$ in $G$ is
$\operatorname{Comm}_{G}(H) = 
\Set{
g \in G \st \ind{H \cap H^{g}}{H} < \infty \text{ and } \ind{H \cap H^{g} }{H^{g}} < \infty
}$ (a subgroup of $G$ containing $H$). 

Here, the idea is to use the fact that the Stallings automaton of a finitely generated infinite index subgroup of a free group always has a deficient vertex to reinterpret finite index in terms of the commensurator as follows.

\begin{lem}
A nontrivial finitely generated subgroup $H \leqslant\fg \Fn$ has finite index in $\Fn$ if and only if $\operatorname{Comm}_{\Fn}(H) = \Fn$.
\end{lem}

\begin{thm}[{\citenr{greenberg_commensurable_1974}}; {\citenr{stallings_topology_1983}}] \label{thm: Greenberg-Stallings}
Let $H$ and $K$ be commensurable finitely generated subgroups of a free group.
Then, the intersection $H \cap K$ has finite index in the subgroup  $\gen{H \cup K}$.
(This property is sometimes called \defin{Property S}.)
\end{thm}

\begin{prop}
Every nontrivial finitely generated subgroup $H$ of $\Fn$ has finite index in $\operatorname{Comm}_{\Fn}(H)$. 
\end{prop}

\subsubsection{\bf Coset enumeration}
A fundamental problem in algorithmic group theory is that of coset enumeration, which, grossly speaking, consists of listing a set of coset representatives for a finitely generated subgroup $H$ of a group $G$ given by a finite presentation. A classical result by \citeauthor{todd_practical_1936} in \cite{todd_practical_1936} states that if the index $\ind{H}{G}$ is finite, then such a list can be always obtained algorithmically (see also \cite{knuth_simple_1970} for the alternative --- and more modern --- Knuth--Bendix algorithm). As originally stated, Todd--Coxeter's algorithm is somewhat complicated, but it turns out that the ideas involved admit a very transparent interpretation in terms of Stallings automata (see~\cite{stallings_todd-coxeter_1987}).

Let us reformulate this interpretation in our automata language. Suppose that
\smash{$\normalcl{R} \xinto{\ } \Free[A] \xonto{\rho\, } G$} is a  finite presentation for a group $G$.
Then, if $S$ is a subset of $\Free[A]$, and $H = \gen{S \rho} \leqslant G$, it is clear that the sets $H \backslash G$ and $G/H$ (of right and left cosets of $H$ in $G$) are in bijection with the set $H \rho \preim \backslash \Free[A]$ (of cosets of the full preimage $H\rho\preim$ in $\Free[A]$); hence the index of $H$ is $G$ is expressible in the background group $\Free[A]$ as $\Ind{H}{G}=\Ind{H \rho\preim}{\Free[A]}= \Ind{\gen{S}\normalcl{R}}{\Free[A]}$. Therefore, it is enough to compute the Schreier automaton of the subgroup $H\rho\preim =\gen{S}\normalcl{R}\leqslant \Free[A]$ to recover the cosets of $H$ in $G$. This approach entails the obstacle that the target subgroup $H\rho\preim = \gen{S}\normalcl{R}$ may not be finitely generated. %
However, it turns out that if we restrict our scope to subgroups $H\leqslant\fin G$ of finite index (equivalently $H\rho\preim \leqslant\fin \Free[A]$), then the obstacle can be overcome.  

So, assume that $\ind{H}{G}<\infty$ (and hence $S$ is finite). Note that then the automaton $\Atiii$  obtained after attaching the flower automaton $\flower(R)$ to every vertex in the (potentially infinite) automaton $\schreier(\gen{S})$ recognizes $\gen{S}$ at the basepoint, and $\gen{R}$ at every vertex;
that is, $\gen{\Atiii} = \gen{S} \normalcl{R} = H \rho\preim $. Now, the key idea is to attach these (theoretically infinitely many) flower automatons $\flower(R)$ sequentially, starting from $\flower(S)$ and folding after each layer addition. Formally, we set $\Ati_{\!0} = \flower(S)$ and, for each $i \geq 0$,
$\Ati_{\!i+1}$ is the automaton obtained after attaching $\flower(R)$ to every vertex in~$\Ati_{\!i}$ and folding.

By construction, every automaton $\Ati_{\!i}$ in this sequence is reduced, and recognizes the subgroup $\gen{S}$ at the basepoint, and $\gen{R}$ at every vertex \emph{except possibly the ones corresponding to the $\flower(R)$'s attached at the $i$-th step}.
Since the subgroup $H \rho \preim$ is, by hypothesis, of finite index in $\Free[A]$, we know that
$\stallings(\Atiii) = \stallings(H \rho \preim) = \schreier(H \rho \preim)$ will be a finite and saturated reduced automaton recognizing $\gen{S}$ at the basepoint and~$\gen{R}$ at every vertex. This means that the sequence $(\Ati_{\!i})_{i \geq 0}$ must stabilize after a finite number of steps precisely at $\schreier(H \rho \preim)$, which is therefore computable.

\begin{thm}[\citenr{todd_practical_1936}]
Let $H$ be a finite index subgroup of a group $G$ given by a finite presentation. Then, the Schreier automaton $\schreier(H)$ (and hence the index and a set of coset representatives of $H$) is computable from any finite generating set for $H$.
\end{thm}

\begin{cor}
The subgroup membership problem restricted to finite index subgroups is always decidable. 
\end{cor}

\subsubsection{\bf The Herzog–Schönheim conjecture}

To close this section, let us survey on an interesting open problem concerning finite index subgroups, in which there has been recent progress using Stallings techniques: the Herzog--Sch\"onheim conjecture. 

The story began with the following arithmetic problem: a \emph{cover} of the set of integers $\mathbb{Z}$ is a finite collection of remainder classes $r_i+d_i\mathbb{Z}$, $i=1,\ldots ,s$, disjoint to each other, and such that $\mathbb{Z}=\bigsqcup_{i=1}^s (r_i+d_i\mathbb{Z})$. Typical examples of covers are $\mathbb{Z}=0+1\mathbb{Z}$ (the trivial one), $\mathbb{Z}=2\mathbb{Z}\sqcup (1+2\mathbb{Z})$, $\mathbb{Z}=4\mathbb{Z}\sqcup (2+4\mathbb{Z})\sqcup (1+2\mathbb{Z})$, $\mathbb{Z}=4\mathbb{Z}\sqcup (2+12\mathbb{Z})\sqcup (6+12\mathbb{Z})\sqcup (10+12\mathbb{Z})\sqcup (1+2\mathbb{Z})$, etc. This concept was introduced by P. Erd\"os in~\parencite{erdhos_integers_1950}, who conjectured that in any such non-trivial cover the largest index $d_s$ must appear at least twice. This conjecture was proved independently by Davenport, Rado, Mirsky, and Newman (see~\parencite{znam_exactly_1970} for the precise story) using analysis of complex functions; furthermore, it was proved that this largest index $d_s$ appears at least $p$ times, where $p$ is the smallest prime dividing $d_s$ (among other related results). See~\parencite{ginosar_tile_2018} for an alternative modern proof using group representations. 

In 1974, M. Herzog and J. Sch\"onheim extended Erd\"os’ conjecture to arbitrary finitely generated groups and launched the following much more general conjecture: if $\{ H_1g_1,\ldots ,H_sg_s\}$, $s\geqslant 2$, is a nontrivial \emph{coset partition} of a finitely generated group $G$ (meaning that $H_i$ are finite index subgroups of $G$, and $G=\bigsqcup_{i=1}^s H_i g_i$) then the list of indices $d_i=\ind{H_i}{G}$ must contain at least a repetition; one usually refers to this fact by saying that the coset partition has \emph{multiplicity}. Although during the following decades several papers appeared providing partial results, the conjecture in its general form remains still open today. Among the known partial results we highlight that of Berger--Felzenbaum--Fraenkel~\parencite{berger_remark_1987}, proving the conjecture for pyramidal groups (a subclass of finite solvable groups), and the recent one given by Margolis--Schnabel~\parencite{margolis_herzog-schonheim_2019}, proving it true for all finite groups of order less than~1440.

The common approach to the Herzog--Sch\"onheim conjecture is through finite groups. Indeed, given a finitely generated group $G$, and a coset partition $G=\bigsqcup_{i=1}^s H_i g_i$, one can consider $N\normaleq G$ to be the intersection of the normal cores of the $H_i$'s (the normal core of $H_i\leqslant G$ is $c(H_i)=\bigcap_{g\in G} H_i^g$, still of finite index in $G$ because of the finite generability of $G$). Clearly, $N$ is a normal finite index subgroup of $G$ and the given coset partition projects to a coset partition of $G/N$, namely $G/N =\bigsqcup_{i=1}^s H_i/N\, (Ng_i)$, with exactly the same list of indices, $\ind{H_i}{G} = \ind{H_i/N}{G/N}$. Hence, the general Herzog--Sch\"onheim conjecture reduces to the case of finite groups. 

In a series of papers by F. Chouraqui
(see~\parencite{chouraqui_herzog-schonheim_2019,chouraqui_about_2019,chouraqui_approach_2020,chouraqui_space_2018}),
the author adopts a completely different approach to attack Herzog--Sch\"onheim conjecture: instead of finite groups, she considers free groups of finite rank $\Free[r]$. Indeed, in~\parencite{chouraqui_herzog-schonheim_2019} she observes that, in a finitely generated group viewed as a quotient of a free group, $G=\Free[r]/N$, any coset partition $G=\bigsqcup_{i=1}^s H_i g_i$ can be lifted up by taking full preimages through the canonical projection $\rho\colon \Free[r] \onto G$, getting a coset partition of the free group, namely $\Free[r]=\bigsqcup_{i=1}^s (H_i\rho\preim)\widetilde{g}_i$ (where $\widetilde{g}_i\in \Free[r]$ is any $\rho$-preimage of $g_i\in G$), with the same list of indices, $\ind{H_i}{G} = \ind{H_i\rho\preim}{\Free[r]}$. Hence, the general Herzog--Sch\"onheim conjecture reduces to the case of finitely generated free groups. 

We are now in the land of Stallings automata: a nontrivial coset partition  $\Free[r]=\bigsqcup_{i=1}^s H_i w_i$, can be codified with the finite list of Stallings automata $\Ati_{\!1},\ldots ,\Ati_{\!s}$ for the cosets $H_1 w_1,\ldots, H_s w_s$, respectively; these are $s\geqslant 2$ saturated and finite Stallings automata, and the goal is to show that the equality $\Free[r]=\bigsqcup_{i=1}^s H_i w_i$ implies that at least two of them must have the same number of vertices, \ie $\card \Verts \Ati_{\!i}=\card \Verts \Ati_{\!j}$, for some $i \neq j$. In the above mentioned series of papers Chouraqui finds several conditions on the $\Ati_{\!i}$'s to ensure this conclusion, proving several special cases of Herzog--Sch\"onheim conjecture. 

Fix a nontrivial coset partition $\Free[r]=\bigsqcup_{i=1}^s H_i w_i$, with $d_i=\ind{H_i}{\Free[r]}$ and $2\leqslant d_1\leqslant \cdots \leqslant d_s$. In~\parencite{chouraqui_herzog-schonheim_2019}, the author considers the corresponding transition monoids $M_1,\ldots, M_s$ (see \Cref{sec: automata}), which in this case are subgroups of the corresponding symmetric group, $M_i\leqslant S_{d_i}$ (because the Stallings automata $\Ati_{\!1},\ldots, \Ati_{\!s}$ are saturated and so, all partial injections of the form $\tau_w$ are real permutations). A relevant contribution to the conjecture from~\parencite{chouraqui_herzog-schonheim_2019} is stated below.

\begin{thm}[\citenr{chouraqui_herzog-schonheim_2019}]
In the above situation, if $M_s$ contains a cycle of length $d_s$ then the coset partition $\Free[r]=\bigsqcup_{i=1}^s H_i w_i$ has mutiplicity. 
\end{thm}

It is mentioned that, in some statistical sense, the condition on $M_s$ is satisfied with probability at least $3/4$, as $d_s \to \infty$.
This reinforces the belief that the Herzog--Sch\"onheim conjecture is true.
Besides this one, some other results are also proved, assuming certain technical relations between the group $M_s$ and the top indices $d_{s-3}\leqslant d_{s-2}\leqslant d_{s-1}\leqslant d_s$ in the coset partition. 

\subsection{\bf Intersections}

Intersections of subgroups of free groups is a research topic with a long-standing and interesting history. We state below, for an arbitrary group $G = \pres{A}{R}$, the motivating concept (first considered by \citeauthor{howson_intersection_1954} over free groups in \parencite{howson_intersection_1954}), and the natural associated algorithmic problem.

\begin{defn} \label{def: Howson group}
A group $G$ is said to satisfy the \defin{Howson property} (or to be \defin[Howson group]{Howson}, for short) if the intersection of any pair of finitely generated subgroups of $G$ is again finitely generated.
\end{defn}

\begin{named}[Subgroup intersection problem, $\SIP(G)$] \label{def: SIP} 
Decide\margin{Subgroup intersection problem (\SIP)}, given two finite sets of words $U,V$ in the generators of $G$, whether the intersection $\gen{U} \cap \gen{V}$ is finitely generated; and, in the affirmative case, compute a generating set for the intersection.
\end{named}

The related natural intersection problem for cosets is stated below. 

\begin{named}[Coset intersection problem, $\CIP(G)$] \label{def: CIP} 
Decide\margin{Coset intersection problem (\CIP)}, given two words $u,v$ and two finite sets of words $U,V$ in the generators of $G$, whether the coset intersection $u\gen{U} \cap v\gen{V}$ is empty; and, in the negative case, compute a representative for the intersection coset.
\end{named}

As stated in \Cref{cor: ranks of subgroups}, the subgroups of a finitely generated free group can have  any (finite or countably infinite) rank. However,
in \citeyear{howson_intersection_1954} %
Howson proved that the intersection of two finitely generated subgroups of the free group must be again finitely generated. Below we see how the Stallings machinery provides a neat and algorithmic-friendly proof for this remarkable fact, and furthermore allows to compute a basis for the intersection, and to prove Hanna Neumann's bound \eqref{eq: Hanna-Neumann bound} on its rank.

The key construction used for the geometric proof of this fact is that of product (or \defin{pull-back}) of automata, defined below.

\begin{defn} \label{def: product of automata}
Let $\Ati_{\!1} = (\Ati_{\!1}, P_1,Q_1)$ and $\Ati_{\!2} = (\Ati_{\!2},P_2,Q_2)$ be $A$-automata. The \defin[product of automata]{(tensor or categorical) product} of $\Ati_{\!1}$ and $\Ati_{\!2}$, denoted by $%
\Ati_{\!1} \times \Ati_{\!2}$, is the $\Alfi$-automaton with:
\begin{itemize}[beginpenalty=10000]
    \item vertex set the Cartesian product $\Verts\Ati_{\!1} \times \Verts \Ati_{\!2}$;
    \item an arc $(\verti_1,\verti_2) \xarc{\alfi\,} (\vertii_1,\vertii_2)$ for every pair of arcs $\verti_1 \xarc{\alfi\,} \vertii_1$ in~$\Ati_{\!1}$, and $\verti_2 \xarc{\alfi\,} \vertii_2$ in~$\Ati_{\!2}$ with the same label $\alfi \in \Alfi$;
    \item set of initial vertices $P = \set{(\verti_1,\verti_2) \st \verti_1 \in P_1 \text{ and } \verti_2 \in P_2}$;
    \item set of terminal vertices $Q = \set{(\vertii_1,\vertii_2) \st \vertii_1 \in Q_1 \text{ and } {\vertii_2 \in Q_2}}$.
\end{itemize}
\end{defn}
As usual, the corresponding versions of product for labelled and unlabelled digraphs, undirected graphs, etc.\ are obtained after omitting in \Cref{def: product of automata} the conditions involving the missing parts.

\begin{figure}[H]
    \centering
\begin{tikzpicture}[shorten >=1pt, node distance=1.25cm and 1.25cm, on grid,auto,>=stealth']
   \node (00) {};
   \node[state,accepting] (0) [below = 0.5 of 00]{};
   \node[] (bp1) [below right =0.1 and 0.15 of 0] {$\scriptscriptstyle{1}$};
   \node[state,accepting] (0') [right = 0.5 of 00]{};
   \node[] (bp1) [below right =0.1 and 0.15 of 0'] {$\scriptscriptstyle{2}$};
   \node[state,accepting,blue] (00') [right = 0.5 of 0]{};
   \node[state] (1) [below = 1 of 0] {};
   \node[state] (2) [below = of 1] {};
   \node[] (G') [left = 0.8 of 1] {$\scriptstyle{\Ati_{\!1}}$};
   \node[state] (1') [right = 1 of 0'] {};
   \node[state] (2') [right = 1.25 of 1'] {};
   \node[] (G) [right = 0.7 of 2'] {$\scriptstyle{\Ati_{\!2}}$};
   \node[state,blue] (11') [below right = 1 and 1 of 00'] {};
   \node[state,blue] (22') [below right = of 11'] {};
    \node[blue] (GG') [right = 0.9 of 22'] {$\scriptstyle{\Ati_{\!1} \times \Ati_{\!2}}$};

   \path[dashed]
       (0') edge[] (1');
   \path[dashed]
       (0) edge[] (1);
   \path[dashed,blue]
       (00') edge[] (11');
   \path[->]
       (1') edge[]
            node[pos=0.5,above=-.1mm] {$\alfi$}
            (2');
   \path[->]
        (1) edge[]
            node[pos=0.5,left=-.1mm] {$\alfi$}
            (2);
   \path[->,blue]
        (11') edge[]
            node[pos=0.5,above right] {$\alfi$}
            (22');
\end{tikzpicture}
    \caption{Scheme of the product (in blue) of two pointed automata (in black)}
    \label{fig doubly enriched arc}
\end{figure}
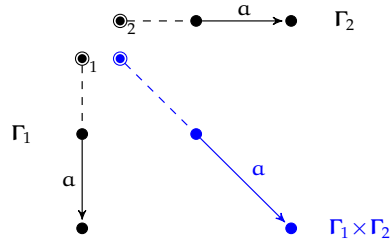

The link between intersections of subgroups of the free group and products of Stallings automata is easily checkable by inspection (just chasing common $\bp$-rounds in the Stallings automata of the intersecting subgroups).

\begin{lem} \label{lem: product intersection}
Let $\Ati_{\!1}$ and $\Ati_{\!2}$ be deterministic $\Alfi$-automata. Then, the product ${\Ati_{\!1} \times \Ati_{\!2}}$ is a deterministic $\Alfi$-automaton recognizing the intersection of the subgroups recognized by $\Ati_{\!1}$ and $\Ati_{\!2}$ respectively; that is, $\gen{\Ati_{\!1} \times \Ati_{\!2}} = \gen{\Ati_{\!1}} \cap \gen{\Ati_{\!2}}$. 
\end{lem}

We note that, in general, the product of two core automata is not necessarily core (not even connected); so it is necessary to take the core of the product in order to reach the Stallings automaton of the intersection.  It is straightforward to see that the same kind of argument allows to extend the result to intersections of arbitrary cosets of the free group.

\begin{prop}\label{prop: St(intersection)}
Let $H,K$ be subgroups of a free group $\Free[\Alfi]$, and let $u,v \in \Free[\Alfi]$. Then, the Stallings automaton of the coset intersection is $\stallings(H u \cap K v) = \core (\stallings(H u) \times \stallings(K v))$. 
In particular, the Stallings automaton of the intersection of subgroups is $\stallings(H \cap K) = \core (\stallings(H) \times \stallings(K))$.
\end{prop}

So, if $H_1$ and $H_2$ are finitely generated, and hence (from \Cref{thm: Stallings bijection}) $\stallings(H_1)$ and $\stallings(H_2)$ are finite and computable, then $\stallings(H_1 \cap H_2)$ is finite and computable as well. This proves Howson property for free groups; once a (finite) Stallings automaton for the intersection $H_1 \cap H_2$ is obtained, a basis for $H_1 \cap H_2$ can be computed using the standard method in~\Cref{prop: generators from tree}.

\begin{thm}[\citenr{howson_intersection_1954}]\label{thm: Howson}
Free groups are Howson.
\end{thm}

\begin{thm}
The Subgroup Intersection Problem $\SIP(\Free)$is computable.
\end{thm}

In a similar way the Coset Intersection Problem is also decidable: given cosets $Hu$ and $Kv$, consider the Stallings automata $\stallings(Hu)$ (with basepoint $\bp_H$ and final vertex $\verti$) and $\stallings(Kv)$ (with basepoint $\bp_K$ and final vertex $\vertii$). It is straightforward to see that the intersection $Hu\cap Kv$ is empty if and only if the vertices $(\bp_H, \bp_K)$ and $(\verti, \vertii)$ lie in different connected components of the product $\stallings(Hu)\times \stallings(Kv)$; and otherwise, the label, say $w$, of a walk from $(\bp_H, \bp_K)$ to $(\verti, \vertii)$ satisfies $Hu=Hw$ and~$Kv=Kw$; therefore $w\in Hu\cap Kv$, and hence $Hu\cap Kv=(H\cap K)w$.

\begin{thm}
The coset intersection problem $\CIP(\Free)$ is computable.
\end{thm}

After Howson's result, the natural quest for bounds for the rank of the intersection of two finitely generated subgroups $H,K \leqslant \Free$ became a popular topic, and ended up being one of the most famous open problems in geometric group theory. The first (partial) answer to this question was given by Howson himself as \margin{Howson's bound}
 \begin{equation} \label{eq: Howson bound}
\rrk (\Sgpi \cap \Sgpii) \,\leq\, 2\, \rrk (\Sgpi)  \, \rrk(\Sgpii) + \rrk(\Sgpi) + \rrk(\Sgpii) \,,
 \end{equation}
where $\rrk(H) = \max \set{\rk(H) -1, 0}$ is the \defin{reduced rank} of $H$, providing also examples of subgroups ${\Sgpi,\Sgpii \leqslant \Free[n]}$ where $\rrk (\Sgpi \cap \Sgpii) = \rrk (\Sgpi) \, \rrk (\Sgpii)$ (see \Cref{fig: pullback}). Shortly after, in 1956, \citeauthor{neumann_intersection_1956} improved Howson's bound and proved that \margin{Hanna Neumann bound}
 \begin{equation} \label{eq: Hanna-Neumann bound} 
\rrk (\Sgpi \cap \Sgpii) \,\leq\, 2\, \rrk (\Sgpi)  \, \rrk(\Sgpii) 
 \end{equation}
which, again, is not hard to derive from the product automaton using elementary graph theoretical arguments. %
The famous conjecture by Hanna Neumann was that the factor `2' in equation \eqref{eq: Hanna-Neumann bound} can be removed, making the following bound tight:
\margin{Hanna Neumann conjecture}
 \begin{equation} \label{eq: Hanna-Neumann conj} 
\rrk (\Sgpi \cap \Sgpii) \,\leq\, \rrk (\Sgpi)  \, \rrk(\Sgpii) \,. 
  \end{equation}

\begin{exm} 
Consider the free group $\Free[\set{\!a,b\!}]$, and its subgroups $H=\gen{ab,a^3, a^{-1}ba}$ and $K=\gen{b, a^3, a^{-1}bab^{-1}a}$. Their intersection $H\cap K$ realizes tightly the Hanna Neumann conjecture bound, namely, $\rrk (\Sgpi \cap \Sgpii) =4=(3-1)(3-1)= \rrk (\Sgpi)  \, \rrk(\Sgpii)$ (see \Cref{fig: pullback}).

\begin{figure}[H]
\centering
\begin{tikzpicture}[shorten >=1pt, node distance=1.2cm and 2cm, on grid,auto,auto,>=stealth']

\newcommand{\dx}{1.3}
\newcommand{\dy}{1.2}
\node[] (0)  {};
\node[state,accepting] (a1) [right = \dy-1/3 of 0] {};
\node[state] (a2) [right = \dx of a1] {};
\node[state] (a3) [right = \dx of a2] {};
\node[state] (a4) [right = \dx of a3] {};

\node[state,accepting] (b1) [below = 2*\dy/3 of 0] {};
\node[state] (b2) [below = \dy of b1] {};
\node[state] (b3) [below = \dy of b2] {};

\foreach \y in {1,...,4}
\foreach \x in {1,...,3} 
\node[state] (\x\y) [below right = (\x-1/3)*\dy and (\y-1/3)*\dx of 0] {};
    
\node[state,accepting] () [below right = 2*\dy/3 and 2*\dx/3 of 0] {};

\path[->]
    (a1) edge[red,loop above,min distance=7mm,in=205,out=155]
    node[above = 0.1] {$b$}
    (a1)
    (a1) edge[blue] node[below] {$a$} (a2)
    (a2) edge[blue] (a3)
    (a3) edge[blue,bend right=30] (a1)
    (a3) edge[red] (a4)
    (a4) edge[blue,loop above,min distance=7mm,in=25,out=-25] (a4)

    (b1) edge[blue,bend right=25] (b2)
    (b2) edge[red,bend right=25] (b1)
    (b2) edge[blue] (b3)
    (b3) edge[blue,bend left=40] (b1) 
    (b3) edge[red,loop left,min distance=7mm,in=295,out=245] (b3);
             
\path[->]   
    (11) edge[blue] (22)
    (22) edge[blue] (33) 
    (33) edge[blue, bend right=25] (11)
            
    (21) edge[blue] (32)
    (32) edge[blue] (13)
    (13) edge[blue] (21)
        
    (12) edge[blue] (23)
    (23) edge[blue] (31) 
    (31) edge[blue] (12)
            
    (14) edge[blue] (24)
    (24) edge[blue] (34)
    (34) edge[blue, bend right=25] (14)
            
    (21) edge[red] (11)
    (23) edge[red] (14)
    (33) edge[red] (34)            
    (31) edge[red,loop left,min distance=7mm,in=250,out=200] (31);
    
\node[] (i) [right = 1 of 24] {$=$};  
\node[state,accepting] (1) [above right = 0.7 and 2.25 of i] {};
\node[state] (1a) [below left = \dy and \dx/3 of 1] {};
\node[state] (1b) [below right = \dy and \dx/3 of 1] {};
\node[state] (2) [right = \dx*0.7 of 1b] {};
\node[state] (2a) [above left = \dy and \dx/3 of 2] {};
\node[state] (2b) [above right = \dy and \dx/3 of 2] {};
\node[state] (3) [right = \dx*0.7 of 2b] {};
\node[state] (3a) [below left = \dy and \dx/3 of 3] {};
\node[state] (3b) [below right = \dy and \dx/3 of 3] {};
\node[state] (0) [left = \dx*0.7 of 1a] {};
\node[state] (0a) [above left = \dy and \dx/3 of 0] {};
\node[state] (0b) [above right = \dy and \dx/3 of 0] {};

\path[->]
    (0b) edge[red] (1)
    (1b) edge[red] (2)
    (3) edge[red] (2b)
    
    (0b) edge[blue] (0)
    (0a) edge[blue] (0b)
    (0) edge[blue] (0a)
    
    (1b) edge[blue] (1)
    (1a) edge[blue] (1b)
    (1) edge[blue] (1a)
    
    (2) edge[blue] (2b)
    (2b) edge[blue] (2a)
    (2a) edge[blue] (2)
    
    (3) edge[blue] (3b)
    (3b) edge[blue] (3a)
    (3a) edge[blue] (3)
    (3b) edge[red,loop,min distance=7mm,in=25,out=-25] (3b);
\end{tikzpicture}
\vspace{-15pt}
\caption{An intersection where $\rrk (\Sgpi \cap \Sgpii)
= \rrk (\Sgpi)  \, \rrk(\Sgpii)$}
\label{fig: pullback}
\end{figure}
\end{exm}

In fact, if one looks at the full product of automata $\stallings(H) \times \stallings(K)$ instead of at its core, one can see that the Hanna Neumann's bound \eqref{eq: Hanna-Neumann bound} affects indeed the intersections of $H$ \emph{with all the conjugates of~$K$}. 
The corresponding claim for Hanna Neumann's conjecture (removing the $2$) was stated by \citeauthor{kovacs_intersections_1990} in \parencite{kovacs_intersections_1990}, and received the name of \defin{Strengthened Hanna Neumann conjecture}:
 \begin{equation} \label{eq: str Hanna-Neumann conj} 
\sum_{x\in H \backslash \Free / K} \rrk (\Sgpi \cap \Sgpii^x) \,\leq\, \rrk (\Sgpi)  \, \rrk(\Sgpii) \,, 
  \end{equation}
where every summand corresponds to a connected component in the product $\stallings(H) \times \stallings(K)$. 

After many unsuccessful attempts and partial results, two independent proofs by
\citeauthor{friedman_sheaves_2015} and \citeauthor{mineyev_submultiplicativity_2012}
appeared almost simultaneously more than fifty years later (see~\parencite{friedman_sheaves_2015,mineyev_submultiplicativity_2012}, and the remarkable respective simplifications by W.~Dicks in \cite[Appendix B]{friedman_sheaves_2015} and the unpublished note~ \parencite{dicks_simplified_2012}).

\begin{thm}[\citenr{friedman_sheaves_2015}; \citenr{mineyev_submultiplicativity_2012}]
The Strengthened Hanna Neumann conjecture holds in free groups.
\end{thm}

In the recent preprint \cite{linton_intersections_2021} \citeauthor{linton_intersections_2021} proves a similar formula in terms of ranks instead of reduced ranks: 
 \begin{equation} 
\sum_{[\Sgpi\cap \Sgpii^x]} \rk (\Sgpi \cap \Sgpii^x) \,\leq\, C(\rk (\Sgpi)  \, \rk(\Sgpii))^2 \,,
  \end{equation}
where $C>0$ is a universal constant, and the sum is over the conjugacy classes of subgroups of the form $\Sgpi\cap \Sgpii^x$.

Among other generalizations, alternative proofs extending the Strengthened Hanna Neumann conjecture to the realms of
right-orderable groups~\parencite{antolin_Kurosh_rank_intersections_2014},
pro-$p$ groups~\parencite{jaikin-zapirain_approximation_2017}, and surface groups \parencite{antolin_hanna_2021} were given, respectively, by
\citeauthor{antolin_Kurosh_rank_intersections_2014};
\citeauthor{jaikin-zapirain_approximation_2017}; and \citeauthor{antolin_hanna_2021} shortly after.

To illustrate the applicability of the pull-back technique, let us give a short and elegant argument to see the following folklore result restricted to the free ambient case (the general proof uses Kurosh Subgroup Theorem).

\begin{prop}
Let $G$ be a group and $H,K,H',K'$ be subgroups of $G$. If $H\leqff K$ and $H'\leqff K'$ then $H\cap H'\leqff K\cap K'$. 
\end{prop}

Suppose $G$ is free. Let us first see that if $H \leqslant\ff K \leqslant G$ and $L \leqslant G$ then $H \cap L \leqslant\ff K \cap L$.
Take a basis $A$ of $K$ extending a basis of $H$, and note that $\stallings(H,A)$ is just a bouquet with some petals (those labelled by elements in $H$). Consider $\stallings(K\cap L, A)$  and let us compute ${H\cap L}=H\cap (K\cap L)$ using the product automata: clearly, $\stallings(H,A)\times \stallings(K\cap L,A)$ is the subautomaton of $\stallings(K\cap L,A)$ determined by the edges with labels in $H$; therefore, $H\cap L\leqff K\cap L$ as claimed.

Now, applying this argument twice we obtain $H\cap H'\leqff K\cap H'\leqff K \cap K'$ and therefore $H\cap H'\leqff K\cap K'$, as we wanted to prove.

\subsubsection{\bf Malnormality}

Recall that a subgroup $H$ of an arbitrary group $G$ is said to be \defin[malnormal subgroup]{malnormal} (in~$G$) if for every $g \in G \setmin H$ the intersection $H \cap H^{g} = \Trivial$. This terminology emphasizes their character opposite to that of normal subgroups (were every intersection $H \cap H^{g} = H$).
We can use our geometric description of intersections within free groups to characterize malnormality in this family.

\begin{prop} \label{prop: malnormality}
Let $H$ be a subgroup of a free group $\Free$. Then $H$ is malnormal in $\Free$ if and only if every connected component of $\stallings(H) \times \stallings(H)$ not containing the basepoint is a tree.
\end{prop}

\begin{rem}
Note that, although the diagonal component of $\stallings(H) \times \stallings(H)$ is obviously $\stallings(H)$ itself (meaning that $H\cap H=H$), there must be other connected components in the product. These are the relevant ones in order to analyze the malnormality of $H$. 
\end{rem}

As a consequence of \Cref{prop: malnormality} we recover another well known algorithmic result on free groups.

\begin{cor}[\citenr{baumslag_malnormality_1999}]
In a free group it is algorithmically decidable whether a finitely generated subgroup is malnormal.
\end{cor}

\begin{exm}
The existence of malnormal subgroups of infinite rank in $\Free_2$ was needed at some point in the literature, and proved with quite technical algebraic arguments (see, \eg \parencite{das_controlled_2015}). Using Stallings automata, this fact can be proved by just inspecting a labelled automaton like the one in \Cref{fig: nonfg malnormal}. 
\begin{figure}[H]
\centering
\begin{tikzpicture}[shorten >=1pt, node distance=1.2cm and 2cm, on grid,auto,auto,>=stealth']

\newcommand{\dx}{0.85}
\newcommand{\dy}{0.8}
\node[] (0)  {};

\foreach \x in {0,...,4}
\foreach \y in {0,...,4} 
{
\node[state] (\x\y) [above right = \y*\dy and \x*\dx of 0] {};
\node[state,accepting] (bp)  {};
}
\node[state] (50) [right = \dx*1.2 of 40] {};
\node[] (60) [right = \dx of 50] {};
\node[state] (05) [above = \dy*1.2 of 04] {};
\node[] (06) [above = \dy of 05] {};
\node[] (444) [above right = \dy*0.6 and \dx*0.6 of 44] {\raisebox{5pt}{\reflectbox{$\ddots$}}};
\node[state] (55) [above right = \dy*1.2 and \dx*1.2 of 44] {};
\node[] (66) [above right = \dy/2 and \dx/2 of 55] {\raisebox{5pt}{\reflectbox{$\ddots$}}};
\draw[->,blue,dotted,thick] (40) edge (50);
\draw[->,blue,dotted,thick] (50) edge (60);
\draw[->,red,dotted,thick] (04) edge (05);
\draw[->,red,dotted,thick] (05) edge (06);
\draw[->,red,dotted,thick] (50) edge (55);
\draw[->,blue,dotted,thick] (05) edge (55);

\path[draw=none,red] (00) edge node[left] {$b$} (01);
\path[draw=none,blue] (00) edge node[above] {$a$} (10);

\foreach \x [count=\xi] in {0,...,3}
{\draw[->,blue] (\x0) edge (\xi0);
\draw[->,red] (0\x) edge (0\xi);}

\foreach \y [count=\yi] in {1,...,4}{
\pgfmathsetmacro{\yy}{\y-1}
\foreach \x [count=\xi] in {0,...,\yy}
{
\draw[->,blue] (\x\y) edge (\xi\y);
\draw[->,red] (\y\x) edge (\y\xi);
}}
\end{tikzpicture}
\caption{A malnormal subgroup of \protect{$\Free[\{a,b\}]$} of infinite rank: ${H = \gen{a^kb^k a^{-k}b^{-k} \st k \geq 1}}$}
\label{fig: nonfg malnormal}
\end{figure}
\end{exm}

\subsection{\bf Relative order and spectra}

In \parencite{delgado_relative_2022}, the authors (together with A. Zakharov) consider natural generalizations of the concepts of root and order of an element in a group. Precise definitions are given below for a group $G$, a subset $S$ of $G$, $g \in G$, and $k \in \ZZ$.

\begin{defn} \label{def: k-root}
If $g^k \in S$ we say that $g$ is a \defin[$k$-root]{$k$-root of $S$}. We denote by $\sqrt[k]{S}$ the set of $k$-roots of $S$ in~$G$.
\end{defn}

\begin{defn}
The \defin[relative order]{(relative) order} of $g \in G$ in  $S$, denoted by $\ord_S(g)$, is the minimum $k\geqslant 1$ such that $g^k \in S$, and zero if there is no such $k$; \ie $\ord_{S}(g) = \max \left\{0 \,,\, \min \,\set{k \geq 1 \st g^k \in S}\right\}$. In particular, $\ord_S(g) = 1$ if and only if $g\in S$.
\end{defn}

\begin{defn}
The \margin{relative preorder} set of elements from $G$ of order $k$ in $S$ is denoted by~$\ofo{G}{k}{S}$; \ie $\ofo{G}{k}{S}=\Set{g \in G \st \ord_{S}(g) = k}$.
\end{defn}

\begin{defn}
The set $\ord_{S}(G)$ (of orders in $S$ of the elements in~$G$) is called the \defin[relative spectrum]{spectrum} of~$G$ in~$S$. That is, $\ord_S(G)=\set{ k \in \NN \st \ofo{G}{k}{S} \neq \varnothing }$.
\end{defn}

Throughout this section we will consider relative order and spectrum in subgroups and cosets of a finitely generated free group $\Free$. Note that if $H$ is a subgroup of $G$, then, for every $k \geq 1$, $\sqrt[k]{H} = \bigsqcup\nolimits_{i \divides k} \ofo{G}{i}{H}$ and $\ofo{G}{k}{H} = \sqrt[k]{H} \setmin \bigsqcup \nolimits_{k \smallneq i \divides k} \sqrt[i]{H}$.

Note also that, by definition, $\sqrt[0]{H} = G$ and $\sqrt[1]{H} = \ofo{G}{1}{H} = H$. We denote by $\sqrt{H}$ the set of elements of strictly positive order
(called \defin{nontrivial roots}) in~$H$, and by $\sqrt[*]{H}$ the set of elements of order at least $2$ (called \defin{proper roots}) in~$H$. Hence, $G=\ofo{G}{0}{H} \sqcup \sqrt{H}$ and $\sqrt{H}=H\sqcup \sqrt[*]{H}$.

The standard notions of order, periodicity, etc., correspond to instances of the introduced relative notions \wrt the trivial subgroup, in which case we omit the reference to the subgroup. For example, $\ord(g)$ denotes the standard order of $g$ (except for elements of infinite order for which we put $\ord(g)=0$), $\ofo{G}{k}{\phantom{H}}$ denotes the set of elements in $G$ of order $k$, and a group $G$ is periodic if ${0 \notin \ord(G)}$.

Another group-theoretical notion which admits a neat description in terms of the introduced language is that of purity.

\begin{defn}
A subgroup $H$ of a group $G$ is said to be \defin[pure subgroup]{pure} (in $G$) if it has no proper roots in $G$; that is, if $\sqrt{H}=H$ or, equivalently, $\Iset{H}{G} \subseteq \set{0,1}$.
\end{defn}

Note that nontrivial free factors are always pure. And the trivial subgroup is pure in a group $G$ if and only if $G$ is torsion-free. On the other hand, the whole group $G$ is obviously pure, and the intersection of pure subgroups is again pure. 

\begin{defn}
The \defin{pure closure} of a subgroup $H$ of $G$, denoted $\pcl{H}$, is the smallest pure subgroup of $G$ containing $H$; that is, the intersection of all pure subgroups of~$G$ containing $H$.

\end{defn}

\medskip
Stallings automata provide a natural approach to the study of relative orders in free groups. The naive starting idea is that an element $w \in \Free$ has order $k \geq 1$ in $H \leqslant \Free$ if and only if 
the sequence $\vvertii = (\bp w^i)_{i\geq 0}$ has its first repetition at $i=k$ (\ie $\bp = \bp w^k$).
Define a \defin{closed trail} (of length $k$) as a
sequence of vertices $\vvertii =(\vertii_0,\vertii_1,\ldots,\vertii_k= \vertii_0)$ in $\schreier(H)$ whose only repetition is $\vertii_k= \vertii_0$.
Then, for each such trail $\vvertii$, one can use \Cref{lem: recognized subgroup} to compute the set $\ofo{\Free}{\vvertii}{H}$ of elements $w$ realizing it (called \defin[preorbit]{the preorbit} of~$\vvertii$). The problem is, of course, that there are, in general, infinitely many trails of length $k$ to consider in $\schreier(H)$. However, it is not difficult to see that if we restrict the preorbits to cyclically reduced words, then one must only consider orbits in the restricted Stallings automaton $\rstallings(H)$, \emph{which is finite if $H$ is finitely generated}.

The final step is to see that one can overcome the gap between cyclically reduced to general preorbits by conjugating by $H$, which can still be done in an algorithmic way. Therefore subgroup preorders are computable in free groups. The precise statement is given below.

\begin{prop}[\citenr{delgado_relative_2022}] \label{prop: sgp preorder}
Let $H$ be a finitely generated subgroup of $\Free$, and let $k\geq 1$. Then, the set of $\Free$-elements of order $k$ in the subgroup $H$ is either empty or
 \begin{equation} \label{eq: order k = conj finite union}
\ofo{\Free}{k}{H} \,=\, \left( \bigsqcup\nolimits_{|\vvertii| = k} \big( \ofo{\Free}{\vvertii}{H} \big)^{\!v_{\vvertii}} \right)^{\!H} ,
 \end{equation}
where the union goes over the (finitely many) trails $\vvertii$ of length $k$ within the restricted Stallings automaton $\rstallings(H)$, and $v_{\vvertii}$ is the label of any walk from $\vertii_0$ to~$\bp$.  
\end{prop}

We  emphasize that it is possible to stretch the previous arguments and extend preorder computability to cosets (see~\parencite{delgado_relative_2022} for proofs and details). Accordingly, below we state the consequences of \Cref{prop: sgp preorder} (and its coset counterpart) in the most general way.  

\begin{cor}\label{cor: SBP Free}
Let $H$ be a finitely generated subgroup of $\Free[n]$. The spectrum of $\Free[n]$ in $H$ (resp., in a coset $Hu$) is bounded above by $\card \Verts\rstallings(H)$ (resp., by $\card \Verts\stallings(Hu)$). 
\end{cor}

\begin{rem}
According to this last result, $\card \ord_H(\Free[n])\leqslant 1+\card \Verts\rstallings(H)$. Since $\ord_H(\Free[n])$ is closed under taking divisors, it also makes sense to consider the set of maximal elements in $\ord_H(\Free[n])$ \wrt divisibility. In this direction, \citeauthor{linton_intersections_2021} obtains a bound for the cardinal of this set in terms of the rank of $H$ (see~\cite{linton_intersections_2021}). 
\end{rem}

Combining \Cref{cor: SBP Free} with the decidability of $\MP(\Free)$ (\Cref{thm: MP}) we obtain the computability of relative orders in $\Free$.

\begin{thm}[\citenr{delgado_relative_2022}]
There is an algorithm that, on input a finitely generated subgroup $H\leqslant \Free$, a coset $Hu$, and a word $w\in \Free$, outputs the relative order $\ord_{Hu}(w)$.
\end{thm}

On the other hand, it is straightforward to see that there exist elements of order $k = 0$ in $H \leqslant \Free$ (and compute them) if and only if the index $\ind{H}{\Free} = \infty$ (a condition again computable, see~\Cref{cor: FIP}).

\begin{thm}[\citenr{delgado_relative_2022}]
There is an algorithm that, on input a finitely generated subgroup $H\leqslant \Free$, an element $u\in \Free$, and an integer $k\geq 0$, it computes the preorder $\ofo{\Free}{k}{Hu}$.
\end{thm}

\begin{cor} \label{cor: k-roots of F computable}
Let  $H$ be a finitely generated subgroup of $\Free$. Then, for every $k \geq 0$, the set $\sqrt[k]{H}$ of $k$-roots of $H$ is computable; and hence the sets $\sqrt{H}$ and $\sqrt[*]{H}$ are computable as well.
\end{cor}

Finally, the algorithmic boundability of the spectra together with the computability of relative preorders provides the computability of spectra in free groups.

\begin{thm}[\citenr{delgado_relative_2022}] \label{thm: spectrum F computable}
There exists an algorithm which, on input a finitely generated subgroup $H\leq\fg \Free$, and an element $u\in \Free$, outputs the (finite) spectrum $\ord_{Hu}(\Free)$. 
\end{thm}

As a corollary we obtain the algorithmic decidability of subgroup purity, one of the first results in this direction.

\begin{cor}[\citenr{birget_pspace-complete_2000}; \citenr{kapovich_stallings_2002}]
There exists an algorithm to decide, given a finitely generated subgroup $H\leq\fg \Free$, whether it is pure.
\end{cor}
 
Note that if we restrict the results in this section relative to the trivial subgroup, we recover some classic facts about free groups.

\begin{cor}
Free groups are torsion-free.
\end{cor}

\begin{cor}
Elements in free groups have unique $p$-roots, if any ($p \geq 1$). That is, if $\Free$ is free, $u,v \in \Free$ and $u^p = v^p$ then $u=v$.
\end{cor}

\begin{cor}
Every element in a free group has only finitely many roots; \ie for every $w \in \Free$, $\sqrt{w}$ is finite.
\end{cor}

We finish this subsection with the computability of pure closures in free groups. It is important to realize that adding to $H$ all its proper roots is not enough to reach the pure closure $\pcl{H}$, because these additions may create new proper roots to be attached (consider \eg $H=\gen{a^2, ab^2}\leq \Free[2]$).
In order to guarantee that this process stabilizes after finitely many steps, we need to use the concept of algebraic extension (see \Cref{ssec: algebraic extensions}) through the lemma below.

\begin{lem}\label{lem: H algebraic}
A finitely generated subgroup $H$ of $\Free$ is always algebraic in $\gen{\sqrt{H}}$, which is again finitely generated.
\end{lem}

Now, the desired result follows easily from the transitivity of algebraic extensions (see~\Cref{ssec: algebraic extensions}).

\begin{prop}\label{prop: closure estabilization}
Let $H$ be a finitely generated subgroup of $\Free$, and for $i\geq 1$, let $H_{i+1}=\gen{\sqrt{H_i} }$, with $H_0=H$. Then, all the $H_i$'s are finitely generated, and the ascending sequence $(H_i)_{i\geq 0}$ stabilizes at~$\pcl{H}$, which is therefore finitely generated and computable. 
\end{prop}

\subsection{\bf Algebraic extensions and Takahasi's theorem} \label{ssec: algebraic extensions}

In \Cref{ssec: hom} we defined the notion of homomorphism between $A$-automata. The Stallings bijection in \eqref{eq: Stallings bijection} behaves well with respect to inclusions in the sense that two subgroups are contained one into the other, $H\leqslant K\leqslant \Free[\Alfi]$, if and only if there exists an $A$-homomorphism from $\stallings(H)$ to $\stallings(K)$ which, in case it exists, is unique and will be denoted by $\theta_{H,K}\colon \stallings(H) \to \stallings(K)$. Further, we have already seen that, in the special case when $\theta_{H,K}$ is injective, \ie when $\stallings(H)$ is a subautomaton~of $\stallings(K)$, then $H$ is a free factor of $K$; see \Cref{cor: incl => ff}. 

To dualize the notion of free factor, \citeauthor{kapovich_stallings_2002} introduced in \parencite{kapovich_stallings_2002} the notion of algebraic extension: a subgroup extension $H\leqslant K\leqslant \Free[\Alfi]$ is called \emph{algebraic}, denoted by $H\leqalg K$, if $H$ is not contained in any proper free factor of $K$; we denote by $\AAEE(H)$ the set of algebraic extensions of $H$ (within $\Free[\Alfi]$). It is straightforward to see transitivity, \ie $H\leqalg K\leqalg L$ implies ${H\leqalg L}$. 

The notions of free and algebraic extension have a close relation, which is nice to compare with a similar behavior happening in other algebraic structures. Free factors are the non-abelian analog of direct summands from commutative algebra. In a vector space, every basis of a subspace $E$ can be extended to a basis of any subspace $F \geqslant E$, i.e., $F=E\oplus E'$ for some complementary subspace $E'$. When we consider, for example, free abelian groups (i.e., free modules over $\mathbb{Z}$) the exact same result is not true,
but still every submodule $H$ `is close' to a direct summand of any submodule $K$ containing $H$:
 every subgroup $H\leqslant \mathbb{Z}^m$ is of finite index only in finitely many subgroups $H=H_0\leqfi H_1, \ldots ,H_n\leqslant \mathbb{Z}^m$ and, for every $K\leqslant \mathbb{Z}^m$ containing $H$, there exists a unique $i$ such that $H\leqfi H_i\leqslant_{\oplus} K\leqslant \mathbb{Z}^m$. Of course, the situation in the free group seems much wilder, starting from the well known fact that $H\leqslant K\leqslant \Free[\Alfi]$ does not even imply $\rk(H)\leqslant \rk(K)$ (to the extreme that $\Free[\hspace{1pt}{\aleph_0}]$ can be viewed as a subgroup of $\Free[\hspace{1pt}2]$, see~\Cref{cor: ranks of subgroups}). However, back in the 1950's, Takahasi~\parencite{takahasi_note_1951} proved that, again, the same result can be adapted to the free group case, after admitting a little bit more of degeneration: we have to restrict ourselves to finitely generated subgroups, and we lose the finite index condition. In modern language it can be stated as follows.
\begin{thm}[\citenr{takahasi_note_1951}]\label{thm: Tak}
Let $H\leqfg \Free[\Alfi]$. Then $\AAEE(H)$ is a finite and computable collection of finitely generated subgroups of $\Free[\Alfi]$; furthermore, every $K\leqslant \Free[\Alfi]$ containing $H$ algorithmically determines a unique $H'\in \AAEE(H)$ such that $H\leqalg H'\leqff K\leqslant \Free[\Alfi]$ (called the \defin[algebraic closure]{$K$-algebraic closure of~$H$}).
\end{thm}
Takahasi's Theorem was proved in the 1950's, in a different language, and using purely combinatorial and algebraic techniques. However, in more recent years, it was rediscovered independently, by Ventura~\parencite{ventura_fixed_1997} in 1997, by Margolis--Sapir--Weil~\parencite{margolis_closed_2001} in 2001, and by Kapovich--Miasnikov~\parencite{kapovich_stallings_2002} in 2002, in slightly different contexts; see also the subsequent paper \parencite{miasnikov_algebraic_2007} by~\citeauthor{miasnikov_algebraic_2007} joining the three points of view. These authors, independently, gave their own proofs of Takahasi's theorem, which happened to be essentially the same; we would say, the natural geometric proof of this result using Stallings automata. Let us sketch it below. 

Given $H\leqfg \Free[\Alfi]$, consider its (finite) Stallings automaton $\stallings(H)$. Identifying certain sets of vertices together, we obtain a new $A$-automaton which may very well be not deterministic; apply then a sequence of foldings until reaching a deterministic one, say $\Ati_{\! 1}$, which will correspond to a finitely generated extension of $H$, say~$H\leqslant H_1=\gen{\Ati_{\! 1}}\leqslant \Free[\Alfi]$. Repeating this operation for each of the \emph{finitely many} partitions of the set of vertices $\Verts\stallings(H)$, we obtain a finite list (with possible repetitions) of finitely generated extensions of $H$, say~$\OO_A (H)=\{H_0, H_1, \ldots ,H_n\}$, called the \defin{$A$-fringe of $H$} in~\parencite{ventura_fixed_1997} (note that $H_0=H$ is always present in the list, corresponding to trivial partition). Observe that the fringe satisfies the following property. Let $K$ be a (not necessarily finitely generated) subgroup with $H\leqslant K\leqslant \Free[\Alfi]$, and consider the corresponding $A$-homomorphism $\theta_{H,K}\colon \stallings(H)\to \stallings(K)$. Looking at the image $\operatorname{Im} (\theta_{H,K})$ as a (finite and deterministic) subautomaton of $\stallings(K)$, we see that: (1) the $A$-homomorphism $\theta_{H,K}\colon \stallings(H)\onto \operatorname{Im}(\theta_{H,K})$ is onto and so, $\gen{\operatorname{Im}(\theta_{H,K})}\in \OO_A(H)$; and (2) $\operatorname{Im}(\theta_{H,K})$ is a subautomaton of $\stallings(K)$ and so, $\gen{\operatorname{Im}(\theta_{H,K})}\leqff K$. This shows  that $\AAEE(H)\subseteq \OO_A(H)$. Finally, as done in~\parencite{kapovich_stallings_2002,miasnikov_algebraic_2007}, it is possible to algorithmically clean up $\OO_A(H)$ in order to gain uniqueness of the middle subgroup $H'$ in \Cref{thm: Tak}: for each pair of distinct subgroups $H_i, H_j\in \OO_A(H)$, if $H_i\leqff H_j$ then remove $H_j$ from the list.

For this last cleaning process to be computable, one needs an algorithm deciding whether a given extension $H\leqslant K$ is free or not; this can be done using classical Whitehead techniques or, alternatively, using more modern algorithms based on Stallings automata.

\begin{thm}[Whitehead, \parencite{lyndon_combinatorial_2001}; Silva--Weil, \parencite{silva_algorithm_2008}; Puder, \parencite{puder_primitive_2014}]
There is an algorithm which, on input an extension of finitely generated subgroups $H\leqfg K\leqfg \Free[\Alfi]$, decides whether $H$ is a free factor of $K$.
\end{thm}

We remark, however, that these three algorithms have exponential complexity on the sum of lengths of the given generators for $H$ and $K$. The best known improvement is a variation on the Whitehead algorithm given by Roig--Ventura--Weil~\parencite{roig_complexity_2007} and working in polynomial time with respect to that sum of lengths.  

The modern version of Takahasi's Theorem in \Cref{thm: Tak} has several applications to the study of the lattice of subgroups of a free group; we survey some of them below. However, before that, we mention a recent development originated by the above graphical proof of \Cref{thm: Tak} (see \parencite{ventura_onto_2021} for the details on this story).

\subsubsection{\bf Onto and fully-onto extensions}

The fringe of $H$ strongly depends on the ambient basis $A$ (reflected in the notation with the subscript $A$ in $\OO_A(H)$), whereas the set of algebraic extensions $\AAEE(H)$ does not, and is canonically associated to the subgroup $H$, since it is defined completely in algebraic terms. To illustrate this fact, see Example~2.5 from Miasnikov--Ventura--Weil~\parencite{miasnikov_algebraic_2007}, where the fringe of $H=\langle ab, acba\rangle \leqslant \Free[\Alfi]$, with $A=\{a,b,c\}$, is computed: $\OO_A(H)=\{ H_0, H_1, H_2, H_3, H_4, H_5\}$, where $H_0=H$, $H_1=\langle ab, ac, ba \rangle$, $H_2=\langle ba, ba^{-1}, cb\rangle$, $H_3=\langle ab, ac, ab^{-1}, a^2 \rangle$, $H_4=\langle ab, aca, acba\rangle$, and $H_5=\langle a, b, c\rangle =\Free[\Alfi]$; however, with respect to the new ambient basis $A\!'=\{ d, e, f\}$, where $d=a$, $e=ab$, and $f=acba$, the $A\!'$-automaton $\Ati_{\!A\!'}(H)$ has a single vertex, and hence the $A\!'$-fringe of $H$ is much simpler, $\OO_{A\!'}(H)=\{H\}$. Of course, in this example, $H\leqff \Free[\Alfi]$ and $\AAEE(H)=\{H\}$. 

One can interpret this fact by thinking that $\AAEE(H)$ is what really carries relevant algebraic information about the subgroup $H$ and its relative position within the lattice of subgroups of $\Free[\Alfi]$; and $\OO_A(H)$ is the same set plus some accidental new members depending on the ambient basis used to draw and work with the graphs. From this point of view, Miasnikov--Ventura--Weil launched in~\parencite{miasnikov_algebraic_2007} a natural conjecture stating that the common subgroups in $\OO_{A\!'}(H)$, when $A\!'$ runs over all ambient bases might be, precisely, the algebraic extensions: for every $H\leqfg \Free[\Alfi]$, 
$\AAEE(H) =\bigcap_{A\!'\,\, \mbox{\scriptsize basis of } \Free[\Alfi]} \OO_{A\!'}(H)$. Somehow contraintuitively, seven years later, in 2014, Parzanchevski--Puder found the counterexample below. 

\begin{prop}[{\citenr{parzanchevski_stallings_2014}}]\label{pp}
Let $A=\{a,b\}$. In $\Free[\Alfi]$, the proper subgroup extension $H=\langle a^2b^2\rangle\leqslant \langle a^2b^2, ab\rangle =K\leqslant \Free[\Alfi]$ is free, $H\leqff K$ (so, it is not algebraic), but it satisfies $K\in \OO_{A\!'}(H)$ for every ambient basis $A\!'$.
\end{prop}

They proposed a couple of possible reformulations for the conjecture, making it more plausible. However, both were disproved again by the following stronger counterexample.

\begin{prop}[{\citenr{kolodner_algebraic_2021}}]\label{kolodner}
Let $A=\{a,b\}$. In $\Free[\Alfi]$, the proper subgroup extension $H=\langle b^2aba^{-1}\rangle \leqff \langle b, aba^{-1}\rangle=K$ is free (and so, not algebraic) but, for an arbitrary alphabet $B$, and for every homomorphism $\varphi\colon \Free[\Alfi]\to \Free[B]$ with $a\varphi, b\varphi \neq 1$, $K\varphi \in \OO_{B}(H\varphi)$.
\end{prop}

See~\parencite{ventura_onto_2021} for another reformulation of the conjecture, this one being true. These counterexamples gave rise to the notion of onto and fully-onto extensions among subgroups of $\Free[\Alfi]$, introduced in~\parencite{ventura_onto_2021}.

\begin{defn}
Let $H\leqslant K\leqslant \Free[\Alfi]$. We say that this is an \defin[onto extension]{onto} extension of subgroups, denoted by $H\leqont K$, if $\theta_{H,K}\colon \stallings(H,A\!')\to \stallings(K,A\!')$ is onto, for every basis $A\!'$ of $\Free[\Alfi]$; in other words, if $K\in \bigcap_{A\!'} \OO_{A\!'}(H)$, where $A\!'$ runs over all possible basis for $\Free[\Alfi]$. Further, we say that $H\leqslant K$ is \defin[fully onto extension]{fully onto}, denoted by $H\leqfont K$, if $\theta_{H,K}\colon \stallings(H,B')\to \stallings(K,B')$ is onto, for every basis $B'$ of every free extension $\Free[\Alfi]\leqff \Free[B]$, where $B\supseteq A$. 
\end{defn}

The relation between these concepts is expressed in the following result. See~\parencite{ventura_onto_2021} for more properties of onto and fully-onto extensions, and of the corresponding closure operators (which do not agree, in general, with the algebraic closure). 

\begin{prop}[\citenr{ventura_onto_2021}]
Let $H\leqslant K\leqslant \Free[\Alfi]$. The following implications hold, while the reverse implications are not true in general: 
 \begin{equation}
H\leqalg K \, \Rightarrow \, H\leqfont K \, \Rightarrow \, H\leqont K \, \Rightarrow \, K\in \OO_A(H).
 \end{equation}
\end{prop}

\subsubsection{\bf Computation of pro-${\mathcal V}$ closures}

A classical interplay between Group Theory and Topology is the construction and study of the so-called \emph{profinite topology} in a given arbitrary group $G$. This can be generalized to any pseudo-variety of finite groups as follows. 

A \defin{pseudo-variety} ${\mathcal V}$ of finite groups is a family of (isomorphism classes of) finite groups closed under taking subgroups, quotients, and finite direct products (the name comes from the notion of \defin{variety} of groups, being the same but not restricted to finite groups, and allowing also infinite direct products). For instance, all finite groups, the family of $p$-groups for a given prime $p$, finite nilpotent groups, finite solvable groups or finite abelian groups, etc. are typical examples of pseudo-varieties of groups. Given such a pseudo-variety ${\mathcal V}$, one can define (metrically) the \defin{pro-${\mathcal V}$ topology} in an arbitrary group $G$ in the following way: given two elements $g,g'\in G$ define the \defin{${\mathcal V}$-distance} between them as $d_{\mathcal V} (g,g')=2^{-v(g,g')}$, where $v(g,g')$ is the smallest cardinal of a group $H\in {\mathcal V}$ for which there is a homomorphism $\varphi \colon G\to H$ separating $g$ and $g'$, i.e., such that $g\varphi \neq g'\varphi$ (take $d_{\mathcal V} (g,g')=0$ if $v(g,g')=\infty$ meaning that there is no such finite group $H\in {\mathcal V}$). It is easy to see that $d_{\mathcal V}$ is a pseudo-metric in $G$, which induces a topology called the pro-${\mathcal V}$ topology. In case the group $G$ is \defin{residually-${\mathcal V}$} (for every two distinct elements $g,g'\in G$ there exists $H\in {\mathcal V}$ and a homomorphism $\varphi \colon G\to H$ separating $g$ and $g'$), the defined pseudo-metric is then a real metric, and the induced topology on $G$ becomes Hausdorff. The above examples of pseudo-varieties give rise to the so-called \defin{pro-finite topology}, the \defin{pro-$p$ topology}, the \defin{pro-nilpotent topology}, the \defin{pro-solvable topology}, the \defin{pro-abelian topology}, etc. in any arbitrary group $G$. It is easy to see that the ${\mathcal V}$-closure of a subgroup $H\leqslant G$ is again a subgroup, $H\leqslant \Cl_{\mathcal V}(H)\leqslant G$; this observation opens the door to nice questions about algebraic or algorithmic properties of these closure operators.  

Particularizing to free groups, and restricting the attention to extension-closed pseudo-varieties~${\mathcal V}$ (${\mathcal V}$ is called \defin[extension-closed $\mathcal V$]{extension-closed} if, for any short exact sequence $1\to A\to B\to C\to 1$ of finite groups, $A,C\in {\mathcal V}$ imply $B\in {\mathcal V}$), Margolis--Sapir--Weil proved the result below. 

\begin{thm}[{\citenr{margolis_closed_2001}}]\label{MSW1}
Let ${\mathcal V}$ be an extension-closed pseudo-variety of finite groups, and consider the pro-${\mathcal V}$ topology on the free group~$\Fn$. Any free factor of a ${\mathcal V}$-closed subgroup of $\Fn$ is again ${\mathcal V}$-closed.
\end{thm} 

This automatically connects with Takahasi's Theorem because it implies that, in the extension-closed case, the pro-${\mathcal V}$ closure of any subgroup $H\leqslant_{fg} \Fn$ is an algebraic extension of $H$ itself; that is, $H\leqalg \Cl_{\mathcal V}(H)$. Using this idea the authors of \parencite{margolis_closed_2001} gave nice algorithms based on Stallings automata to compute several closures.

\begin{thm}[{\citenr{margolis_closed_2001}}]
There is an algorithm which, on input a finite set of generators for $H\leqfg \Fn$, computes a basis for its pro-finite, its pro-$p$, and its pro-nilpotent closures. 
\end{thm}

The computability of the pro-solvable closure is still an open problem. Being an extension-closed pseudo-variety of finite groups, \Cref{MSW1} is valid for the solvable case and so, $H\leqalg \Cl_{\operatorname{sol}}(H)$ for any $H\leqfg \Fn$; however, no criterion is known yet to distinguish, among the computable candidates $\AAEE(H)$, which one is the solvable-closure of $H$ (it would suffice to find an algorithm to decide whether a given $H\leqfg \Fn$ is solvable-closed). 

\subsubsection{\bf Fixed subgroups}

One of the first modern applications of Takahasi's Theorem was about fixed subgroups of endomorphisms of free groups (see~\parencite{ventura_fixed_1997}). This is one out of many results in a rich line of research that we quickly survey in this subsection; see~\parencite{ventura_fixed_2002} for more details. 

The \defin{fixed subgroup} of an endomorphism $\varphi \colon G\to G$ of an arbitrary group $G$ is  $\Fix(\varphi)=\{ g\in G \st g\varphi=g\}\leqslant G$. For the free case $G=\Fn$, Dyer and Scott proved in~\parencite{dyer_periodic_1975} that, when $\varphi$ is an automorphism of finite order, $\Fix(\varphi)$ is a free factor of the ambient $\Fn$. After checking that this is not the case in general (for example, the fixed subgroup of $\varphi\colon \Free_2\to \Free_2$, $a\mapsto a$, $b\mapsto ab$ is $\Fix(\varphi)=\langle a, b^{-1}ab\rangle$), they conjectured that any automorphism $\varphi$ of $\Fn$ should at least satisfy the inequality $\rk(\Fix(\varphi))\leqslant n$. This opened a very fruitful line of research for several decades, proving particular cases of the conjecture, until the celebrated result by Bestvina--Handel solving it completely (and opening an even more active line for the following years, giving rise to what is today known as the \emph{train track} theory for graphs).

\begin{thm}[{\citenr{bestvina_train_1992}}]
For every $\varphi \in \Aut(\Fn)$, we have that $\rk(\Fix(\varphi))\leqslant n$.
\end{thm}

Far from exhausting the research activity in this direction, this result stimulated more research around fixed point subgroups. Already before its publication, Imrich--Turner~\parencite{imrich_endomorphisms_1989}, using a short and elegant algebraic argument, proved the same result for endomorphisms (by reducing it to the automorphism case). And Collins--Turner in~\parencite{collins_all_1996} gave an explicit description of the fixed subgroups of maximal rank, \ie with $\rk(\Fix \varphi)=n$. Using this description, Ventura~\parencite{ventura_fixed_1997} showed the following result which, as far as we know, is the first modern application of Takahasi's Theorem.

\begin{thm}[{\citenr{ventura_fixed_1997}}]
Let $\varphi$ be an automorphism of $\Fn$ with $H=\Fix(\varphi)$ of maximal rank, and let $H<K\leqslant \Fn$. Then either $\rk(H)<\rk(K)$ or, $\rk(H)=\rk(K)$ and~$\rk(H\rho)<\rk(K\rho)$, where $\rho\colon \Fn \twoheadrightarrow \ZZ^n$ is the ambient abelianization. 
\end{thm}

As a consequence, the following two corollaries were obtained.

\begin{cor}[{\citenr{ventura_fixed_1997}}]
Among the strictly ascending chains of maximal rank fixed subgroups of~$\Fn$, the maximum length is exactly $n$.
\end{cor}

\begin{cor}[{\citenr{ventura_fixed_1997}}]
Let $\Free_2$ be the free group of rank two, and $S$ a non-empty set of non-identity endomorphisms of $\Free_2$ such that $\rk(\Fix S)=2$. Then, $S\subseteq \Aut(\Free_2)$, and $\Fix(S)=\Fix(\varphi)=\langle a, b^{-1}ab\rangle$, for some $\varphi\in S$ and some basis $\{a, b\}$ of $\Free_2$. 
\end{cor}

Here, by the fixed subgroup of a family of endomorphisms we mean the set of points fixed by each member of the family, $\Fix(S)=\bigcap_{\varphi\in S} \Fix (\varphi)$. This last result solved the rank 2 case of a conjecture which is still open today: for an arbitrary family $S$ of endomorphisms of $\Fn$, there exists $\varphi\in \gen{S}$ (the submonoid of $\End(\Fn)$ generated by $S$) such that $\Fix(S)=\Fix(\varphi)$. The intuition behind this conjecture is clear: given $\varphi,\psi\colon \Fn \to \Fn$, it is obvious that any word on them (and their inverses if they are automorphisms),  $w(\varphi, \psi)$, satisfies the inclusion $\Fix (\varphi)\cap \Fix (\psi)\leqslant \Fix (w(\varphi, \psi))$. The opposite inclusion is not true in general because it could be the case that $\varphi$ moves $x$ to $y\neq x$, and $\psi$ moves back $y$ to~$x$, with then $x$ belonging to $\Fix (\varphi\psi)$ but not to $\Fix (\varphi)\cap \Fix (\psi)$; the conjecture proposes that one should be able to avoid all these accidental situations by using a complicated enough word $w(\varphi, \psi)$. 

This conjecture has been resolved in the positive for free groups of rank 3 by A. Martino in~\parencite{martino_intersections_2004}, but it remains still open in general, as far as we know. One of the few results obtained in this direction is, again, a consequence of Takahasi's Theorem.

\begin{thm}[{\citenr{martino_automorphism-fixed_2000}}]
Let $\Fn$ be the free group of rank $n$, and $S$ be a set of endomorphisms of $\Fn$. There exists $\varphi\in \gen{S}$ (the submonoid of $\End(\Fn)$ generated by $S$) such that $\Fix(S)$ is a free factor of $\Fix(\varphi)$. 
\end{thm}

We can briefly explain the core idea in the proof of this last result, in order to highlight the role played by Takahasi's Theorem. Assume $\varphi, \psi\in \Aut(\Fn)$ and assume also that $\Fix(\psi^r)=\Fix(\psi)$ for every $r\geqslant 1$ (technical arguments can be made to justify this assumption from the general case). Then, consider the inclusion $\Fix(\varphi)\cap \Fix(\psi)\leqslant \Fix(\varphi\psi^r)$, for every $r\in \NN$. Since $\Fix(\varphi)\cap \Fix(\psi)$ is finitely generated, Takahasi's Theorem tells us that, for every $r$, there exists $M_r\in \AAEE(\Fix(\varphi)\cap \Fix(\psi))$ such that $\Fix(\varphi)\cap \Fix(\psi)\leqalg M_r \leqff \Fix(\varphi\psi^r)$. But there are infinitely many natural numbers $r$, and only finitely many algebraic extensions of $\Fix(\varphi)\cap \Fix(\psi)$ so, we must have $M_r=M_s$ for some $r<s$ (in fact, there must be some $M_r$ repeated infinitely many times). From this equality we deduce that 
 \begin{equation*}
\Fix(\varphi)\cap \Fix(\psi) \,\leqslant\, M_r\cap M_s \,\leqslant\, \Fix(\varphi\psi^r)\cap \Fix(\varphi\psi^s) \,\leqslant\, \Fix(\varphi) \cap \Fix(\psi)
 \end{equation*}
and so, all the inclusions must be equalities; this implies $\Fix(\varphi)\cap \Fix(\psi)=M_r\leqff \Fix(\varphi\psi^r)$, as we wanted to see. 

Finally, let us mention yet another nice application of Takahasi's Theorem. Dualizing the concept of fixed subgroup, we can define the \defin{auto-stabilizer} (resp., \defin{endo-stabilizer}) of a subgroup $H\leqslant G$ as $\aStab(H)=\{ \varphi\in \Aut(G) \st H\leqslant \Fix(\varphi) \}$ (resp., $\eStab(H)=\{ \varphi\in \End(G) \st H\leqslant \Fix(\varphi) \}$), a subgroup of $\Aut(G)$ (resp., submonoid of $\End(G)$). Also, we can define the \defin{auto-closure} (resp., \defin{endo-closure}) of~$H$ as $\aCl(H)=\Fix \big( \aStab(H)\big)$ (resp., $\eCl(H)=\Fix \big( \eStab(H)\big)$, the smallest subgroup fixed by all automorphisms (resp., endomorphisms) of $G$ fixing $H$. As an elementary example note that, in the free group $\Fn$, the equation $X^2=a^2$ has only one solution, namely $X=a$; this means that no endomorphism of $\Fn$ can fix~$a^2$ moving $a$ or, in other words, $\eCl(\langle a^2\rangle)=\aCl(\langle a^2\rangle)=\langle a\rangle$. In this setting, natural questions arise about finite generation and computability of these stabilizers and closures. Some of them are classical results, and some others were solved more recently, again as an application of Takahasi's Theorem.  

\begin{thm}[\citenr{mccool_finitely_1975}]
Let $H\leqfg \Fn$ be given by a finite set of generators. Then the auto-stabilizer, $\aStab(H)$, is finitely presented, and a finite presentation can be algorithmically computed.
\end{thm}

\begin{thm}[{\citenr{bogopolski_algorithm_2016}}; {\citenr{feighn_algorithmic_2018}}]
Let
$\varphi$ be an automorphism of a free group~$\Free[n]$.
Then a basis for $\Fix(\varphi)$ is computable. 
\end{thm}

As a consequence of these two results one can compute auto-closures of finitely generated subgroups of $\Free[n]$: given $H\leq \Free[n]$ (by a finite set of generators) compute finitely many automorphisms $\varphi_1, \ldots ,\varphi_m\in \Aut(\Free[n])$ such that $\aStab(H)=\gen{\varphi_1, \ldots ,\varphi_m}$, compute bases for $\Fix(\varphi_1),\ldots,\allowbreak \Fix(\varphi_m)$, and intersect them all with the pull-back technique; clearly,
$
\aCl(H)=\bigcap\nolimits_{\varphi\in \aStab(H)} \Fix{\varphi} =\Fix(\varphi_1)\cap \cdots \cap \Fix(\varphi_m)
$.

\begin{prop}[\citenr{ventura_computing_2010}]
If $H\leqslant \Free[n]$ is finitely generated, then the auto-closure $\aCl(H)$ is again finitely generated, and a basis is computable from a given set of generators for $H$.
\end{prop}

The corresponding problem for endomorphisms is more interesting and tricky. The strategy used in the automorphism case fails because, even when $H$ is finitely generated as a subgroup, $\eStab(H)$ needs not be finitely generated as submonoid (see \parencite{ciobanu_two_2006} for an explicit counterexample). Despite this obstruction, the problem was solved in \parencite{ventura_computing_2010}, again with an argument where Takahasi's Theorem played a crucial role. 

\begin{thm}[\citenr{ventura_computing_2010}]\label{thm: endoclosure}
Let $H\leqfg \Free[n]$ be a finitely generated subgroup of a free group, given by a finite set of generators. Then, the endo-closure $\eCl(H)$ of $H$ is again finitely generated, and a basis is algorithmically computable, together with a set of $m\leqslant 2n$ endomorphisms $\varphi_1,\ldots ,\varphi_m\in \End(\Free[n])$, such that $\eCl(H) = \Fix(\varphi_1)\cap \cdots \cap \Fix(\varphi_m)$.
\end{thm}

The relation with Takahasi's theorem is the following: compute the algebraic extensions of $H$ which are retracts, say $\AAEE_{ret}(H)=\{H_1, \ldots ,H_r\}$, and consider the auto-closure of $H$ as a subgroup of each $H_i$, \ie $\aCl^{H_i}(H)$. A technical argument using Takahasi's Theorem shows that $\eCl(H)=\aCl^{H_1}(H)\cap \cdots \cap \aCl^{H_r}(H)$. This completes the proof. 

\subsubsection{\bf Inertia and compression}

It is worth mentioning the first (of several) refinements and improvements of Bestvina--Handel development, giving rise to the following stronger result. Recall that in a free group, the intersection of two subgroups $H\cap K$ can have a greater rank than those of $H$ and~$K$.

\begin{thm}[\citenr{dicks_group_1996}]\label{thm: DV}
For every injective endomorphism $\varphi \colon \Fn \to \Fn$, $\Fix(\varphi)$ is inert in $\Fn$; in particular, $\rk(\bigcap_{\varphi\in S} \Fix (\varphi))\leqslant n$, for any family $S$ of injective endomorphisms of $\Fn$.
\end{thm}

Recall that a subgroup is \defin[inert subgroup]{inert} if $\rk(H\cap K)\leqslant \rk(K)$, for every $K\leqslant G$. Observe that if $H$ is inert then $\rk(H)=\rk(H\cap G)\leqslant \rk(G)$ so, the case $\card S=1$ of \Cref{thm: DV} is already an improvement of the 
Bestvina--Handel theorem. It is immediate from the definition that the intersection of two inert subgroups is inert as well (and, with a technical argument about chains of subgroups, the same is true for infinite intersections in the case $G=\Fn$).  

The same fact was conjectured to be true as well for general sets of endomorphisms. After several partial results (see \parencite{martino_fixed_2004,zhang_fixed_2015}) this conjecture has recently been proved (and extended to surface groups as well).

\begin{thm}[\citenr{antolin_hanna_2021}]
Let $G$ be a free or surface group. Then, for every endomorphism $\varphi \colon G\to G$, $\Fix(\varphi)$ is inert in $G$; in particular, $\rk(\bigcap_{\varphi\in S} \Fix (\varphi))\leqslant n$, for any family $S$ of endomorphisms of~$G$.
\end{thm}

With the idea of quantifying how far is a subgroup from being inert, the paper \parencite{roy_degrees_2021} introduced the concept of degree of inertia: in a group $G$, the \defin{degree of inertia} of a subgroup $H\leqslant G$ is the supremum 
 \begin{equation*}
\di_G(H)=\sup_{K\leqfg G} \frac{\rrk(H\cap K)}{\rrk(K)},
 \end{equation*}
where the supremum is taken over all finitely generated subgroups ${K\leqslant G}$ (and~$0/0$ is understood to be 1); in general, $\di_G(H)\geqslant 1$, with equality if and only if $H$ is inert. In \parencite{roy_degrees_2021}, using the Stallings machinery, some technical results were developed to help understand this notion, both in free groups $\Free[n]$, and in groups of the form $\Free[n] \times \mathbb{Z}^m$. In general, it is not known whether the degree of inertia is computable in terms of a given set of generators for $H$, even in the free ambient.

It is worth noting that, in the free group case, it is not even known whether the supremum in the definition of degree of inertia is a maximum, \ie whether it is always achieved in a particular subgroup $K$. Interestingly, \citeauthor{ivanov_intersection_2018} showed in  \parencite{ivanov_intersection_2018} that this is the case if we replace the numerator $\rrk(H\cap K)$ by the corresponding disconnected version $\sum_{x\in H \backslash \Free / K} \rrk (\Sgpi \cap x^{-1} \Sgpii x)$ (in \parencite{ivanov_intersection_2018}, this related notion is called the \defin{Walter Neumann coefficient} of $H$, in clear connection with the Strengthened Hanna Neumann inequality). This result is proved using quite unusual techniques: the author codifies the Stallings automata of all the subgroups $K$ that must be taken into account in the supremum as points in some abstract simplex, and then shows that the function to be maximized takes its maximum at some point in the boundary by using linear programming techniques. A crucial point here is that, with the disconnected version in the numerator, the arguments are more symmetric because the counting is directly related to the whole pullback $\stallings(H)\times \stallings(K)$, instead of with the particular connected component corresponding to $H\cap K$. It is not clear whether one can overcome this technical point and prove a similar result for the degree of inertia.  

A relaxed version of the notion of inertia is that of compression: a subgroup $H\leqslant G$ is  \defin[compressed subgroup]{compressed} if $\rk(H)\leqslant \rk(K)$, for every $H\leqslant K\leqslant G$. Note that inert subgroups are compressed, whereas the converse is not true in general. Using Takahasi's \Cref{thm: Tak} it is easy to see that, for $G$ free, compression can be algorithmically decided; deciding inertia is an open problem.

\begin{prop}
There is an algorithm which, given $u_1,\ldots ,u_k\in \Free$, it decides whether the subgroup $\gen{u_1,\ldots ,u_k}$ is compressed.  
\end{prop}

On the other hand, directly from the definition, if $H_1$ and $H_2$ are inert then $H_1\cap H_2$ is again inert. The corresponding result for compression is not known (compression of $H_1$ and $H_2$ does not seem to contain enough information about the eventual subgroups $K$ satisfying ${H_1\cap H_2 \leqslant K}$, but ${H_1\not\leqslant K}$ and $H_2 \not\leqslant K$).

Despite this seemingly contrasting behaviour between the notions of compression and inertia, no single example is known of a subgroup of a free group $H\leq \Free$ which is compressed but not inert. This brought some authors to ask whether compressed implies inert in a free ambient.  

Finally, following the same motivation, the notion of degree of compression was also introduced in \parencite{roy_degrees_2021}: the \defin{degree of compression} of a subgroup $H\leqslant G$ is the supremum 
\begin{equation*}
\dc_G(H)\,=\,\max_{H\leqslant K\leqfg G} \frac{\rrk(H)}{\rrk(K)} \,=\,\frac{\rrk(H)}{\min_{H\leqslant K\leqfg G} \rrk(K)} \,,
\end{equation*}
where $0/0$ is understood to be $1$; here, it is clear that the supremum is always a maximum; and again $\dc_G(H)\geqslant 1$ with equality if and only if $H$ is compressed. A straightforward application of Takahasi's \Cref{thm: Tak} shows that the degree of compression is computable in free groups:

\begin{prop}
There is an algorithm which, given $u_1,\ldots ,u_k\in \Free$, it computes the degree of compression $\dc_{\Free}(\gen{u_1,\ldots ,u_k})$. 
\end{prop}

\subsection{\bf Asymptotic behavior}\label{sec: asymptotic}
The asymptotic behavior of infinite groups has aroused increasing interest in the last decades,
especially after the famous announcement by Gromov that `generically'
all finitely presented groups are hyperbolic (see~\parencite{gromov_hyperbolic_1987}).

In a finite environment, there is a natural default meaning for this kind of claims; namely the one given by the uniform distribution, which reduces probabilistic questions to essentially counting proportions among interesting or meaningful subsets. In order to analyze an infinite population, the first serious obstacle is that one needs a rigorous way to compare the size of infinite subsets. To estimate `how frequent' a given property of subgroups of $\Free[r]$ is, or to compute, for example, their average rank, one has to previously specify in which sense those claims should be understood, \ie we need to specify a
measure on the set of subgroups of $\Free[r]$; this is what is called a \defin{probabilistic model}. It has been customary to call \defin[generic property]{generic} a property which happens with probability one \emph{in some well-defined model}. Therefore, \emph{the meaning of the word `generic' may differ between statements} and must be clarified in each case.

An standard way to extend probabilistic claims to an infinite population~$\Omega$ is to stratify using finite layers; that is, to distinguish an ascending sequence ${\Omega_0 \subseteq \Omega_1 \subseteq \cdots}$ of finite subsets such that $\bigcup_{i=0}^{\infty} \Omega_i =\Omega$, and (try to) define the probabilistic notion on $\Omega$ as the limit (\emph{if it exists}) of the corresponding notions in the finite subsets~$\Omega_i$. Then, given a property (subset) $\mathcal{P} \subseteq \Omega$ of the objects under study, we say that $\mathcal{P}$ is \defin[generic property]{generic} (\resp \defin[negligible property]{negligible}) \emph{\wrt the used model} if $\lim_{n\to \infty} \card \,(\Omega_n \cap \mathcal{P})/\card\Omega_n = 1$ (\resp $\lim_{n\to \infty} \card \,(\Omega_n \cap \mathcal{P})/\card\Omega_n = 0$). If $\mathcal{F}$ is a class of functions tending to $0$ and closed under $\max$, we  say that $\mathcal{P}$ is \defin[$\mathcal{F}$-negligibility]{$\mathcal{F}$-negligible} if $\card \,(\Omega_n \cap \mathcal{P})/\card\Omega_n = \mathcal{O}(f(n))$ for some $f \in \mathcal{F}$, and \defin[$\mathcal{F}$-genericity]{$\mathcal{F}$-generic} if $\Omega \setmin \mathcal{P}$ is $\mathcal{F}$-negligible. Most of the interest so far is in \defin{exponential genericity}.

For example, in~\parencite{hall_jr_coset_1949}, \citeauthor{hall_jr_coset_1949} proved the nice result below. %

\begin{thm}[\citenr{hall_jr_subgroups_1949}] \label{thm: Hall tally}
The number $N(k,n)$ of subgroups of index $k$ in~$\Free[n]$ satisfies the following recursive formula: 
 \begin{equation}
 N(k,n) \,=\, k(k!)^{n-1}-\sum\nolimits_{i=1}^{k-1}((k-i)!)^{n-1}N(i,n) \,.
 \end{equation}
\end{thm}

Formulas of this kind open the door to the study of asymptotics on subgroups using the previous approach. One could try to count, for example, how many of such subgroups (of index $k$) are normal: if $\alpha(k)$ denotes the proportion of normal subgroups among subgroups of $\Free[n]$ of index $k$, the limit $\lim_{k\to \infty} \alpha(k)$ (if it exists!) provides a precise \emph{and specific} quantitative measure of the frequency of the property of normality among finite index subgroups of~$\Free[n]$.

\subsubsection{\bf Word-based models}

The mentioned claim on generic hyperbolicity (hinted by \citeauthor{gromov_hyperbolic_1987} in his seminal monograph \parencite{gromov_hyperbolic_1987}, and formally  proved by \citeauthor{olshanskii_almost_1992} and
\citeauthor{champetier_phd_1991}) also fits in the stratified scheme, in this case based on word lengths.

\begin{thm}[\citenr{gromov_hyperbolic_1987}; \citenr{olshanskii_almost_1992}; \citeauthor{champetier_phd_1991} \cite{champetier_phd_1991,champetier_proprietes_1995}]  \label{thm: Gromov}
Generically,  a  finitely  presented  group $\pres{A}{R}$ is  non-trivial  and hyperbolic.
\end{thm}

In this model, the universe is that of (cyclically reduced) finite presentations %
with a fixed number of generators and relators (\ie the set of \mbox{$m$-tuples} of cyclically reduced words in $\Free[n]$). Concretely, the claim is that for every fixed $n\geq 2$ and $m\geq 1$, ${\lim_{l \to \infty} N_h/N = 1} $, where~$N$ is the number of presentations $\pres{a_1,\ldots,a_n}{r_1,\ldots,r_m}$ such that the $r_j$'s are (cyclically\footnote{It is not difficult to see that the result can be stated equivalently for presentations with reduced or cyclically reduced relators (see~\cite{olshanskii_almost_1992}).}) reduced with length $|r_j|=l_j$, $N_h$ is the number of them being hyperbolic, and $l=\min (l_1, \ldots, l_m)$.

Ol'shanskii's proof of \Cref{thm: Gromov} is strongly geometrical, it starts from the interpretation of hyperbolicity  in terms of Van Kampen diagrams in order to perform a graph-theoretical analysis to test small cancellation (and hence hyperbolicity) in the conditions of the statement.

A related (but different) model of genericity, now called the \defin{Arzhantseva-Ol’shanskii model}, was introduced shortly after in~\parencite{arzhantseva_generality_1996,arzhantseva_class_1996}, and frequently used since then. It turns out that in many of these works Stallings automata play, again, an essential role (see \eg \parencite{arzhantseva_class_1996,arzhantseva_on_groups_1997,arzhantseva_phd_1998,arzhantseva_generic_1998,arzhantseva_property_2000,kapovich_genericity_2005,kapovich_random_2009}). A remarkable result in this direction is stated below.

\begin{thm}[\citenr{arzhantseva_class_1996}] \label{thm: generically free}
Exponentially generically, a $m$-tuple of reduced s in $\Free[n]$ generates a subgroup of rank $m$ (\ie is freely independent).
\end{thm}

In the Arzhantseva--Olshanskii model, the universe is again that of $m$-tuples of cyclically reduced words in $\Free[n]$ (\ie  $\Omega = {\Cred_n}^{\,m}$) but the stratification is in terms of \mbox{$l$-balls}, that is $\Omega_l = \set{w\in \Cred_n \st {|w|\leq l}}^m$. Hence, \Cref{thm: generically free} means that, for any fixed $m\geq 1$, and among the set of \mbox{$m$-tuples} of cyclically reduced words of length at most $l$, the proportion of those being freely independent tends to 1, exponentially fast, when~${l\to \infty}$. 

Another interesting consequence in this model is an alternative viewpoint for the genericity of hyperbolic groups, which follows easily from~\cite[Lemma~3]{arzhantseva_class_1996}.

\begin{thm}[\citenr{arzhantseva_class_1996}]
In the Arzhantseva--Olshanskii model, exponentially generically, a finitely presented group is non-trivial and hyperbolic. 
\end{thm}

A slight variation of the Arzhantseva--Olshanskii model consists in considering reduced instead of cyclically reduced words; that is, taking ${\Omega = {\Red_n}^{\hspace{-5pt}m}}$ and $\Omega_l = \set{w\in \Red_{n} \st  |w|\leq l}^{m}$. This is sometimes called the \defin{word-based model} of random groups, and was used by \citeauthor{jitsukawa_malnormal_2002} to
prove the result below, again using Stallings automata. 

\begin{prop}[\citenr{jitsukawa_malnormal_2002}]\label{prop: jitsukawa}
Malnormal (and hence pure) finitely generated subgroups of $\Free[n]$ are exponentially generic.
\end{prop}

Again, the precise meaning of this statement is that, for a fixed $n$, and among the set of \mbox{$m$-tuples} of reduced elements of length at most $l$ in $\Free[n]$, the proportion of those which generate a malnormal (resp., pure) subgroup tends to 1, exponentially fast, when $l \to \infty$.

We would like to point out, however, some nuances on word-based models: (1) the size of the tuples $m$ is fixed, so the statements are not really about \emph{all} finitely generated subgroups of $\Free[n]$, but about those which are \emph{$m$-generated}, \ie those of rank up to $m$; (2) we are really counting $m$-tuples and so neglecting the fact that different $m$-tuples could generate the same subgroup;
(3) different $m$-tuples $R$, $R'$ generating different subgroups, $\gen{R}\neq \gen{R'}$, could generate the same normal closure $\normalcl{R}=\normalcl{R'}\leq \Free[r]$ and so, still present the same group $\pres{A}{R}=\pres{A}{R'}$; and finally (4) even different $k$-tuples $R$, $R'$ generating different normal subgroups, $\normalcl{R}\neq \normalcl{R'}$, could still present isomorphic groups, $\pres{A}{R}\isom \pres{A}{R'}$.
 
One can try to overcome some of the previous objections arguing that, when ${l\to \infty}$, the corresponding redundancy happens with probability tending to $0$. As far as we know, there are no published results formalizing this idea.

\subsubsection{\bf Graph-based model}

An appealing attempt to partially overcome the previous objections was done by \citeauthor{bassino_random_2008} in \parencite{bassino_random_2008} (see also~\parencite{bassino_statistical_2013,bassino_generic_2016,bassino_genericity_2016}), where the authors used for the first time Stallings automata to stratify the family of finitely generated subgroups of the free group, and study them asymptotically.

In this model the stratification uses the number of vertices in the Stallings automaton, called the size of the corresponding subgroup.

\begin{defn}
The \emph{$\Alfi$-size} (\defin[subgroup size]{size} for short) of a given finitely generated subgroup $H\leq \Free[A]$ is the number of vertices in its Stallings automaton, denoted by $|H|=\card{\Verts\stallings(H,A)}$.
\end{defn}

Of course, for any $k\geqslant 1$, $\Free[n]$ has finitely many subgroups of size $k$. Counting the proportion $\alpha(k)$ of them which satisfies a certain property, and studying the asymptotic behaviour of this function when $k\to \infty$, we obtain some precise measure of how frequent is the property under study, among \emph{all} finitely generated subgroups of $\Free[n]$.

The probabilistic model obtained in this way is known as the \defin{graph-based model}, and was introduced in~\parencite{bassino_random_2008} (in contrast with the \emph{word-based model} mentioned above). It has the disadvantage of being combinatorially more complicated than just counting reduced words, but the advantage that we are directly counting the objects we are interested in, namely subgroups, and not tuples generating them (with possible repetitions). The precise approach in~\parencite{bassino_random_2008} is as follows: Fix an alphabet $A=\{a_1, \ldots ,a_n\}$ and a set of $k$ elements, say $\Verts=\{1, \ldots ,k\}$. Giving a Stallings automaton over the set of vertices $\Verts$ (and with basepoint $\bp=1$) is the same as giving an $A$-tuple of partial injections on the set $V$, $(\tau_a)_{a\in A}$, where $\tau_a\in \Verts^{\subseteq \Verts}$ for $a\in A$ (see~\Cref{ssec: involutive automata}). Such an $A$-tuple uniquely defines a deterministic automaton with vertex set $\Verts$ (and basepoint $\bp=1$), which happens to be a Stallings automaton if and only if it is additionally connected and core. A couple of technical lemmas from~\parencite{bassino_random_2008} guarantee that these last two properties happen with probability tending to 1, as $k\to \infty$. So, the problem of randomly generating a Stallings automaton reduces (via an efficient rejection algorithm, see~\parencite{bassino_random_2008}) to the problem of efficiently generating random partial injections on the set~$\Verts$. Based on this nice idea, \parencite{bassino_random_2008, bassino_statistical_2013} contain the following results. 

\begin{prop}[\citenr{bassino_random_2008}]
In the graph-based model, the expected rank of a randomly chosen size $k$ subgroup of $\Free[n]$ is asymptotically equivalent~to $(n-1)k-n\sqrt{k}+1$.
\end{prop}

\begin{prop}[\citenr{bassino_statistical_2013}]
In the graph-based model, pure (and hence malnormal) finitely generated subgroups of $\Free[n]$ are negligible; more precisely, the probability that a random subgroup of $\Free[n]$ of size $k$ is pure (\resp malnormal) is ${\mathcal O}(k^{-n/2})$.
\end{prop}

Note the contrast with \Cref{prop: jitsukawa}: malnormal and pure subgroups are generic in the word-based distribution, whereas they are negligible in the graph-based distribution. This is not a contradiction, but just a reflex of the fact that these two distributions are quite different: they provide two different points of view to observe the (infinite) set of finitely generated subgroups of $\Free[n]$, two different ways to tend to infinity within this set. 

Finally, using the graph-based distribution to analyze finitely presented groups, one also gets the surprisingly contrasting result below. 

\begin{thm}[\citenr{bassino_statistical_2013}]
In the graph-based model, generically, a finitely presented group $\pres{A}{\Ati}$ is trivial (\ie generically, the normal closure of a randomly chosen subgroup of $\Free[n]$, is $\Free[n]$ itself).
\end{thm}

Here, the idea is the following. Partial injections decompose as a product of cycles and oriented paths (and the decomposition is all in cycles if and only if the partial injection is a permutation). It can be proved that the probability that the decomposition of a partial injection on a set of $k$ elements contains at least one cycle tends to 1 when $k\to \infty$. This translates into saying that, in a random presentation $\pres{A}{\Ati}$, any letter represents a torsion element. Strengthening this argument, it can also be proved that, in fact, even the probability that the decomposition of a partial injection on a set of $k$ elements contains at least two cycles with coprime lengths also tends to 1 when $k\to \infty$. And this translates into saying that, in a random presentation $\pres{A}{\Ati}$, \emph{any} letter represents the trivial element, \ie $\pres{A}{\Ati}=1$.

We note that the above concerns (3) and (4) (about different $\Ati$, $\Ati'$ generating the same normal closure, or different normal closures but isomorphic quotients) are equally alive here. But part of concern (1), and concern (2) disappear: we are counting $n$-generated groups, with no restriction on the (finite) number of relations, and \emph{genuine} subgroups, not just tuples of words generating them (with possible repetitions).

\vskip 0.4 true cm
\begin{center}{\textbf{Acknowledgments}}
\end{center}
We are grateful to Goulnara Arzhantseva for some historical remarks related to~\Cref{sec: asymptotic}. The authors acknowledge partial support from the Spanish Agencia Estatal de Investigación, through grant \mbox{MTM2017-82740-P} (AEI/FEDER, UE), and also from the Barcelona Graduate School of Mathematics. The first named author was partially supported by MINECO grant PID2019-107444GA-I00 and the Basque Government grant \mbox{IT974-16}.

%\cleardoublepage
%
%
\renewcommand*{\bibfont}{\footnotesize}
\printbibliography

\clearpage
\printindex

%
%

% \newpage
\bigskip
\bigskip

\noindent{\footnotesize {\bf Jordi Delgado}\; \\ {Department of
Mathematics}, {University of the Basque Country, Faculty of Science and Technology, Barrio Sarriena, s/n} {48940 Leioa, Spain}\\
{\tt Email: jdelgado@crm.cat}\\

\noindent{\footnotesize {\bf Enric Ventura}\; \\ {Departament de Matemàtiques}, {Universitat Politècnica de Catalunya and Institut de Matemàtiques de la UPC-BarcelonaTech,} {Barcelona, Catalunya}\\
{\tt Email: enric.ventura@upc.edu}\\
\end{document}